\documentclass[sn-mathphys,Numbered,iicol]{sn-jnl}

\usepackage{graphicx}%
\usepackage{multirow}%
\usepackage{amsmath,amssymb,amsfonts}%
\usepackage{amsthm}%
\usepackage{mathrsfs}%
\usepackage[title]{appendix}%
\usepackage[table]{xcolor}%
\usepackage{textcomp}%
\usepackage{manyfoot}%
\usepackage{booktabs}%
\usepackage{algorithm}%
\usepackage{algorithmicx}%
\usepackage{algpseudocode}%
\usepackage{listings}%
\usepackage{bm}%
\usepackage{diagbox}%
\usepackage[subrefformat=parens,labelformat=parens]{subcaption}%

\definecolor{MatlabPurple}{rgb}{0.49,0.19,0.63}
\definecolor{MatlabPurple2}{rgb}{0.6,0.2,0.6}
\definecolor{MatlabCyan}{rgb}{0,1,1}
\definecolor{MatlabDarkGreen}{rgb}{0,0.5,0}
\definecolor{MatlabGreen}{rgb}{0,1,0}
\definecolor{MatlabCol1}{rgb}{0,0.2353,0.6824}
\definecolor{MatlabCol2}{rgb}{0.7804,0.0039,0.0549}
\definecolor{MatlabCol3}{rgb}{0.2353,0.6078,0.2745}
\definecolor{MatlabCol4}{rgb}{0.6000,0.3020,0.6000}
\definecolor{MatlabCol5}{rgb}{1.0000,0.6078,0}
\definecolor{MatlabCol6}{rgb}{0.0039,0.5490,0.7569}
\definecolor{Gray1}{rgb}{0.6,0.6,0.6}
\definecolor{Gray2}{rgb}{0.5,0.5,0.5}
\definecolor{schRed}{rgb}{0.7,0,0}
\definecolor{schBlue}{rgb}{0,0,0.7}

\theoremstyle{thmstyleone}%
\theoremstyle{thmstyletwo}%

\theoremstyle{thmstylethree}%

\begin{document}

\title[Experimental continuation in nonlinear dynamics: recent advances and future challenges]{Experimental continuation in nonlinear dynamics: recent advances and future challenges}


\author*{\fnm{Ghislain} \sur{Raze}}\email{g.raze@uliege.be}

\author{\fnm{Gaëtan} \sur{Abeloos}}

\author{\fnm{Gaëtan} \sur{Kerschen}}

\affil*{\orgdiv{Aerospace and Mechanical Engineering Department}, \orgname{University of Liège}, \orgaddress{\street{Allée de la Découverte, 9}, \city{Liège}, \postcode{4000}, \country{Belgium}}}



\abstract{
Experimental continuation encompasses a set of methods that combine control and continuation to obtain the full bifurcation diagram of a nonlinear system experimentally, including responses that would be unstable in the system without feedback control. Such control-based methods thus allow the experimenter to directly and exhaustively explore the dynamics of the system without the need for a good mathematical model. The objective of this paper is twofold, namely (i) to review and present the state-of-the-art methods in a unified manner and (ii) to introduce a novel experimental derivative-free arclength continuation procedure, termed arclength control-based continuation. These methods are also demonstrated and compared on an electronic Duffing oscillator and a clamped thin plate featuring geometrical nonlinearity. Finally, the current state of the art is reflected upon, and the challenges lying ahead for this growing field are discussed.
}




\maketitle

\section{Introduction}

The theory of nonlinear dynamical systems was developed based on Poincaré's seminal work. It witnessed extraordinary advances during the 20th century thanks to Arnold,  Holmes, Lorenz, Ruelle and Smale to name a few. A couple of decades ago, the theory spread across engineering thanks to a series of reference monographs thoroughly discussing the different phenomena and attractors nonlinear systems can exhibit. Nayfeh and Mook~\cite{Nayfeh1985} applied perturbation methods to study nonlinear phenomena in single- and multi-degree-of-freedom systems. Guckenheimer and Holmes~\cite{Guckenheimer1984} adopted a different, geometrical viewpoint, appearing as an ideal companion to the perturbation approach. Kuznetsov~\cite{Kuznetsov2004} published a complete treatise on bifurcations.

Since the 1970s, impressive progress was also realized in computational nonlinear dynamics with the development of  nonlinear finite element methods~\cite{Wriggers} and of numerical continuation (NC)~\cite{Keller1977,allgower2003introduction,Dhooge2003,govaerts2000numerical}. Numerical continuation amounts to calculating a bifurcation diagram, a key concept in nonlinear dynamics~\cite{Nayfeh1995} representing the different branches of responses of a nonlinear system together with their stability and bifurcations as one or more parameters are varied. Today, the calculation of bifurcation diagrams of complex numerical models is within reach~\cite{Knowles,Detroux2015,Dankowicz2011,Freour}. Noteworthy are also the advances in the area of nonlinear model reduction~\cite{Haller2016,Vizzaccaro2021,Jain2022,Touze2021}. We also note that there now exists the possibility to perform direct, data-driven reduction onto invariant manifolds~\cite{Cenedese2022}.

Although not all challenges have been overcome yet, the theoretical understanding of nonlinear dynamical phenomena and their prediction using numerical models have reached a high level of maturity. Conversely, experimental nonlinear dynamics \cite{Virgin} received comparatively far less attention and is the main thrust of this paper. In this context, nonlinear system identification (NSI) which focuses on the development of mathematical models directly from experimental data is the most popular approach to date in the literature. There exists a wide variety of techniques see, e.g.,~\cite{Kerschen2006,Noel2017}. Once a good model of the system has been obtained thanks to NSI, numerical continuation techniques can then be exploited to calculate a bifurcation diagram, as carried out, e.g., in~\cite{Detroux2015,Anastasio2023}.  

Despite the recent advances in the field, NSI remains a particularly difficult exercise and imposes a number of challenges and constraints to the experimenter. An elegant, model-less alternative to NSI is experimental continuation
which leverages feedback control for the identification of a bifurcation diagram. This novel experimental strategy represents a true paradigm shift in nonlinear dynamics. There are two key differences with numerical continuation, namely (i) bifurcation diagrams have to be obtained in a model-less manner and (ii) feedback control must be applied to stabilize unstable periodic orbits which otherwise cannot be measured experimentally. The first experimental continuation method, control-based continuation (CBC), pioneered by Sieber et al~\cite{Sieber2008,Sieber2008b}, can be seen as an experimental equivalent to pseudo-arclength continuation. The fact that CBC requires to calculate derivatives using finite differences from experimental data, a computationally-intensive and highly sensitive-to-noise process, motivated researchers to develop methods which do not rely on experimental derivatives. The three most effective derivative-free experimental continuation methods to date, namely the simplified CBC (SCBC) method~\cite{Barton2013}, phase-locked loops (PLLs)~\cite{Peter2017} and response-controlled testing (RCT)~\cite{Karaagacl2020}, resort to natural parameter continuation, which simply increases the bifurcation parameter using the previous solution as an initial guess for the next solution. Thanks to these advances, the frequency responses and backbone curves of nonlinear systems around their primary resonances can now be identified in a robust and effective manner. 

The objective of this feature article is to offer the first extensive review of the developments in the EC area since its invention, to present the three derivative-free methods in an unified manner and to introduce a novel experimental continuation method, termed arclength control-based continuation (ACBC)~\cite{Abeloos2022Thesis}, the first derivative-free arclength continuation scheme.

The paper is organized as follows. Bifurcation diagrams and the means to compute them are introduced in Section~\ref{sec:BifurcationDiagrams}. Section~\ref{sec:Feedback} provides a brief introduction to basic concepts in feedback control. Section~\ref{sec:ECPioneers} presents the core ideas of EC and details its first forms. Section~\ref{sec:ENPC} then reviews the most recent methods and presents them under a unified framework. The novel ACBC approach is presented in Section~\ref{sec:ACBC}. All the methods are compared experimentally with an electronic Duffing oscillator and a clamped plate featuring geometrical nonlinearity in Sections~\ref{sec:ElectronicDuffing} and~\ref{sec:Plate}, respectively. Finally, the conclusions of the present study and the challenges lying ahead of us are summarized in Section~\ref{sec:Conclusion}.

\section{Bifurcation diagrams: a central tool in nonlinear dynamics}
\label{sec:BifurcationDiagrams}

Nonlinear dynamical systems usually behave in a complex manner, requiring adequate tools to study their behavior. Bifurcation diagrams are ubiquitously used in nonlinear dynamics for that purpose. These diagrams typically present how one quantity characterizing the response of the system under study is influenced by a parameter, often called a \textit{bifurcation parameter}. In this work, we focus on the periodic solutions of nonlinear systems and their stability (but fixed points can be seen as a special subcase).

In this section, the main features of analytical and numerical approaches to draw out bifurcation diagrams are summarized, and used to illustrate the challenges faced in experimental nonlinear dynamics.

\subsection{Periodic solutions}

 We consider a generic, non-autonomous system governed by the set of first-order ordinary differential equations (ODEs)
\begin{equation}
     \dot{\mathbf{x}}(t) = \mathbf{f}(\mathbf{x}(t),t),
     \label{eq:ODE}
\end{equation}
where $\mathbf{x}$ is a vector gathering the states of the system under consideration,  $t$ is the time, and $\mathbf{f}$ is a state evolution law. We note that this framework encompasses the case of a set of second-order ODEs generally used in vibration problems
\begin{equation}
    \mathbf{M} \ddot{\mathbf{q}}(t) + \mathbf{C} \dot{\mathbf{q}}(t) + \mathbf{K} \mathbf{q}(t) + \mathbf{f}_{\mathrm{nl}}(\mathbf{q}(t),\dot{\mathbf{q}}(t)) = \mathbf{f}_{\mathrm{ext}}(t),
\end{equation}
where $\mathbf{M}$, $\mathbf{C}$ and $\mathbf{K}$ are linear structural mass, damping and stiffness matrices, respectively, $\mathbf{f}_{\rm{nl}}$ is the vector of nonlinear forces and $\mathbf{f}_{\mathrm{ext}}$ is the vector of external forces. Indeed, this set of second-order ODEs can be recast into a set of first-order ODEs (provided that the mass matrix is invertible) using
\begin{equation}
    \begin{array}{rll}
    \displaystyle
    \dot{\mathbf{x}} & = \displaystyle\begin{bmatrix}
        \dot{\mathbf{q}} \\ \ddot{\mathbf{q}}
    \end{bmatrix} = & \displaystyle\begin{bmatrix}
        \mathbf{0} & \mathbf{I} \\ -\mathbf{M}^{-1}\mathbf{K} & -\mathbf{M}^{-1}\mathbf{C}
    \end{bmatrix} \begin{bmatrix}
        \mathbf{q} \\ \dot{\mathbf{q}}
    \end{bmatrix} \\ & &\displaystyle + \begin{bmatrix}
        \mathbf{0} \\ \mathbf{M}^{-1} \left(\mathbf{f}_{\mathrm{ext}}-\mathbf{f}_{\mathrm{nl}}(\mathbf{q},\dot{\mathbf{q}})\right)
    \end{bmatrix} \\& \displaystyle = \mathbf{f}(\mathbf{x},t). &
    \end{array}
\end{equation}

When the right-hand side of Equation~\eqref{eq:ODE} is periodic in time, periodic solutions $\mathbf{x}$ generally exist, satisfying
\begin{equation}
    \mathbf{x}(t+T) = \mathbf{x}(t), \qquad \forall t \in \mathbb{R}
\end{equation}
where $T$ is a period of the solution, being related to its angular frequency $\omega$ by $\omega T = 2\pi$.

\subsection{Analytical solutions}

    Equation~\eqref{eq:ODE} with given initial conditions has a unique solution provided that $\mathbf{f}$ satisfies mild smoothness conditions given by Picard's existence theorem. In some cases, analytical insight into these solutions can be gained.

    \subsubsection{Exact solutions}

        Due to the very nature of nonlinear systems, exact analytical solutions seldom exist. A few special or nongeneric cases can be treated, such as the unforced, undamped pendulum or Duffing oscillator with the use of Jacobi elliptic functions~\cite{Kovacic2011}.
        
        For illustration, we consider the free response $q$ of an undamped Duffing oscillator, a typical academic example qualitatively representative of geometrically nonlinear structures, which is governed by the following ODE:
        \begin{equation}
            m\ddot{q}(t) + kq(t) + k_3q^3(t) = 0,
            \label{eq:DuffingUnforcedODE}
        \end{equation}
        where $m$, $k$ and $k_3$ are the mass, linear stiffness and cubic stiffness coefficients, respectively, of the oscillator and an overdot denotes a derivation with respect to time $t$. Dividing Equation~\eqref{eq:DuffingUnforcedODE} by $m$,
        \begin{equation}
            \ddot{q}(t) + \omega_0^2 q(t) + \alpha_3 q^3(t) = 0
            \label{eq:DuffingUnforcedODENormalized}
        \end{equation}
        with $\omega_0^2 = k/m$ and $\alpha_3 = k_3/m$. The exact solution of Equation~\eqref{eq:DuffingUnforcedODENormalized} is given by
        \begin{equation}
            q(t) = A \text{sn} (\Omega t + b |\kappa),
        \end{equation}
        where $A$ is the amplitude of the free oscillations, $\text{sn}$ is the elliptic sine function, $b$ is a phase-like variable to be adjusted to the initial conditions (together with $A$), $\Omega$ is a frequency-like variable and $\kappa$ is the elliptic parameter~\cite{Kovacic2011}. These two latter parameters are related to $A$ by
        \begin{equation}
            \Omega = \sqrt{\omega^2_0 + \dfrac{\alpha_3 A^2}{2}},\qquad \kappa = -\dfrac{\alpha_3 A^2}{2\Omega^2}.
        \end{equation}


    \subsubsection{Perturbation methods}

        Perturbation methods constitute a means to obtain approximate analytical solutions to nonlinear problems. They rely on a small parameter, usually denoted $\epsilon$. In most cases, the external forcing, response amplitude, nonlinear coefficients and/or damping coefficients are assumed to scale with $\epsilon$. When $\epsilon =0$, a (generally known) solution forms the basis around which a power series expansion in terms of $\epsilon$ is performed.

        Perturbation methods are a powerful tool to study nonlinear systems analytically. The most famous ones are the method of multiple scales~\cite{Nayfeh1985} and the averaging method~\cite{Guckenheimer1984}.
        These methods can provide insight into the system's behavior, but are inherently restricted by the limitations associated with power series expansions, and generally break down when one or several assumed small parameters become too large. 

\subsection{Numerical procedures}

If analytical solutions are unavailable, intractable or inaccurate, periodic solutions for a fixed $\omega$ can be found by numerically integrating Equation~\eqref{eq:ODE} until the transients die out. This approach is generally inefficient, and more sophisticated approaches can be used, such as the shooting~\cite{seydel2009practical}, harmonic balance~\cite{Detroux2015,Krack2019} and orthogonal collocation~\cite{dankowicz2013recipes} methods. These methods cast the ODEs into a set of algebraic equations that can be solved by well-known numerical routines such as the Newton-Raphson method. This usually yields accurate results. 

To gain a broader understanding of the system's behavior, a bifurcation diagram can be built by varying a bifurcation parameter. To find the family of periodic solutions that form a complete bifurcation diagram, numerical continuation can be used, as will be presented next.

\subsection{Numerical continuation}
\label{sec:numericalContinuation}

    Numerical continuation (also known as branch tracing or path-following method) denotes a class of numerical methods used to solve a system of algebraic equations with one (or more) free parameter. The solutions of these equations generally form a set whose dimension is equal to the number of bifurcation parameters. Since the most common case considers a single parameter, this set is generically a curve. We treat this case without loss of generality, because the case with multiple bifurcation parameters can be treated, e.g., by sampling the values of all but one parameter, and using a one-parameter continuation for all samples, thereby effectively "slicing" the solution set.
    
\begin{figure}[!ht]
        \centering
        \includegraphics[scale=1]{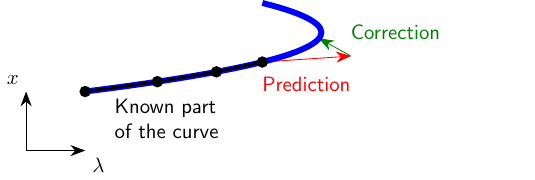}
        \caption{Schematics of a general continuation method.}
        \label{fig:continuationSchematics}
    \end{figure}

    Numerical continuation is an incremental procedure. The main idea behind it is to build a branch of solutions point by point (or piece by piece), using a constructive approach where the information at the previous point (or piece) of the curve is known to be well estimated and is used to find a suitable initial guess for the next one. This defines the prediction step. This initial guess is then refined to ensure that it lies on the sought curve within a prescribed tolerance, which defines the correction step. Figure~\ref{fig:continuationSchematics} schematically depicts the main ideas behind this predictor-corrector scheme, where $x$ is one of the unknowns, and $\lambda$ is a bifurcation parameter.

    Numerical continuation is a mature tool that is now treated in a variety of reference books and exploited in a large number of bifurcation analysis software, such as AUTO~\cite{DOEDEL1991}, MATCONT~\cite{Dhooge2003}, COCO~\cite{Dankowicz2011} and MANLAB~\cite{Cochelin2009}. In this work, we will focus on the two most popular variants of numerical continuation, namely natural parameter and arclength continuations. There also are numerous books discussing this topic, such as~\cite{govaerts2000numerical,allgower2003introduction,krauskopf2007numerical,seydel2009practical,dankowicz2013recipes}.

    \subsubsection{Natural parameter continuation}
    \label{ssec:NPC}

\begin{figure}[!ht]
            \centering
        	\includegraphics[scale=1]{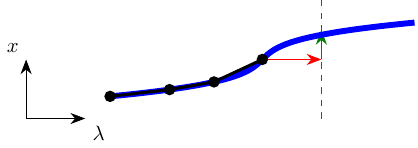}
            \caption{Schematics of the natural parameter continuation method.}
            \label{fig:seqContinuationSchematics}
        \end{figure}

        Natural parameter continuation (also called sequential continuation) is the simplest form of continuation. It uses the bifurcation parameter of the problem to parametrize the solution curve. The implicit function theorem then gives conditions under which this choice is guaranteed to be sensible, i.e., when the solution curve is indeed (locally) parametrizable by the bifurcation parameter. When the solution is known for a given value of the parameter, the prediction step simply consists in increasing the parameter value, using the previous solution as initial guess for the next point. The correction phase generally consists in using an iterative solver such as the Newton-Raphson method. Figure~\ref{fig:seqContinuationSchematics} represents a general schematics of a natural parameter continuation method.

\begin{figure}[!ht]
                \centering
                \begin{subfigure}{.49\textwidth}
                    \centering
                    \includegraphics[scale=1]{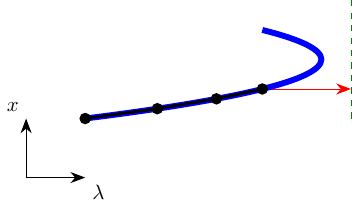}
                    \caption{}
                    \label{sfig:seqContinuationIssues1}
                \end{subfigure}
                \begin{subfigure}{.49\textwidth}
                    \centering
        			\includegraphics[scale=1]{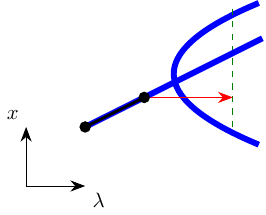}
                    \caption{}
                    \label{sfig:seqContinuationIssues2}
                \end{subfigure}
                \caption{Typical scenarii where the natural parameter continuation fails or has an unpredictable outcome: folds~\subref{sfig:seqContinuationIssues1} and branch points~\subref{sfig:seqContinuationIssues2}.}
                \label{fig:seqContinuationIssues}
            \end{figure}

        Natural parameter continuation fails at points where the implicit function theorem does not hold when one parametrizes a curve with the continuation parameter. Two notable cases are limit points and branch points, represented in Figures~\ref{sfig:seqContinuationIssues1} and~\ref{sfig:seqContinuationIssues2}, respectively. In the presence of the former, natural parameter continuation fails most of the time because there is no nearby solution at the imposed parameter value, as shown in Figure~\ref{sfig:seqContinuationIssues1}. With branch points, it is generally hard to predict and control which path the method will follow close to the bifurcation; this depends on the step length and on the problem at hand (cf. Figure ~\ref{sfig:seqContinuationIssues2}). To address these issues, more elaborate continuation approaches were developed.

    \subsubsection{Arclength continuation}

        The first ideas behind arclength continuation can be traced back to the seminal work of Haselgrove in the early sixties~\cite{Haselgrove1961}, and were extensively developed and popularized by Keller~\cite{Keller1977}. It was (seemingly independently) popularized nearly twenty years later in the field of static mechanics with the famous works of Riks~\cite{Riks1979} and Crisfield~\cite{CRISFIELD1981}, who respectively proposed what is now known as the pseudo-arclength and arclength continuation procedures. Adaptations to dynamics came ten years later using the harmonic balance~\cite{Lewandowski1992} and shooting~\cite{Sundararajan1997} methods. The underlying idea of these methods is to parametrize the curve with an (approximation of its) arclength.
        
\begin{figure}[!ht]
            \centering
            \begin{subfigure}{.49\textwidth}
                \centering
        		\includegraphics[scale=1]{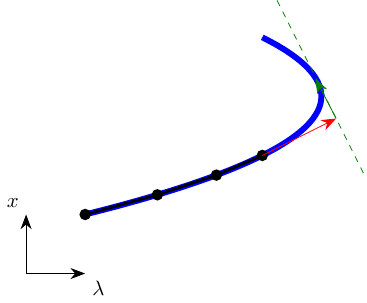}
                \caption{}
                \label{sfig:pseudoArclenghtContinuation}
            \end{subfigure}
            \begin{subfigure}{.49\textwidth}
                \centering
        		\includegraphics[scale=1]{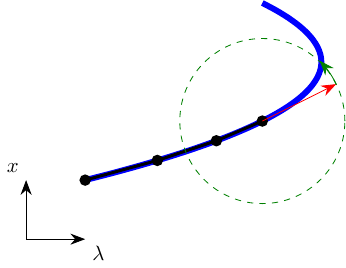}
                \caption{}
                \label{sfig:arclengthContinuation}
            \end{subfigure}
            \caption{Schematics of the pseudo-arclength~\subref{sfig:pseudoArclenghtContinuation} and arclength~\subref{sfig:arclengthContinuation} methods.}
            \label{fig:arclengthContinuation}
        \end{figure}

        These methods both use a more refined prediction step with either a tangent or secant predictor. The correction is then performed either on a hyperplane whose normal is the prediction vector (Figure~\ref{sfig:pseudoArclenghtContinuation}), or on a hypersphere centered around the previous point and whose radius is equal to the prediction step length (Figure~\ref{sfig:arclengthContinuation}). Mathematically, $\lambda$ is let free, and the correction locus adds another constraint to the problem to make it well-posed. Clearly, these methods offer better robustness than natural parameter continuation when the curve folds. In addition, the tangent predictor generally makes the method more likely to stay on the same branch when a branching behavior is encountered. These are some of the reasons for their current popularity and ubiquitous presence in numerical continuation packages.

\subsection{Stability}
\label{ssec:stability}

Equipped with the possibility to draw bifurcation diagrams, we now have to wonder whether the obtained family of solutions will be observable in the nonlinear system under study. By contrast with linear time-invariant (LTI) systems for which solutions are either all stable or all unstable, nonlinear systems can possess solutions of different stability. In addition, while unstable solutions of LTI systems grow unbounded with time, trajectories repelled by a nearby unstable nonlinear solution have the additional possibility to converge to other stable solutions.

In an experimental context, the unstable solutions of a system are impossible to observe in steady-state conditions without additional control. As for numerical aspects, unstable solutions can be computed, but an important question is to know whether this solution is stable or not. To answer this question, we look at the evolution of infinitesimally small perturbations around the equilibrium, which gives us results about the asymptotic stability of a solution. For more elaborate discussions on stability, we refer to, e.g.,~\cite{Guckenheimer1984,slotine1991applied}.

We now assume that we know a periodic solution of Equation~\eqref{eq:ODE} (obtained analytically or numerically) and denote it by $\mathbf{x}_*$. A small perturbation $\mathbf{s}$ is added to the equilibrium $\mathbf{x}_*$, i.e.,
\begin{equation}
    \mathbf{x}(t) = \mathbf{x}_*(t) + \mathbf{s}(t), \qquad |\mathbf{s}(t)| << |\mathbf{x}_*(t)|
\end{equation}
Inserting this perturbed solution into Equation~\eqref{eq:ODE} and expanding the result in Taylor series around $\mathbf{s}(t) = 0$, we obtain
\begin{equation}
    \begin{array}{rl}
     \displaystyle\dot{\mathbf{x}}_*(t) + \dot{\mathbf{s}}(t) = & \displaystyle\mathbf{f}(\mathbf{x}_*,t) + \left.\dfrac{\partial}{\partial\mathbf{x}}\mathbf{f}(\mathbf{x},t)\right|_{\mathbf{x}=\mathbf{x}_*}\mathbf{s}(t) \\ & \displaystyle+ O(|\mathbf{s}(t)|^2).
     \end{array}
\end{equation}
Neglecting terms of orders higher than linear in $|\mathbf{s}|$ and accounting for the fact that $\mathbf{x}_*$ is a solution of Equation~\eqref{eq:ODE}, one gets the first-order variational problem 
\begin{equation}
     \dot{\mathbf{s}}(t) = \left.\dfrac{\partial}{\partial\mathbf{x}}\mathbf{f}(\mathbf{x},t)\right|_{\mathbf{x}=\mathbf{x}_*(t)}\mathbf{s}(t).
\end{equation}
This is a linear, \textit{time-varying} equation featuring the Jacobian matrix $\partial \mathbf{f}/\partial\mathbf{x}$ and governing the evolution of the (small) perturbations. The solution of such an equation cannot be known in closed-form in general, and requires numerical integration. Nevertheless, the time-varying part of this equation is periodic, and this periodicity entails special consequences which are described by the Floquet theory~\cite{Nayfeh1995}. The variational equation can be solved (using numerical time integration) with its principal fundamental matrix $\mathbf{\Phi}(t)$ determined by
\begin{equation}
     \dot{\mathbf{\Phi}}(t) = \left.\dfrac{\partial}{\partial\mathbf{x}}\mathbf{f}(\mathbf{x},t)\right|_{\mathbf{x}=\mathbf{x}_*(t)}\mathbf{\Phi}(t), \quad \mathbf{\Phi}(0) = \mathbf{I}.
     \label{eq:principalFundamentalMatrix}
\end{equation}
Due to the linearity and periodic nature of the first-order variational problem, it can be shown that the growth or decay of \textit{any} small perturbation is related to the eigenvalues of this principal fundamental matrix evaluated at $t=T$ (also known as the monodromy matrix) $\mathbf{\Phi}(T)$, called the Floquet (or characteristic) multipliers. Specifically, if any Floquet multiplier is greater than one in modulus, the solution is asymptotically unstable, and it is stable otherwise.

The monodromy matrix is thus a determining factor for stability. Unfortunately, there generally exists no explicit expression of this matrix from $\partial \mathbf{f}/\partial\mathbf{x}$, but the key point to remember is that \textit{this Jacobian matrix plays a pivotal role for the asymptotic stability of a solution}. 

\subsection{Illustration with a Duffing oscillator}

For illustration, we will consider the damped, forced Duffing oscillator. This example is well-known for the nonlinear dynamicist. However, it can exhibit rather complex dynamics associated with a wide variety of secondary resonances~\cite{Parlitz1985,Marchionne2018}.

The response $q$ of the considered oscillator to an external excitation $f$ is governed by a generalization of Equation~\eqref{eq:DuffingUnforcedODE}:
\begin{equation}
    m\ddot{q}(t) + c\dot{q}(t) + kq(t) + k_2 q^2(t) + k_3q^3(t) = f(t),
    \label{eq:DuffingODE}
\end{equation}
where $c$ is the damping coefficient of the oscillator, and we added a quadratic stiffness force with a coefficient $k_2$ that can be used to represent asymmetries in the system. When $k_2\neq 0$, the oscillator characterized by Equation~\eqref{eq:DuffingODE} is often called the Helmholtz-Duffing oscillator. Upon using the normalizations $\bar{t} = \omega_0 t$ and $\bar{q}=\sqrt{k_3/k}q$, Equation~\eqref{eq:DuffingODE} can be put into the dimensionless form 
\begin{equation}
    \bar{q}''(\bar{t}) + 2\zeta_0 \bar{q}'(\bar{t}) + \bar{q}(\bar{t}) + \beta_2 \bar{q}^2(t) + \bar{q}^3(\bar{t}) = \bar{f}(\bar{t}),
    \label{eq:DuffingODENorm}
\end{equation}
with $\zeta_0=c/(2\sqrt{km})$, $\beta_2 = k_2/\sqrt{kk_3}$, $\bar{f} = \sqrt{k_3/k^3}f$, and a prime denotes derivation with respect to the dimensionless time $\bar{t}$. For the examples shown hereafter, we chose $\zeta_0 = 0.05$ and $\beta_2=0$ (hence we consider a Duffing oscillator), unless stated otherwise.

\begin{figure*}[!ht]
    \centering
    \includegraphics[width=\linewidth]{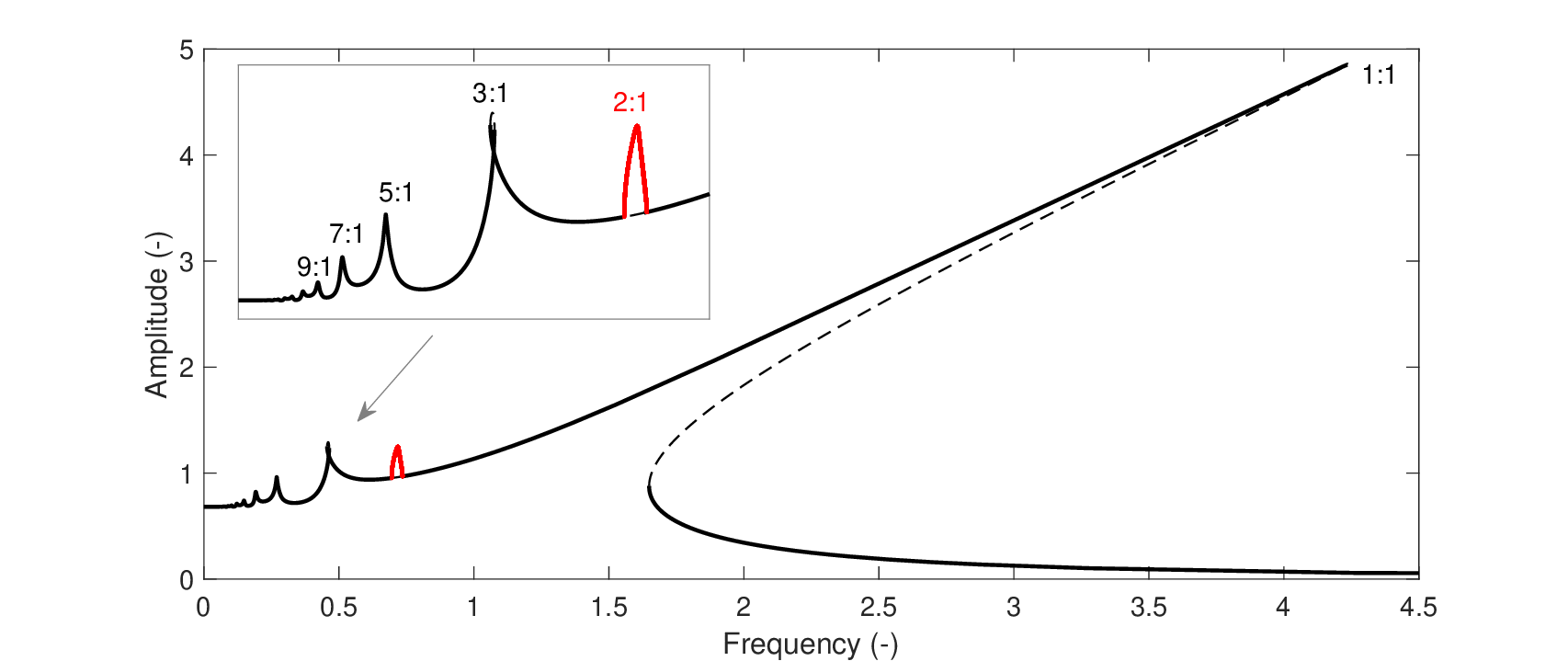}
    \caption{Bifurcation diagram of the Duffing oscillator with a dimensionless forcing amplitude of 1: main branch of the NFR (\rule[.2em]{1em}{.2em}) and bifurcating even superharmonic resonances (\textcolor{red}{\rule[.2em]{1em}{.2em}}). Full and dashed lines represent stable and unstable solutions, respectively.}
    \label{fig:Duffing_Res_Feature_Low}
\end{figure*}

Figure~\ref{fig:Duffing_Res_Feature_Low} features a bifurcation diagram obtained by coupling a harmonic balance procedure with numerical continuation~\cite{Detroux2015,Krack2019}. This diagram represents the response amplitude of the Duffing oscillator under a harmonic forcing of amplitude $\bar{f}(\bar{t})=\sin(\bar{\omega} \bar{t})$ (with a dimensionless frequency such that $\omega=\omega_0\bar{\omega}$), taking $\bar{\omega}$ as bifurcation parameter. This bifurcation diagram is also often called the nonlinear frequency response (NFR). In addition to the well-known primary (1:1) resonance that bends rightward due to the hardening nonlinearity of the oscillator, a series of peaks appear throughout the bifurcation diagram and are associated to non-primary resonances~\cite{stoker1950nonlinear}. Specifically, odd superharmonic ($m$:1, $m\in \mathbb{N}_0$) resonances appear on the main branch of the NFR at frequencies lower than that of the primary resonance, whereas the even superharmonic resonances bifurcate from it through symmetry-breaking bifurcations.

\begin{figure*}[!ht]
    \centering
    \includegraphics[width=\linewidth]{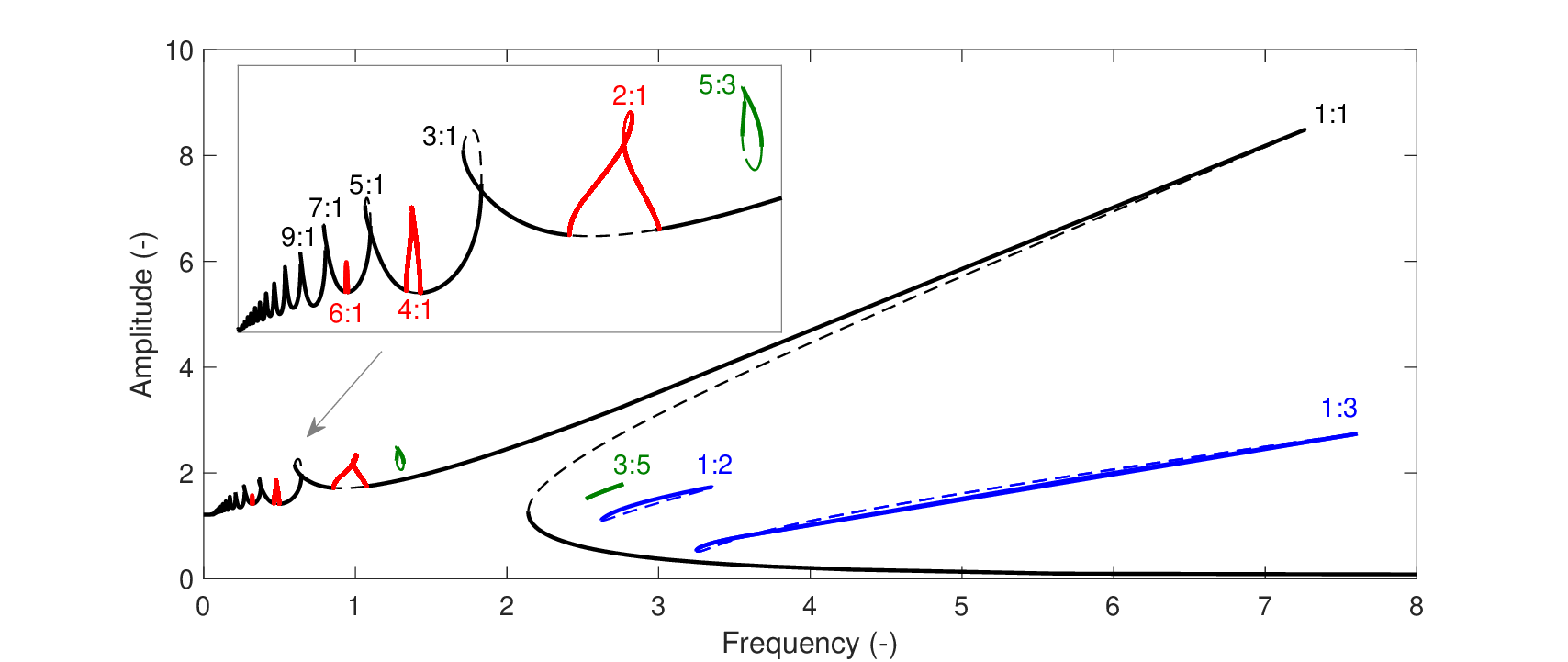}
    \caption{Bifurcation diagram of the Duffing oscillator with a dimensionless forcing amplitude of 3: main branch of the NFR (\rule[.2em]{1em}{.2em}) bifurcating even superharmonic resonances (\textcolor{red}{\rule[.2em]{1em}{.2em}}), isolated subharmonic resonances (\textcolor{blue}{\rule[.2em]{1em}{.2em}}) and isolated ultrasubharmonic resonances (\textcolor{MatlabDarkGreen}{\rule[.2em]{1em}{.2em}}). Full and dashed lines represent stable and unstable solutions, respectively.}
    \label{fig:Duffing_Res_Feature}
\end{figure*}

Figure~\ref{fig:Duffing_Res_Feature} represents the NFR of the Duffing oscillator under a higher-amplitude harmonic forcing, namely $\bar{f}(\bar{t})=3\sin(\bar{\omega} \bar{t})$. In addition to the primary and superharmonic resonances, subharmonic (1:$l$, $l\in\mathbb{N}_0$) and ultrasubharmonic ($m$:$l$, $m,l\in\mathbb{N}_0$, $\gcd(m,l)=1$) resonances appear, and they are isolated from the main NFR. We note that these isolated resonances were not obtained directly from numerical continuation of the main branch; stochastic procedures were used to find a set of initial conditions resulting in a subharmonic or ultrasubharmonic response. Once a suitable initial point was found, numerical continuation was used to trace out the complete isola. Thus, even in this simple oscillator, complex nonlinear features can appear.

\subsection{Consequences for open-loop testing}
\label{ssec:openLoop}


We now illustrate how the complexity featured by nonlinear systems as deceptively simple as the Duffing oscillator impacts the experimental results obtained while testing a system.

During a test, a drive signal is generated by the experimenter and fed to the system. In vibration testing, this drive signal is fed to an actuator (generally through an amplifier) whose purpose is to excite the structure under test, either via a force or via base excitation. In both cases, the drive is seldom proportional to the forcing signal because the amplifier and the actuator have their dynamics~\cite{Pacini2022}, hence we distinguish the drive from the actual forcing. For a harmonic excitation, the two parameters that the user can control are the drive amplitude and frequency. The response of a nonlinear system is generally characterized for multiple values of these parameters. For instance, the drive amplitude may be fixed and the frequency varied. This variation can be continuous or discrete, representing swept-sine (SWS) or stepped-sine (STS) excitations, respectively. The set of points obtained from different frequencies can be interpreted as the NFR. Alternatively, one can fix the frequency and vary the drive amplitude and obtain a so-called S-curve~\cite{Barton2013}.

\begin{figure*}[!ht]
    \centering
    \includegraphics[width=\linewidth]{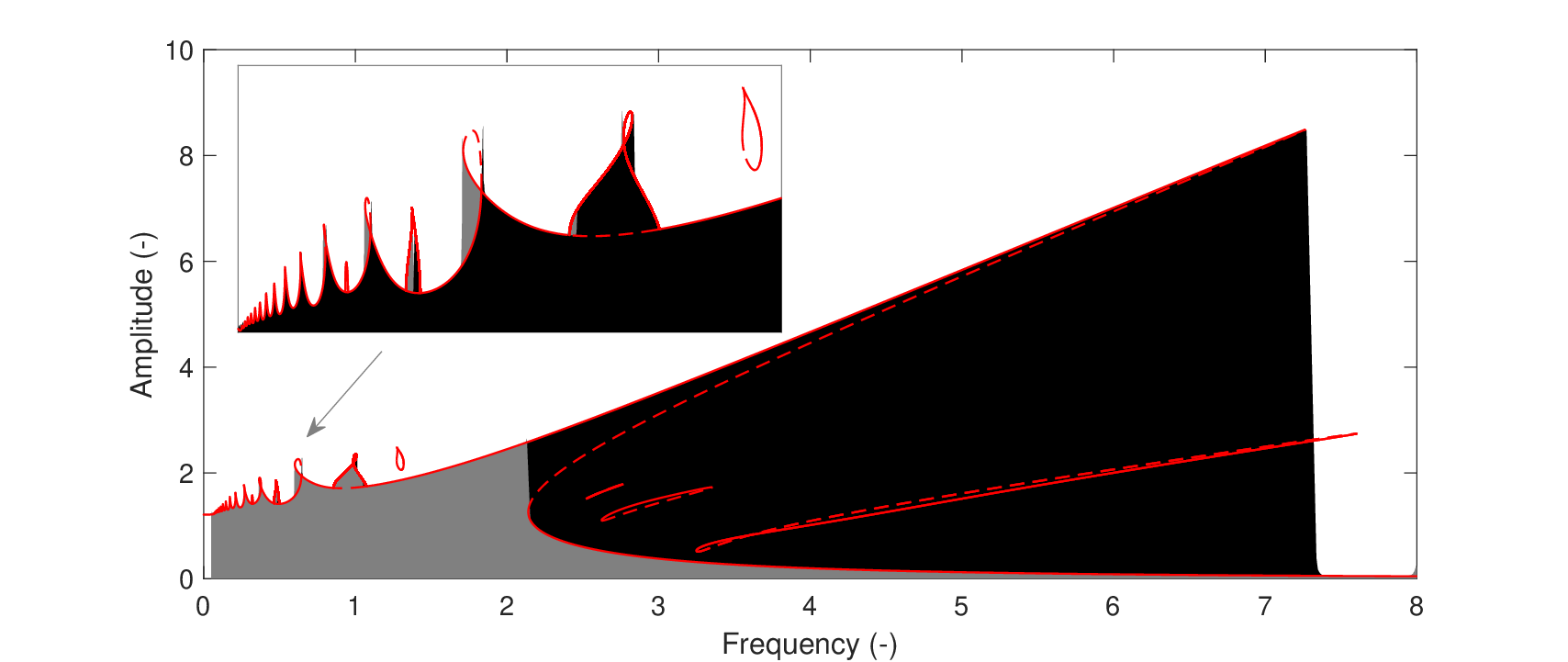}
    \caption{Bifurcation diagram of the Duffing oscillator with a dimensionless forcing amplitude of 3: swept-up (\rule[.2em]{1em}{.2em}) and swept-down (\textcolor{Gray1}{\rule[.2em]{1em}{.2em}}) excitation frequency.
    The results from numerical continuation are given for reference (\textcolor{red}{\rule[.2em]{1em}{.2em}}), where full and dashed lines represent stable and unstable solutions, respectively.}
    \label{fig:DuffingSweep}
\end{figure*}

We apply a SWS with logarithmically varying frequency to the Duffing oscillator, as shown in Figure~\ref{fig:DuffingSweep}. As the frequency increases, a smooth change in amplitude is observed, except in the vicinity of resonances. The most noticeable feature is the occurrence of 
a sudden jump-down at a frequency $\bar{\omega} \approx 7.2$. Conversely, by progressively decreasing the frequency, a jump-up occurs around $\bar{\omega} \approx 2.1$. The jump-up and jump-down phenomena occur where the NFR \textit{folds}, as seen in Figure~\ref{fig:DuffingSweep}. Indeed, SWS can be seen as an experimental equivalent to natural parameter continuation. 

Closing up on the low-frequency region, we can observe most odd superharmonic resonances, whereas the even ones are sometimes partially missed due to the transient nature of the excitation. Sub- and ultrasubharmonic resonances are completely missed by the (simulated) testing procedure because of their isolated nature. 



\subsection{Challenges for experimental nonlinear dynamics}

To summarize this section, we have seen that a combination of theoretical and numerical techniques can be used to obtain bifurcation diagrams that provide a great deal of information about the response of nonlinear systems. Thanks to the maturity of these approaches, and in particular the numerical continuation procedure, stable and unstable parts of a bifurcation diagram can be obtained.

As for experiments, the classical open-loop testing procedures cannot offer the same picture. There are in particular specific challenges associated with nonlinear dynamics:
\begin{enumerate}
    \item \textbf{Unstable solutions} are by nature impossible to follow with open-loop approaches.
    \item \textbf{Jumps} are thus susceptible to occur when a stability change is encountered through bifurcations.
    \item \textbf{Several solutions are missed} during a test, resulting in an incomplete bifurcation diagram. 
    \item \textbf{Isolated solutions} are challenging to reach with classical continuation approaches. In an open-loop setting, they can be obtained in some cases by using a stochastic approach~\cite{Virgin1998}, or leveraging jumps~\cite{Ludeke1951}. 
\end{enumerate}

A novel class of testing methods, herein called control-based methods, has emerged slightly before the beginning of the last decade~\cite{SIEBER2007,Sieber2008}. The salient capabilities of these methods are that, through suitable automatic adjustments of the drive signal, they (\textit{i}) allow the experimenter to prescribe given parameters of the system to extract complete bifurcation diagrams and (\textit{ii}) can stabilize the equilibria of the system under test. These methods are thus paradigm-shifting for nonlinear vibration testing and allow the experimenter to perform experimental continuation.

\section{Feedback control: a brief introduction}
\label{sec:Feedback}

    Before delving into control-based methods, we review a few essential aspects about feedback control, which is a key enabler for experimental continuation.

    \subsection{Linear feedback control}

   Humans have always wanted to control their environment, i.e., to stir the systems that surround them toward desired states of operation. A natural way to perform this is to use feedback control, wherein the output of a system to be controlled (usually called plant) is sensed and used to adapt its input, eventually (and hopefully) reaching a desired output through the action of that input.  When feedback is devoid of human interaction, one speaks of automatic control. We shall refer a system with and without feedback to as closed- and open-loop systems, respectively.

   Feedback is naturally used by most animals, but it can be surprising to realize that this concept was purposefully exploited rather late and formalized even later in human history~\cite{Mayr1970}. Starting from the twentieth century, the study of feedback control has known a great evolution, mostly motivated by technological needs and advances, leading to the ubiquitous presence of controlled systems around us today. For a historical perspective on this subject, we refer, e.g., to~\cite{Fuller1963,Bennett1996}.

\begin{figure*}[!ht]
        \centering
        \includegraphics[scale=.75]{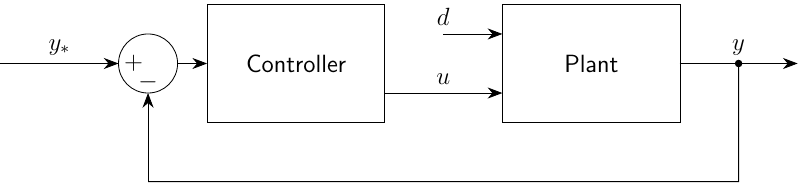}
        \caption{Schematics of a typical control system.}
        \label{fig:controlSchematics}
    \end{figure*}

    Figure~\ref{fig:controlSchematics} presents the typical architecture of a linear control system. The plant, excited by the drive (or control signal) $u$ and unknown disturbances $d$, is monitored through its response $y$ (typically an element of $\mathbf{x}$ or a linear combination of them). This response is compared to a desired reference $y_*$ (in the sequel, a subscript with a star indicates a target quantity) and fed to a controller, that automatically adapts the drive $u$. The controller can be designed to force $y$ to follow the reference $y_*$ as closely as possible, to reject the effect of the external disturbance $d$, or to do both these things at the same time.
    
    In this linear time-invariant setting, the plant response can be expressed as the superposition of the effect of the drive $u$ and the disturbance $d$. We use a Laplace transform formalism to express this as
    \begin{equation}
        Y(p) = H_{yu}(p) U(p) + H_{yd}(p) D(p), 
        \label{eq:linFeedbackOL}
    \end{equation}
    where $Y$, $U$ and $D$ are the Laplace transforms of $y$, $u$ and $d$, respectively, $H_{yu}$ and $H_{yd}$ are the transfer functions of the plant, and $p$ is Laplace's variable. The drive is the controller output given by
    \begin{equation}
        U(p) = C(p)  (Y_*(p) - Y(p)),
        \label{eq:linFeedbackCtr}
    \end{equation}
    where $C$ is the controller's transfer function, and $Y_*$ is the Laplace transform of $y_*$. Injecting Equation~\eqref{eq:linFeedbackCtr} into Equation~\eqref{eq:linFeedbackOL} yields the closed-loop relation
    
    \begin{equation}
        \begin{array}{l}
        \displaystyle (1+H_{yu}(p)C(p)) Y(p) \\
        \displaystyle = H_{yu}(p)C(p) Y_*(p) + H_{yd}(p) D(p) 
        \end{array}
        \label{eq:linFeedbackCL}
    \end{equation}
    In general, the controller transfer function $C$ is chosen to be large compared to the other transfer functions. Hence,
    \begin{equation}
        \begin{array}{l}
        \displaystyle \lim_{|C(p)|\rightarrow\infty} H_{yu}(p)C(p) Y(p) \\
        \displaystyle = \lim_{|C(p)|\rightarrow\infty} H_{yu}(p)C(p) Y_*(p)
        \end{array}
    \end{equation}
    Hence, $Y(p)\approx Y_*(p)$, and $y(t) \approx y_*(t)$ if the controller has a large enough gain. It is thus theoretically possible to achieve perfect tracking between $y_*$ and $y$. In practice, finite gains are chosen for several reasons. First, closed-loop stability may require finite gains (and in general, unknown, unmodeled dynamics will impose a finite gain). Second, sensor noise can pollute the sensed signal $y$ and transmit fully to (and even be amplified by) the closed-loop system. Third, a large controller gain can require a large drive $u$, exceeding the actuator capabilities. There may also be other practical reasons depending on the application at hand.

    \subsection{Nonlinear control}
    \label{sec:NLControl}

        Nonlinear control encompasses feedback control of systems where either the plant or the controller is nonlinear, or both. Naturally, its development was closely related to linear feedback control and nonlinear dynamics theories. The inherent difficulties associated with nonlinear dynamical systems unsurprisingly arise (and are arguably even more problematic) in nonlinear control systems, making them a contemporary active research field. For historical perspectives, we refer to~\cite{Fuller1976a,Fuller1976,Atherton1996,Kokotovic2001,Iqbal2017}. There exists a variety of reference books on nonlinear control~\cite{slotine1991applied,khalil2002control,sepulchre2012constructive}, with one specifically aimed toward nonlinear vibration control~\cite{Wagg2015}.

        One of the challenges associated with nonlinear (and linear) control (except for adaptive and optimal nonlinear control, to a certain extent) stems from the need for a representative model of the plant for the controller design. Although substantial progress has been made in NSI, there still remains inherent challenges associated with this endeavor~\cite{Kerschen2006,Noel2017}, and the complexity associated with modeling a nonlinear system generally reverberates into nonlinear control. This motivated the use of model-less control approaches, which emerged in the context of chaos control.

        In the early nineties, a seminal article written by Ott, Grebogi and Yorke proposed a new model-less method to control the behavior of a chaotic system~\cite{Ott1990}. This OGY algorithm sparked a lot of interest in the late twentieth and early twenty-first centuries; reviews on chaos control methods can be found in~\cite{Shinbrot1993,Boccaletti2000,Fradkov2005}. A notable contribution by Pyragas later proposed a simplified continuous delayed feedback law to achieve the same goal~\cite{Pyragas1992}. By feeding back the difference between an observable and its time-delayed self, this control law was also shown to be able to stabilize the targeted periodic orbit. We note that choosing this delay to be equal to the period of the target makes the control \textit{non-invasive} on that orbit, an important notion that we shall come back to and define later. Chaos control can be considered as the ancestor of the first control-based methods discussed in this work.

\section{Experimental continuation: the pioneers and first efforts}
\label{sec:ECPioneers}

    Control-based methods combine control and continuation to trace out a complete bifurcation diagram of a nonlinear system. This idea was first proposed for chaotic systems (or chaotic at least in some region of the investigated parameter range) by combining the OGY algorithm (or its variants) to a natural parameter continuation and was coined under the term "tracking"~\cite{Carroll1992,Gills1992}. Methods combining continuous feedback and continuation to trace out bifurcation diagrams made up by the collection of periodic orbits of nonchaotic systems were introduced following the seminal works of Sieber and Krauskopf~\cite{SIEBER2007,Sieber2008} and are the focal point of the present work.

\subsection{Control-based approaches: generalities}

    \subsubsection{Why are closed-loop tests necessary?}

        Besides the academic thrill of tracing out complete experimental bifurcation diagrams, one may legitimately wonder what interest there is in attempting to follow unstable solutions in a general setting, given that they will be practically impossible to observe in the bare system. There are multiple reasons for this, and some of them are given below, the most important arguably being the first one.
        
\begin{figure}[!ht]
            \centering
            \includegraphics[width=0.539\textwidth]{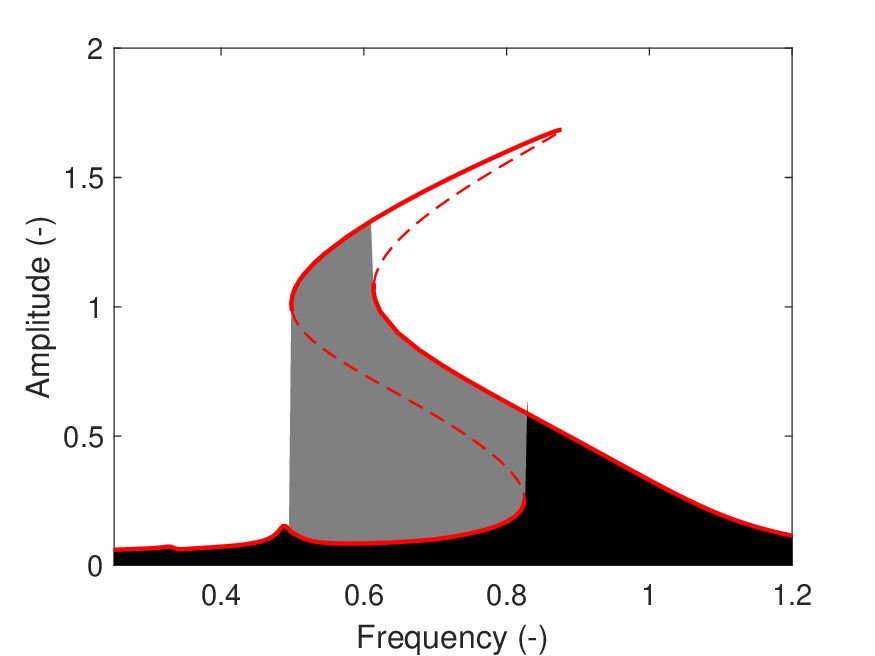}
            \caption{Bifurcation diagram of the Helmholtz-Duffing oscillator with $\zeta_0 = 0.05$, $\beta_2 = 1.9$ and a dimensionless forcing amplitude of 0.05: swept-up (\rule[.2em]{1em}{.2em}) and swept-down (\textcolor{Gray1}{\rule[.2em]{1em}{.2em}}) excitation frequency.
            The results from numerical continuation are given for reference (\textcolor{red}{\rule[.2em]{1em}{.2em}}), where full and dashed lines represent stable and unstable solutions, respectively.}
            \label{fig:HelmholtzDuffingSweep}
        \end{figure}

        \begin{enumerate}
            \item \textbf{Continuation:} unstable branches are connecting various stable branches, which might not (or not fully) be reached with open-loop system testing methods. For instance, in a Helmholtz-Duffing oscillator, SWS tests may systematically miss the highest amplitude of the NFR, as illustrated in Figure~\ref{fig:HelmholtzDuffingSweep}. There can also be stable responses that can only be accessed through unstable branches (so-called islands of stability~\cite{Neville2020}).
            \item \textbf{Global dynamics:} unstable equilibria of saddle type carry important information about the global dynamics. They lie on the separatrices of the system and provide valuable information about the basins of attraction of the system, and especially their boundaries. Loosely speaking, if a stable and an unstable equilibria are "close", one would expect the basin of attraction of the former to be small. On a related note, due to the hyperbolic nature of the dynamics around such sets, saddles also play an important role in transient dynamics.
            \item \textbf{Qualitative changes through bifurcations:} branching on a stable path through bifurcations might lead the system to a completely different type of behavior, whereas following an unstable path may yield a system response that is similar to its pre-bifurcation behavior, easing the analysis and interpretation of the results.
            \item \textbf{Data diversity:} in an NSI or model updating framework, measuring unstable solutions offers more data points for building an accurate model of the system under test.
            \item \textbf{System interconnections:} finally, when one considers nonlinear substructuring or nonlinear control, due to the feedback action enacted upon the system by its surroundings, unstable equilibria might become stable, or play a prominent role in the appearance of new stable equilibria, and thus become obviously practically relevant.
        \end{enumerate}

    \subsubsection{Common architecture and non-invasiveness}
    \label{sec:CBCommon}

\begin{figure*}[!ht]
            \centering
        	\includegraphics[scale=.6]{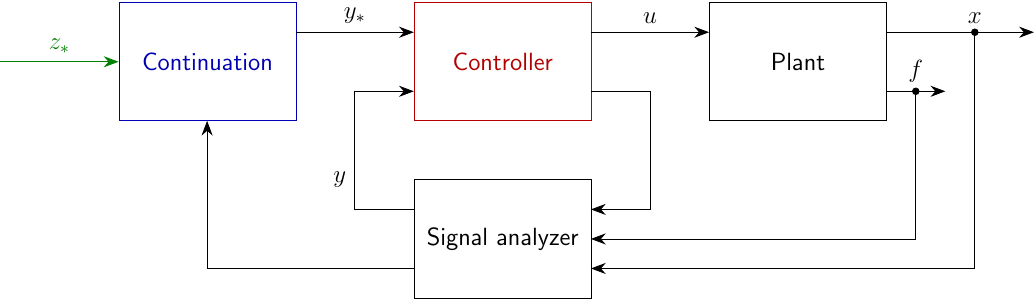}
            \caption{General schematics of a control-based method.}
            \label{fig:controlBasedSchematics}
        \end{figure*}

        Figure~\ref{fig:controlBasedSchematics} presents the general architecture of a control-based testing method (note the similarity with Figure~\ref{fig:controlSchematics}). As in open-loop testing, the plant (e.g., in the vibration testing case, a structure under test, and a set of actuators and sensors) is excited by a drive signal $u$. The response of the plant is then analyzed through any of its observables (e.g., in the vibration testing case, a displacement $q$ or a force $f$) to obtain a representative signal $y$ which is fed back to a controller. This controller compares this representative signal to a reference $y_*$, and adapts the drive signal to enforce as much as possible an equality between the reference and actual signals.
        
        The target (or a set of targets) $z_*$ is predetermined by the experimenter depending on the goals of the tests to be performed, and the selected experimental continuation method. In some methods, the continuation procedure will adapt $y_*$ through a continuation procedure. In other methods (discussed in Section~\ref{sec:ENPC}), one simply sets $y_*=z_*$, and the continuation block is trivial.
        

        A feedback controller can and generally will completely modify the dynamics of the system. Since the goal of experimental continuation is to test the plant, such a modification should be avoided. A requirement put on the control (in addition to its two main goals) is thus that it should be \textit{non-invasive}, i.e., the equilibria of the plant in the closed-loop system should be identical to those of the open-loop system. A sufficient condition to enforce this is to subject the closed-loop system to the same excitation as in open-loop conditions. Hence, the total input to the plant in closed-loop can be seen as the superposition of an open-loop excitation signal used to investigate the dynamics of the plant, and a stabilizing signal which vanishes once the system reaches the targeted open-loop equilibrium. 

\subsubsection{Stabilizing effect}
\label{ssec:stabilization}

As outlined in Section~\ref{sec:CBCommon} and particularly in Figure~\ref{fig:controlBasedSchematics}, the main difference between an open-loop test and a test with a control-based method is the presence of a feedback signal. We thus reconsider the results of Section~\ref{ssec:stability} with this difference, staying at an abstract level for generality. A way to account for this feedback is to modify Equation~\eqref{eq:ODE} as
\begin{equation}
    \dot{\mathbf{x}}(t) = \tilde{\mathbf{f}}(\mathbf{x}(t),\mathbf{u}(t;\mathbf{x}(t)),t), 
    \label{eq:controlODE}
\end{equation}
where the right-hand side $\tilde{\mathbf{f}}$ now depends on $\mathbf{u}$, the control signal. Among other possibilities, $\mathbf{u}$ may either be determined from a prescribed continuous-time evolution
\begin{equation}
    \dot{\mathbf{u}}(t;\mathbf{x}(t)) = \mathbf{g}_1(\mathbf{x}(t),\mathbf{u}(t;\mathbf{x}(t)),t),
\end{equation}
be updated by a discrete-time law, or be a direct function of the states and time
\begin{equation}
    \mathbf{u}(t;\mathbf{x}(t)) = \mathbf{g}_2(\mathbf{x}(t),t).
\end{equation}

In our setting, the control is non-invasive if, at an equilibrium of the open-loop system $\mathbf{x}_*(t)$,
\begin{equation}
    \mathbf{u}(t;\mathbf{x}_*(t)) = \mathbf{0}
    \label{eq:nonInvasiveness1}
\end{equation}
and
\begin{equation}
    \tilde{\mathbf{f}}(\mathbf{x}(t),\mathbf{0},t) = \mathbf{f}(\mathbf{x}(t),t), \qquad \forall \mathbf{x}, t. 
    \label{eq:nonInvasiveness2}
\end{equation}
If these two conditions hold, we thus have
\begin{equation}
    \tilde{\mathbf{f}}(\mathbf{x}_*(t),\mathbf{u}(t;\mathbf{x}_*(t)),t) =
    \mathbf{f}(\mathbf{x}_*(t),t) = \dot{\mathbf{x}}_*(t),
\end{equation}
i.e., $\mathbf{x}_*$ is an equilibrium of both the open- and closed-loop systems. We note that in experiments, being exactly on the equilibrium $\mathbf{x}_*$ is not possible due to, e.g., noise. Equation~\eqref{eq:nonInvasiveness1} can therefore be seen as purely theoretical and can be replaced by the condition that the control signal should be smaller and smaller as the system approaches its equilibrium state. Instead of Equation~\eqref{eq:nonInvasiveness1}, non-invasiveness can therefore be satisfied if the control is required to satisfy the Lipschitz continuity condition 
\begin{equation}
    || \mathbf{u}(t,\mathbf{x}(t)) || \leq L || \mathbf{x}(t)-\mathbf{x}_*(t)||
\end{equation}
around $\mathbf{x}_*(t)$, were $L$ is a finite positive constant and $||\cdot||$ is the Euclidean norm.

Following the same developments as in Section~\ref{ssec:stability}, we may study the stability of this equilibrium by introducing small perturbations, and linearizing the dynamics to reveal the Jacobian of the controlled system 
\begin{equation}
     \left.\dfrac{\mathrm{d} \tilde{\mathbf{f}}}{\mathrm{d} \mathbf{x} }\right|_{\mathbf{x}=\mathbf{x}_*(t)} =\left. \left(\dfrac{\partial \tilde{\mathbf{f}}}{\partial \mathbf{x}}+ \dfrac{\partial \tilde{\mathbf{f}}}{\partial \mathbf{u}} \dfrac{\mathrm{d} \mathbf{u}}{\mathrm{d} \mathbf{x} }\right)\right|_{\mathbf{x}=\mathbf{x}_*(t)}
\end{equation}
(where the total derivative is taken with fixed $t$), and, using Equations~\eqref{eq:nonInvasiveness1} and~\eqref{eq:nonInvasiveness2}, we find
\begin{equation}
    \left.\dfrac{\mathrm{d} \tilde{\mathbf{f}}}{\mathrm{d} \mathbf{x} }\right|_{\mathbf{x}=\mathbf{x}_*(t)} = \left. \left(\dfrac{\partial \mathbf{f}}{\partial \mathbf{x}}+ \dfrac{\partial \tilde{\mathbf{f}}}{\partial \mathbf{u}} \dfrac{\mathrm{d} \mathbf{u}}{\mathrm{d} \mathbf{x} }\right)\right|_{\mathbf{x}=\mathbf{x}_*(t)}.
\end{equation}
We thus observe that the Jacobian of the original system is modified by the effect of the control signal $\mathbf{u}$. Consequently, according to the discussion in Section~\ref{ssec:stability}, the stability of $\mathbf{x}_*$ may be altered by the action of the controller. Unfortunately, going beyond this discussion while staying in a general framework may prove complicated, first because of the non-trivial relation between the Jacobian and the Floquet multipliers, but also because of the different nature of the control signals used by the different approaches.

\begin{figure*}[!ht]
    \centering
     \begin{subfigure}{\textwidth}
        \centering
        \includegraphics[width=\textwidth]{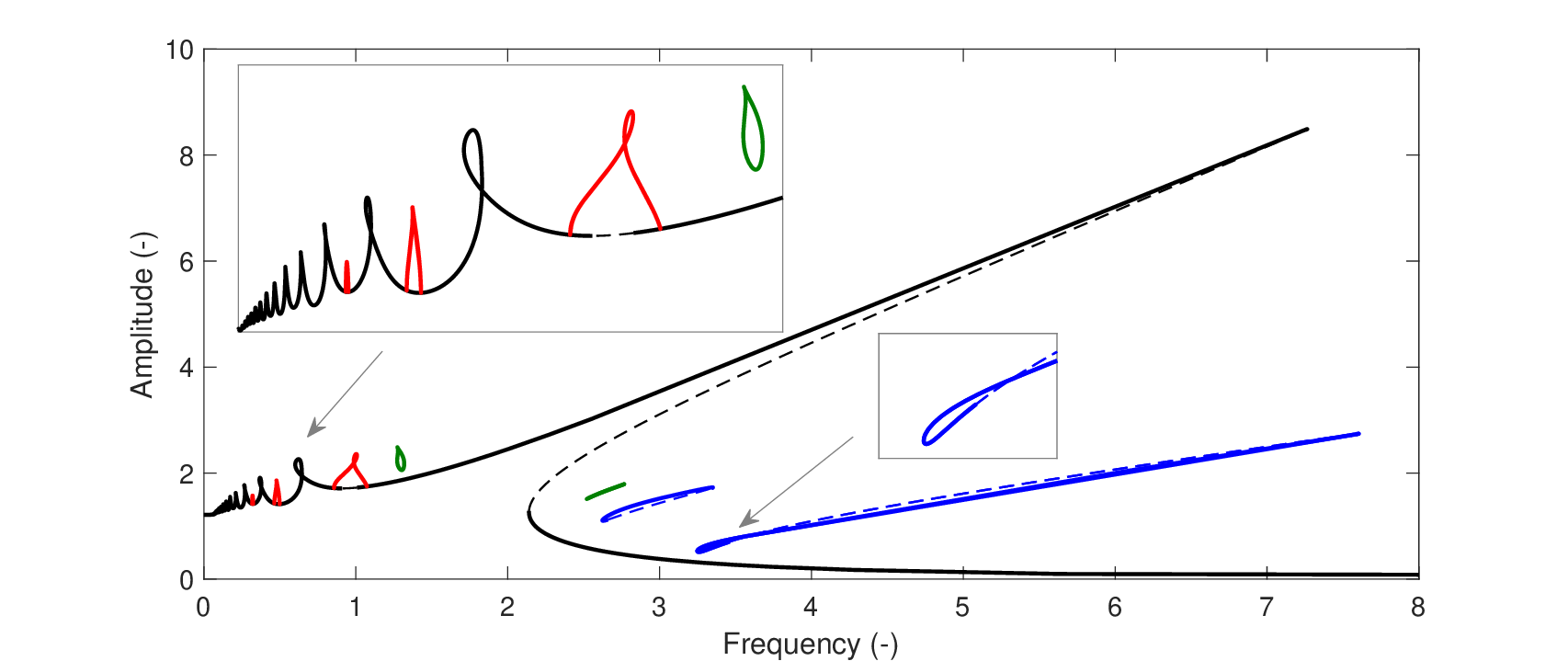}
        \caption{}
        \label{sfig:ControlledDuffing1}
    \end{subfigure}
    \begin{subfigure}{\textwidth}
        \centering
        \includegraphics[width=\textwidth]{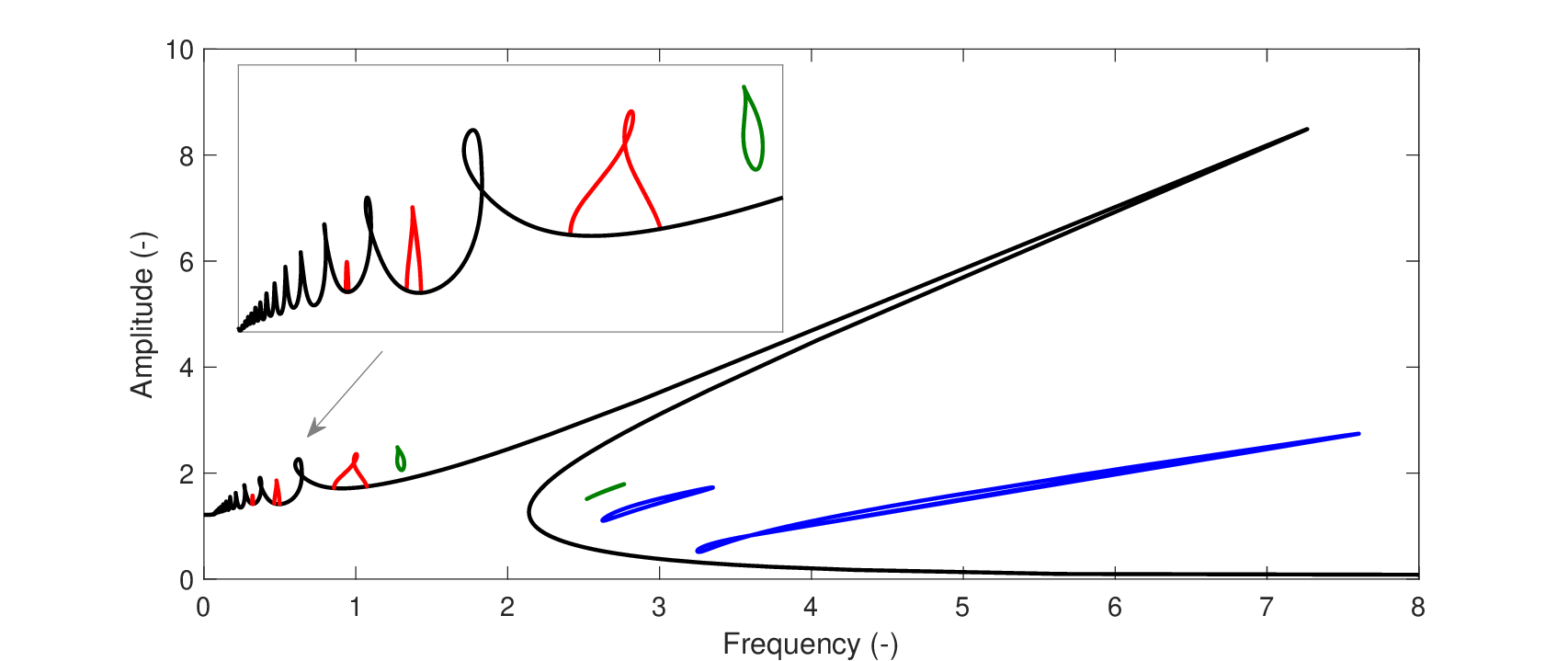}
        \caption{}
        \label{sfig:ControlledDuffing2}
    \end{subfigure}
    \caption{Bifurcation diagram of the Duffing oscillator with a derivative controller with $\bar{k}_\mathrm{d} = 0.2$~\subref{sfig:ControlledDuffing1} and $\bar{k}_\mathrm{d}=2$~\subref{sfig:ControlledDuffing2}, with a dimensionless forcing amplitude of 3: main branch of the NFR (\rule[.2em]{1em}{.2em}) bifurcating even superharmonic resonances (\textcolor{red}{\rule[.2em]{1em}{.2em}}), isolated subharmonic resonances (\textcolor{blue}{\rule[.2em]{1em}{.2em}}) and isolated ultrasubharmonic resonances (\textcolor{MatlabDarkGreen}{\rule[.2em]{1em}{.2em}}). Full and dashed lines represent (closed-loop) stable and unstable solutions, respectively.}
    \label{fig:ControlledDuffingBoAs}
\end{figure*}

We may nevertheless illustrate this section with a specific example. Figure~\ref{fig:ControlledDuffingBoAs} presents the counterpart of Figure~\ref{fig:Duffing_Res_Feature} when the Duffing oscillator is stabilized with a differential controller with a fixed reference signal $\bar{q}_*$ correctly ensuring the non-invasiveness condition, i.e., $\bar{f} = \bar{k}_\mathrm{d} (\bar{q}'_* - \bar{q}')$, where $\bar{f}$ is sinusoidal and $\bar{q}'$ corresponds to an open-loop (stable or unstable) equilibrium solution to that forcing (as done in the control-based continuation method, CBC, introduced in the next section). With the differential gain $\bar{k}_\mathrm{d}$, the open-loop unstable branches can be partially (Figure~\ref{sfig:ControlledDuffing1}) or totally (Figure~\ref{sfig:ControlledDuffing2}) stabilized. When the gain is not large enough, the remaining unstable solution co-exists with other stable attractors for which the control is invasive (keeping $\bar{q}_*$ fixed to make the control non-invasive on the unstable solution). Interpretations of the stabilizing effects of the CBC in terms of folding and closed-loop NFR can be found in~\cite{Abeloos2022Thesis,Tatzko2023} and~\cite{Hayashi2024}, respectively.

The stabilizing and non-invasive properties of the control signal may intuitively seem like conflicting requirements to the newcomer: how can this signal stabilize anything if it is supposed to vanish? The key element of the answer lies in understanding that while it does indeed vanish when the system is \textit{exactly} at the equilibrium, it does not when the system is \textit{near} this equilibrium. To illustrate this, a simple and evocative example is that of an inverted pendulum on a cart, where the cart motion is controlled to ensure that the pendulum stays in its upright position. This control therefore stabilizes what would normally be an unstable equilibrium and acts against any perturbation. It is also non-invasive, because the cart stays still when the pendulum is upright in the (theoretical) absence of perturbations. Another illustrative example is given in~\cite{Bureau2013,Bureau2014Thesis}.

    \subsection{Control-based continuation}
    \label{ssec:CBC}

        \begin{figure*}[!ht]
            \centering
        	\includegraphics[scale=.6]{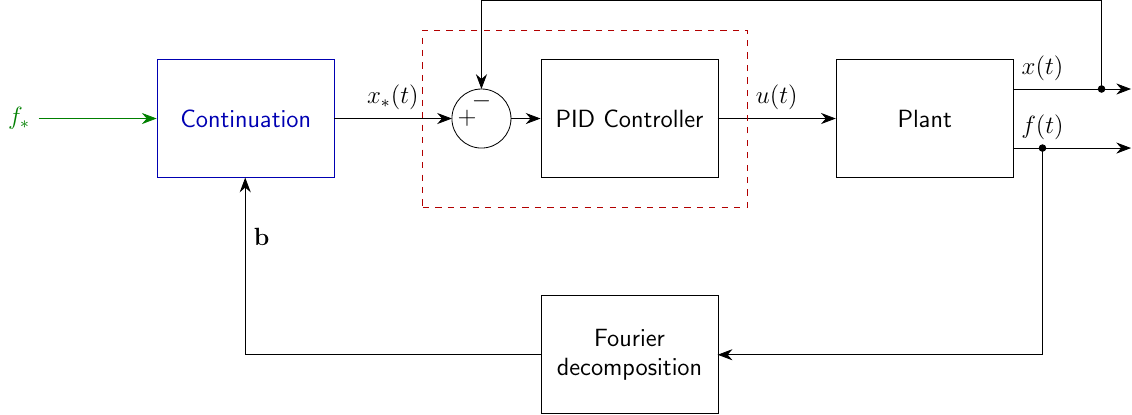}
            \caption{Schematics of the CBC-FD method; the red dashed frame indicates a controller block.}
            \label{fig:CBCSchematics}
        \end{figure*}

        We now introduce the first experimental continuation method, namely CBC. Inspired by previous works~\cite{Ott1990,Pyragas1992,SIETTOS2004}, it was proposed by Sieber and Krauskopf~\cite{SIEBER2007,Sieber2008} in a very general setting by seeking an experimental equivalent to pseudo-arclength continuation. The two key innovations of this work were (\textit{i}) to use a PID-type controller\footnote{We note that other types of controller can be used, as long as they satisfy the non-invasiveness requirement.} to stabilize the unstable responses of the system under test and (\textit{ii}) to formulate a zero problem with measurable functions that indicate a non-invasive control. Figure~\ref{fig:CBCSchematics} schematizes this original version of the CBC method.

        This method is organized around the reference signal $x_*$, which is set so as to make the control non-invasive by prescribing the force applied on the structure $f$ to a preset profile $f_*$ (alternatively, there can be versions where the drive $u$ is specified instead of $f$), i.e., 
        \begin{equation}
            f(t;\mathbf{w}_*) - f_*(t) \approx 0.
            \label{eq:CBCFDNoninvasiveness}
        \end{equation}
        This non-invasiveness requirement is dealt with by using truncated Fourier series of the signals of interest,
        \begin{equation}
            x_*(t) \approx \mathbf{h}(t)\mathbf{w}_*,
            \label{eq:xstarFourier}
        \end{equation}
        where the column vector $\mathbf{w}_*$ gathers the Fourier coefficients of $x_*$, and $\mathbf{h}(t)$ is the associated vector of $2h+1$ harmonic functions
        \begin{equation}
            \mathbf{h}^T(t) = \begin{bmatrix}
                1 \\ \sin(\omega t) \\ \vdots \\ \sin(h\omega t) \\ \cos(\omega t) \\ \vdots \\ \cos(h\omega t)  
            \end{bmatrix},
            \label{eq:harmonicFunctions}
        \end{equation}
        where superscript $T$ denotes the transposition operator. By adjusting the Fourier coefficients $\mathbf{w}_*$ of the reference signal fed to the controller, the difference between the excitation $f$ and a desired excitation profile $f_*$ is sought to be as close to zero as possible. A well-posed problem can be obtained if the force profile is also discretized with Fourier series, i.e.,
        \begin{equation}
            f(t;\mathbf{w}_*) \approx \mathbf{h}(t)\mathbf{b}(\mathbf{w}_*,\omega),
        \end{equation}
        and Equation~\eqref{eq:CBCFDNoninvasiveness} becomes
        \begin{equation}
            \mathbf{b}(\mathbf{w}_*,\omega) - \mathbf{b}_* = \mathbf{0},
            \label{eq:CBCFDNoninvasivenessFourier}
        \end{equation}
        with $\mathbf{b}_*$ the desired force profile in a Fourier basis. For instance, if $f_*(t) = \sin(\omega t)$, $\mathbf{b}_* = \mathbf{e}_2$, where $\mathbf{e}_i$ is the $i^{th}$ canonical basis vector of $\mathbb{R}^{2h+1}$. Equation~\eqref{eq:CBCFDNoninvasivenessFourier} defines a zero problem with a nonlinear function $\mathbf{b}(\mathbf{w}_*,\omega)$, where the number of unknowns ($\mathbf{w}_*$ and $\omega$) exceeds the number of equations by one. To solve this problem, a numerical pseudo-arclength continuation algorithm can be used. Since this algorithm relies on Newton-type methods to find the zeros of a residual function, it requires derivatives, which are not known in an experimental setup, but can be approximated using finite differences. Hence, in this work, this original version of the CBC is called CBC-FD. From the output of this continuation method, $x_*(t)$ can be reconstructed from $\mathbf{w}_*$ using Equations~\eqref{eq:xstarFourier} and~\eqref{eq:harmonicFunctions}.

\begin{figure*}[!ht]
            \centering
            \includegraphics[width=0.6\textwidth]{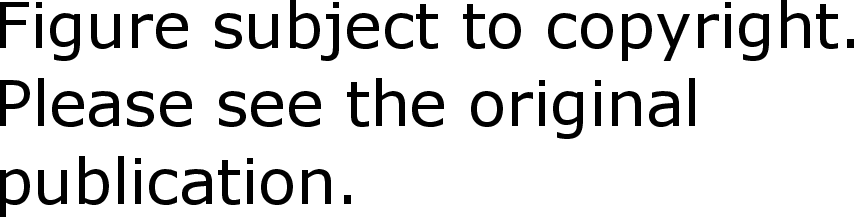}
            \caption{Primary resonance of a nonlinear energy harvester obtained experimentally with the CBC-FD method~\cite{Barton2012}. Reprinted from~\cite{Barton2012}, Copyright (2012), with permission from SAGE.}
            \label{fig:Barton2012}
        \end{figure*}

        In a vibration testing context, as the CBC-FD method prescribes the forcing profile, so is its amplitude, and in general this method provides the NFR of a system, although variants can be conceptualized to provide other types of response curves. Figure~\ref{fig:Barton2012} presents the main resonance of an energy harvester from~\cite{Barton2012}, showing that this method is indeed able to obtain a full branch of an NFR.

        \subsubsection{CBC-FD applications}
            
            The first experimental validation of CBC-FD was with a vertically forced pendulum~\cite{Sieber2008b,Sieber2011}. An alternative way to compute the residual function using a time delay (similarly to~\cite{Pyragas1992}) to handle non-fundamental harmonics was also proposed in~\cite{SIEBER2007}, and later used on energy harvesting systems~\cite{Barton2011,Barton2012}. This extended time delay approach is simpler than the CBC-FD algorithm presented herein, but has less stabilization capabilities. The CBC-FD method was validated on a system with impacts~\cite{Bureau2013} and on a geometrically nonlinear beam~\cite{Kleyman2024}. Interestingly, similar ideas have recently been introduced in quasi-static tests to avoid snapthrough and snapback phenomena~\cite{VanIderstein2019,Neville2020,Shen2021}, in real-time hybrid testing~\cite{Beregi2023p}, and in testing procedures for microwave circuits~\cite{Melville2017}.

            \begin{figure*}[!ht]
                \centering
                \includegraphics[width=0.6\textwidth]{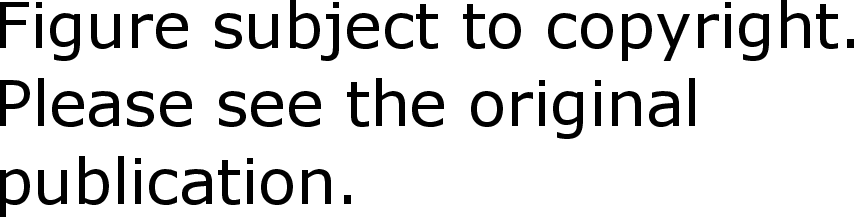}
                \caption{First complete isola (1:3 subharmonic resonance of an impact oscillator) obtained experimentally thanks to two runs of the CBC-FD method (full lines). Data from SWS-up ($+$) and SWS-down ($\circ$) is indicated for reference. Stable and unstable responses are indicated in black and gray, respectively. Reprinted from~\cite{Bureau2014}, Copyright (2014), with permission from Elsevier.}
                \label{fig:Bureau2014}
            \end{figure*}

            Strategies to determine the stability of the found periodic orbit relying on temporarily switching the control off were investigated in~\cite{Bureau2014}, and the effect of noise was thoroughly studied in~\cite{Schilder2015}. To the authors knowledge, Bureau et al~\cite{Bureau2014} were also the first to measure a complete isola experimentally, as shown in Figure~\ref{fig:Bureau2014}. The possibility to use B-splines (instead of Fourier decomposition) was recently explored by~\cite{Blyth2023} for slow-fast systems. An in-depth dissertation on CBC-FD can be found in~\cite{Bureau2014Thesis}, and a software package, CONTINEX, is available in COCO~\cite{COCO}.

        \subsubsection{Limitations}

            The CBC-FD method leverages pseudo-arclength continuation. Folding is almost guaranteed to be absent, which makes this method inherently more robust (continuation-wise) than traditional SWS tests. Unfortunately, this robustness comes at the expense of substantial complexity and longer tests due to the operations required by the method, such as the computation of a Jacobian through finite differences. This motivates the derivative-free methods that are treated in the next sections.

\section{Derivative-free experimental natural parameter continuation}
\label{sec:ENPC}

We now treat the set of state-of-the-art methods that can be classified as "derivative-free", in the sense that they do not rely (explicitly) on a Jacobian of the residual function, and thus do not require finite differences and the likes (e.g. Broyden updates). We first present them in a unified framework, showing that they can be seen as derivative-free experimental natural parameter continuation methods (cf. Section~\ref{ssec:NPC}), before reviewing them.

\subsection{Continuation}
\label{ssec:CBcontinuation}

Similarly to its numerical counterpart presented in Section~\ref{sec:numericalContinuation}, experimental continuation can be performed in two main variants. The first one, arclength-type continuation, is used in the CBC-FD method as presented in Section~\ref{ssec:CBC}. The second one is natural parameter continuation, for which a series of relevant parameters can be chosen, generally based on how representative they are of the excitation or response of the system.

We now specialize the discussion to nonlinear vibrations, given that most of the methods treated herein were developed in this area. The discussion can however be translated to other disciplines by considering relevant variables. To represent the state of a vibrating structure, four distinct variables are generally used:
\begin{itemize}
    \item $f$, an amplitude of the external forcing applied to the system,
    \item $\omega$, the frequency of this forcing,
    \item $a$, the amplitude of some harmonic of the response at some point of the structure (typically, the same harmonic as the forcing of the collocated degree of freedom) and
    \item $\theta$, the phase of this harmonic.
\end{itemize}
There are of course more variables than this, but these four are of particular interest in control-based methods. In the most basic set-ups described in Section~\ref{ssec:openLoop}, the operator controls $f$ and $\omega$, and $a$ and $\theta$ are naturally determined by the evolution of the system (or, in a mathematical setting, they are the unknowns of the problem, whereas $f$ and $\omega$ are given parameters). For linear systems, this evolution is often represented with a fixed $f$ and a varying $\omega$ in Bode diagrams. Due to the invariance of the linear properties with respect to $f$, the results are usually normalized with the latter. For nonlinear systems however, variations of $f$ retain their full importance, resulting in more than a mere scaling of the response.

Control-based methods change the paradigm of Section~\ref{ssec:openLoop}, in the sense that $f$ and/or $\omega$ can be automatically adjusted to prescribe $a$ and/or $\theta$ as parameters set by the experimenter. But in all cases, there are always two parameters among the four variables that can be set by the operator, while the two remaining variables adapt according to the dynamics of the system and can be used to represent its response.

\begin{figure*}[!ht]
    \centering
    \begin{subfigure}{.45\textwidth}
        \centering
        \includegraphics[width=\textwidth]{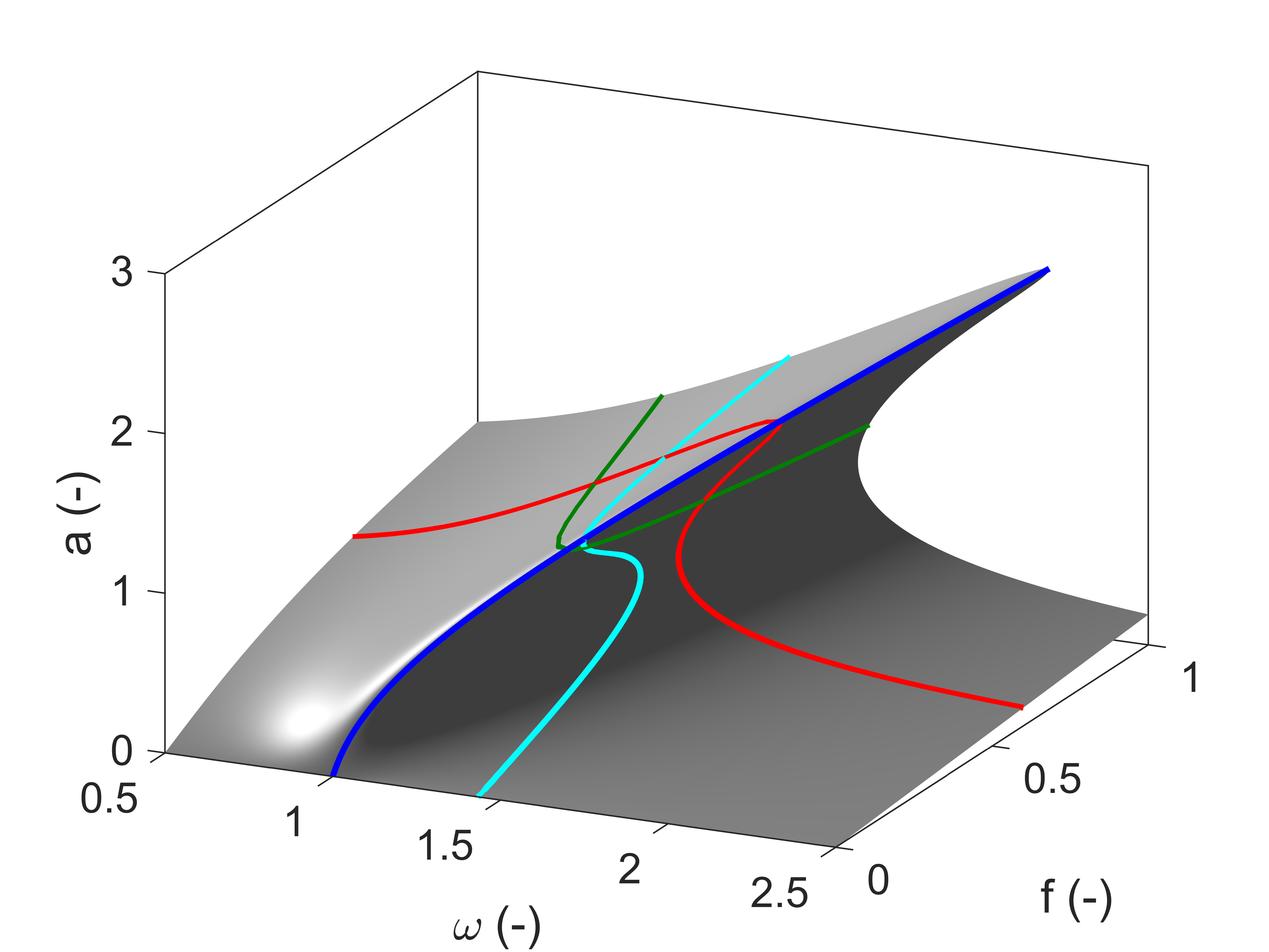}
        \caption{}
        \label{sfig:DuffingManifold}
    \end{subfigure}
    \begin{subfigure}{.45\textwidth}
        \centering
        \includegraphics[width=\textwidth]{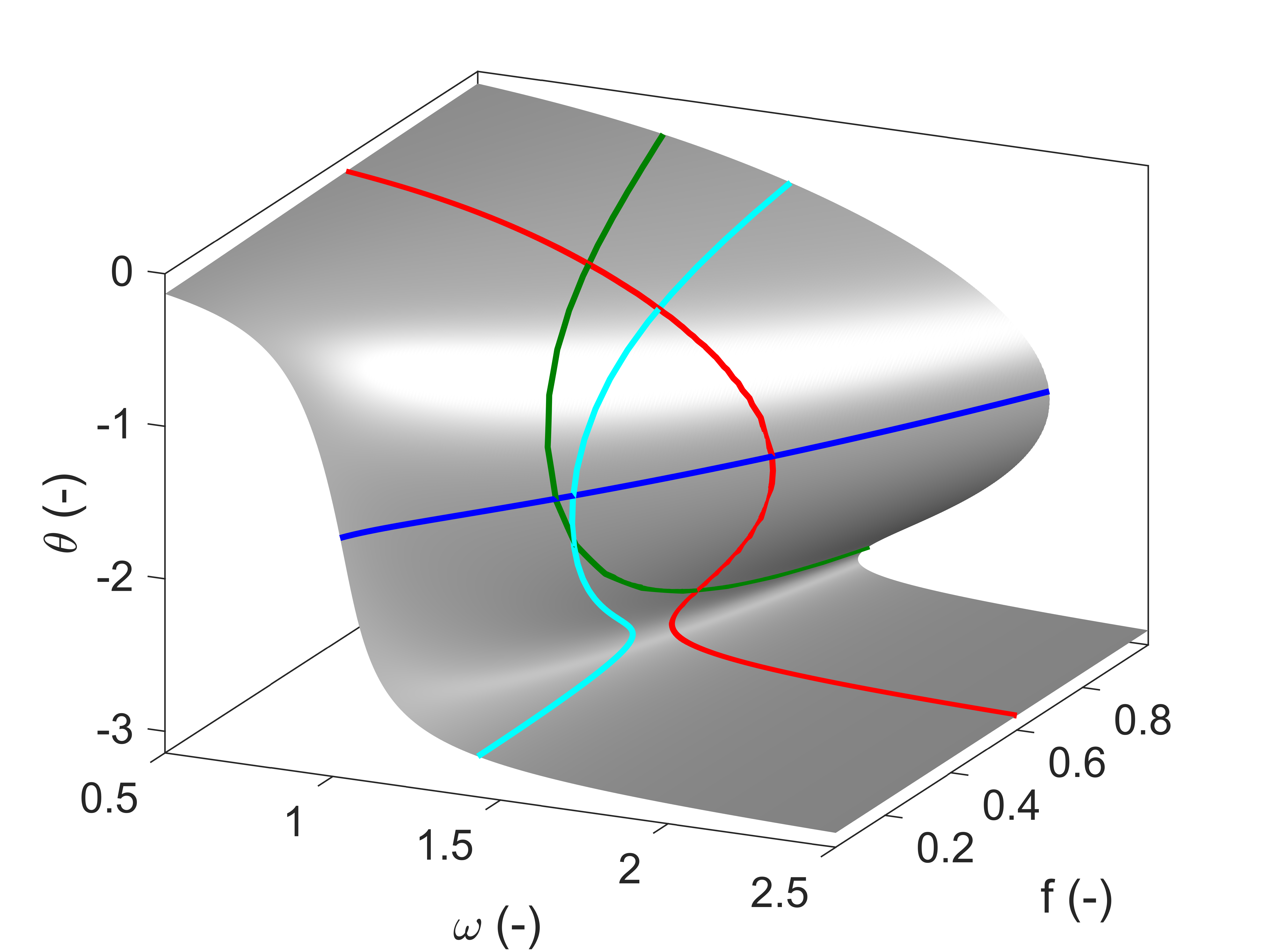}
        \caption{}
        \label{sfig:DuffingManifold2}
    \end{subfigure}
    
    \begin{subfigure}{.45\textwidth}
        \centering
        \includegraphics[width=\textwidth]{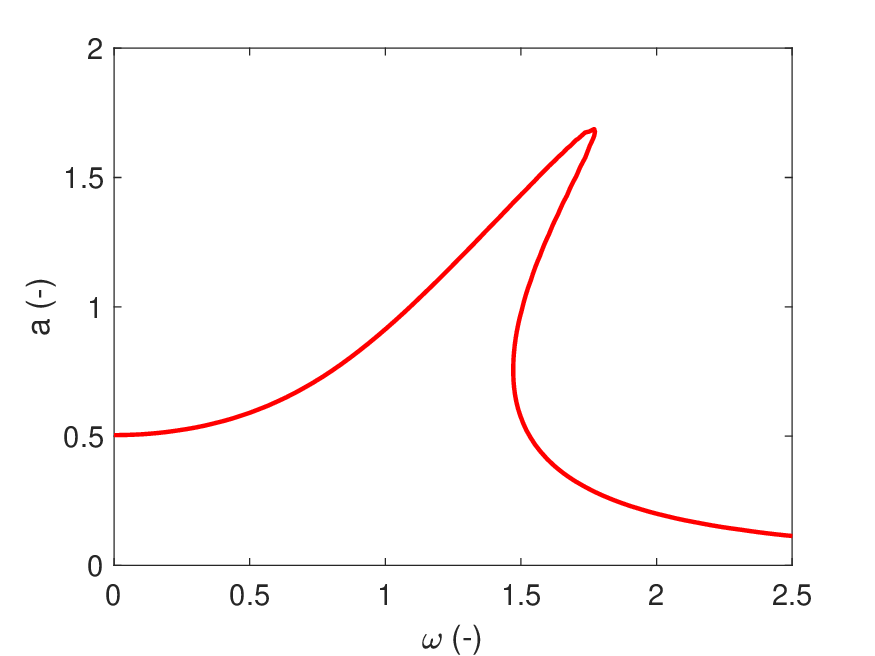}
        \caption{}
        \label{sfig:DuffingCurve1}
    \end{subfigure}
    \begin{subfigure}{.45\textwidth}
        \centering
        \includegraphics[width=\textwidth]{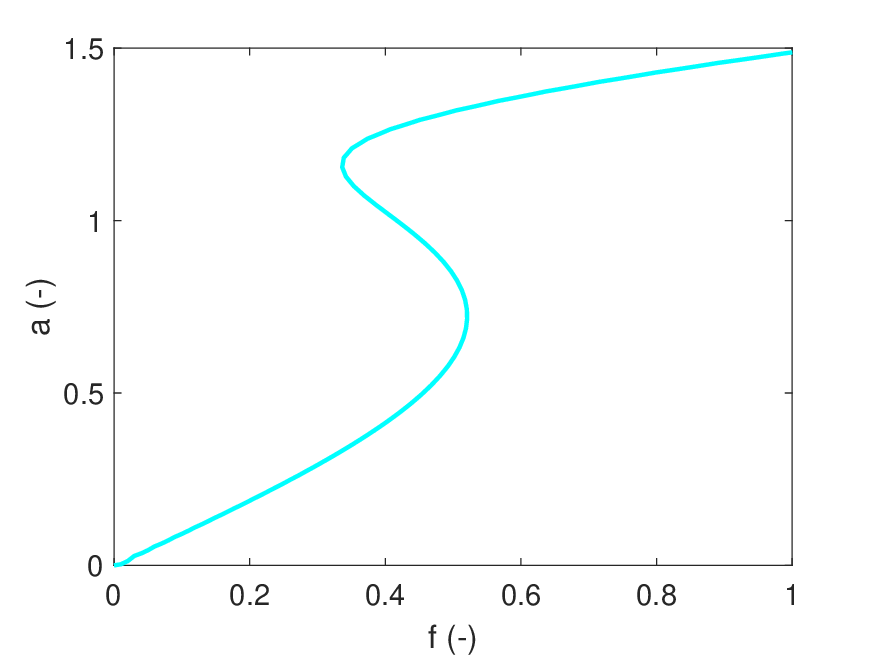}
        \caption{}
        \label{sfig:DuffingCurve2}
    \end{subfigure}
    
    \begin{subfigure}{.45\textwidth}
        \centering
        \includegraphics[width=\textwidth]{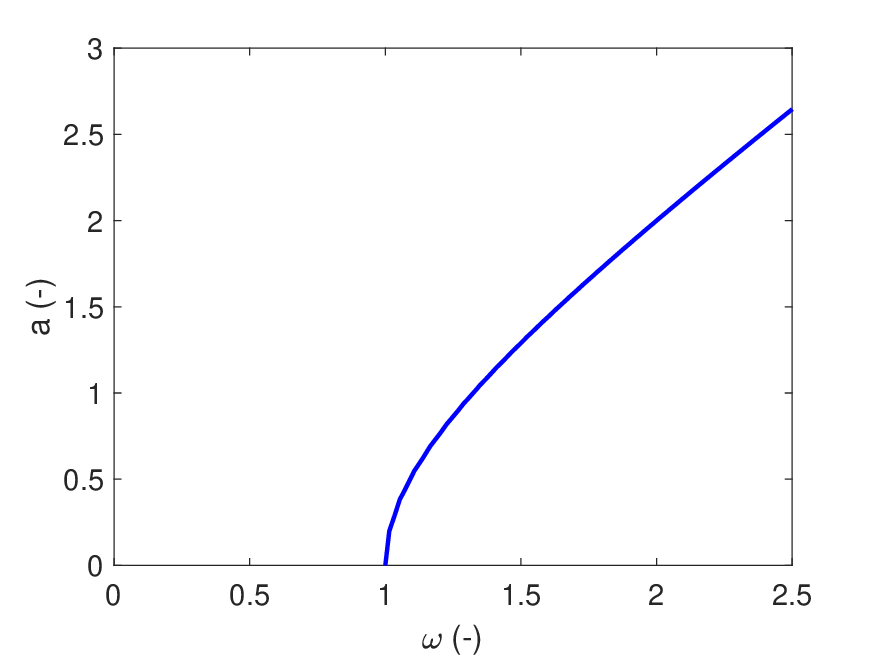}
        \caption{}
        \label{sfig:DuffingCurve3}
    \end{subfigure}
    \begin{subfigure}{.45\textwidth}
        \centering
        \includegraphics[width=\textwidth]{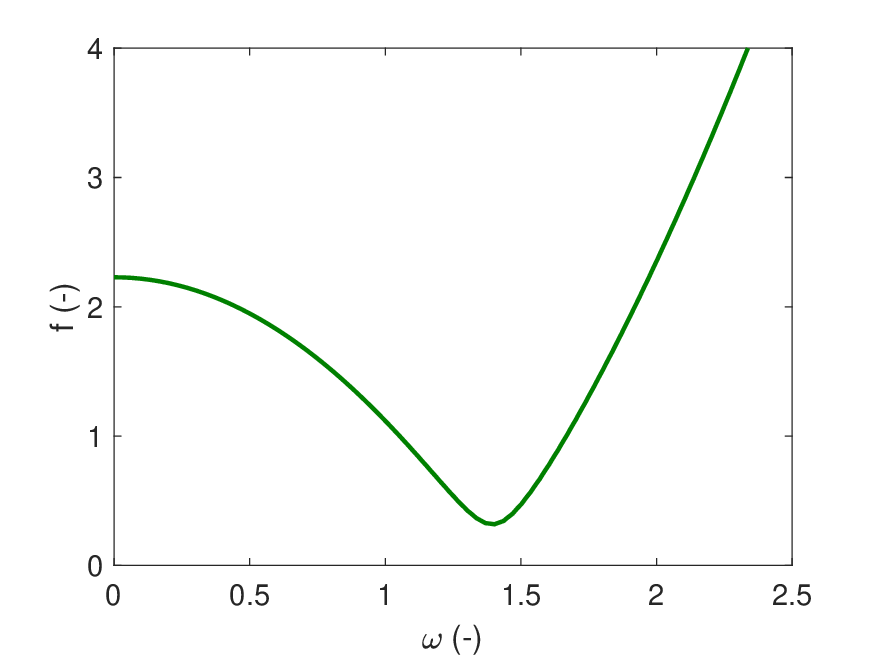}
        \caption{}
        \label{sfig:DuffingCurve4}
    \end{subfigure}
    \caption{Response manifold of a Duffing oscillator under tonal excitation in the ($\omega,f,a$)~\subref{sfig:DuffingManifold} and ($\omega,f,\theta$)~\subref{sfig:DuffingManifold2} spaces, with constant-$f$ (NFR, \textcolor{red}{\rule[.2em]{1em}{.2em}},~\subref{sfig:DuffingCurve1}) constant-$\omega$ (S-curve, \textcolor{MatlabCyan}{\rule[.2em]{1em}{.2em}},~\subref{sfig:DuffingCurve2}), constant-$\theta$ (backbone curve, \textcolor{blue}{\rule[.2em]{1em}{.2em}},~\subref{sfig:DuffingCurve3}) and constant-$a$ (harmonic force curve, \textcolor{MatlabDarkGreen}{\rule[.2em]{1em}{.2em}},~\subref{sfig:DuffingCurve4}) curves.}
    \label{fig:DuffingManifold}
\end{figure*}

In addition, one of the two chosen parameters is usually kept constant during a test, while the other is varied. The fixed parameter determines the type of response, whereas the changing parameter is most of the time used as the abscissa when plotting the response of the system. For illustration, Figure~\ref{fig:DuffingManifold} represents all possible responses of a Duffing oscillator under tonal excitation as a two-dimensional manifold in the three-dimensional ($\omega$, $f$, $a$) and ($\omega$, $f$, $\theta$) spaces, together with response curves obtained when one of these parameters is fixed (which can be seen as "slices" of that manifold). We note that the terminology "backbone curve" refers herein to a special case of a constant-$\theta$ curve where $\theta=-\pi/2$, which is often close to the amplitude maximum of the NFR~\cite{Cenedese2020}.

\begin{table*}[h]
\caption{Classification of the existing experimental  natural parameter continuation methods based on the imposed (fixed and parametrization) parameters. Cells in white, light and dark gray require zero, one and two control procedures, respectively, and cells in black are impossible. The brackets indicate which type of curve is obtained.}\label{tab:continuationClassification}%
\def\arraystretch{2}
\footnotesize \sf
\begin{tabular}{@{}|c|m{3cm}|m{3cm}|m{3cm}|m{3cm}|@{}}
\hline
 \diagbox{Varied}{Fixed}  & $f$ & $\omega$ & $\theta$ & $a$\\
\hline
 $f$ & \cellcolor{black}   & SWS (S-curve) & \cellcolor{black!20!white} PLL (Backbone curve)  & \cellcolor{black!20!white}   \\
\hline
 $\omega$ & SWS (NFR)     & \cellcolor{black}   &  \cellcolor{black!20!white}  & \cellcolor{black!20!white} RCT (Harmonic force curve)  \\
\hline
 $\theta$ & \cellcolor{black!20!white} PLL (NFR) & \cellcolor{black!20!white} & \cellcolor{black} &   \cellcolor{black!50!white} \\
\hline
 $a$ & \cellcolor{black!20!white} & \cellcolor{black!20!white} SCBC (S-curve) & \cellcolor{black!50!white} CBC, CBPLL, PLL-A  (Backbone curve) & \cellcolor{black} \\
\hline
\end{tabular}
\end{table*}

With this view in mind, we can classify most of the experimental continuation methods (that shall be presented hereafter in this section) following Table~\ref{tab:continuationClassification}. In this table, we can distinguish methods that require no, one or two control procedures. By control procedures, we mean, in a broad sense, any procedure that will adjust the input to the systems, using continuous feedback control (for SCBC and PLL) or discrete iterative methods (for SCBC and RCT). Methods simultaneously imposing $f$ and $\omega$ require no control procedure, and are thus always open-loop. Most of the existing control-based methods require just one control procedure, with the exception of the SCBC method applied to the search of specific phase resonance~\cite{Renson2017}, SCBC combined with PLL (CBPLL)~\cite{Abeloos2022Thesis}, or PLL with amplitude control (PLL-A)~\cite{Woiwode2024}.

Fixing a parameter while varying another one to parametrize the response makes the methods listed in Table~\ref{tab:continuationClassification} inherently equivalent to natural parameter continuation. As such, they are all susceptible to folding issues. This is well-known for the open-loop procedures, but should also be kept in mind when using the control-based methods listed in Table~\ref{tab:continuationClassification}. 

We can now present in more detail the main derivative-free experimental natural parameter continuation methods, namely, SCBC, PLL and RCT.

        \subsection{Simplified control-based continuation}

        \begin{figure*}[!ht]
            \centering
        	\includegraphics[scale=.6]{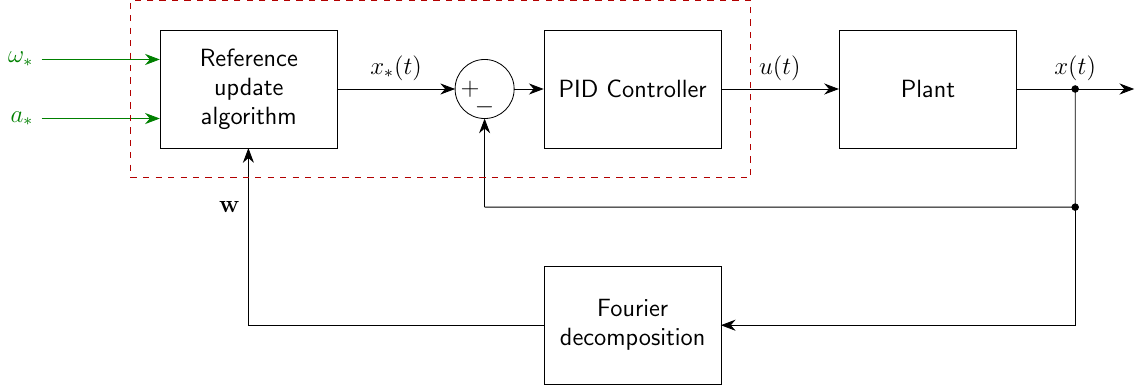}
            \caption{Schematics of the SCBC method; the red dashed frame indicates a controller block.}
            \label{fig:CBCASchematics}
        \end{figure*}
        
        A substantial simplification of the original CBC algorithm was proposed by Barton et al~\cite{Barton2013} and is schematized in Figure~\ref{fig:CBCASchematics}. They noted that, since the PID controller is linear and time-invariant, a single-harmonic input $x_*-x$ to the PID yields a single-harmonic drive $u$. They thus proposed to set the fundamental amplitude of the reference signal to a target value $a_*$,
        \begin{equation}
            x_*(t) = a_* \sin(\omega t) + \mathbf{h}_{\mathrm{nf}}(t)\mathbf{w}_{*,\mathrm{nf}}
        \end{equation}
        and to set its other harmonic coefficients, $\mathbf{w}_{*,\mathrm{nf}}$, equal to those of $x$, thus satisfying the non-invasiveness requirement. Thus, if $\mathbf{w}$ are the Fourier coefficients of $x$ (and $\mathbf{w}_{\mathrm{nf}}$ its non-fundamental Fourier coefficients), Equation~\eqref{eq:CBCFDNoninvasivenessFourier} is changed into
        \begin{equation}
            \mathbf{w}_{*,\mathrm{nf}} - \mathbf{w}_{\mathrm{nf}}(\mathbf{w}_*,\omega) = 0.
            \label{eq:CBCANoninvasiveness}
        \end{equation}
        They also simplified the continuation algorithm by noting that Equation~\eqref{eq:CBCANoninvasiveness} can be solved with a simple Picard iteration scheme, thereby removing the need for experimental finite differences. We refer to this simplified implementation as simplified CBC (SCBC).

        The main reason for using $a_*$ as a bifurcation parameter instead of $\omega$ is that it removes the aforementioned folding issue in most cases of interest, provided the controller's gain are properly chosen. Indeed, it can be observed in Figure~\ref{sfig:DuffingCurve2} that S-curves can be parametrized with the amplitude $a$. It should be noted that the actual amplitude $a$ will not be exactly equal to the target amplitude $a_*$, but will tend to it as the controller gains are increased. We finally note that, thanks to its broad theoretical framework, the SCBC can be performed with another bifurcation parameter than $a_*$, but the latter is preferred in the vast majority of cases.
        
\begin{figure*}[!ht]
            \centering
            \includegraphics[width=0.6\textwidth]{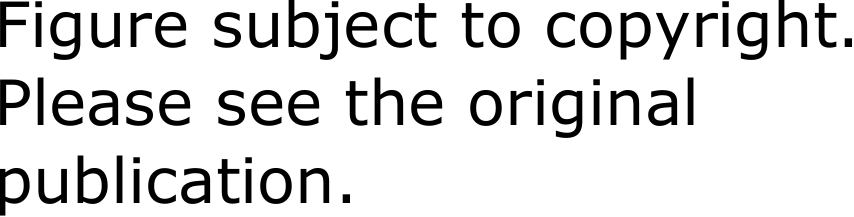}
            \caption{Manifold of the response of a nonlinear energy harvester obtained with the SCBC method. Black and gray curves represent S-curves and NFRs, respectively. The dark gray part of the manifold is open-loop unstable. Reprinted figure with permission from~\cite{Barton2013}. Copyright (2013) by the American Physical Society.}
            \label{fig:Barton2013}
        \end{figure*}

        A SCBC test is performed at a constant frequency $\omega$, by progressively varying the value of the target amplitude $a_*$, obtaining so-called S-curves. Collecting these curves for different target amplitudes gives a surface (in the $(\omega,f,a_*)$ space) that can be sliced at a constant force level to obtain an NFR. This procedure essentially amounts to fitting experimental data and performing interpolation on the fit. Figure~\ref{fig:Barton2013}, reproduced from~\cite{Barton2013}, illustrates such a process. 

        \subsubsection{SCBC applications}
            
            The SCBC method was applied to systems with friction~\cite{Kleyman2020} and beams~\cite{Hayashi2024}, and advanced versions were developed by Renson et al, enabling the tracking of backbone curves~\cite{Renson2016}, fold bifurcations~\cite{Renson2017} and internal resonances~\cite{Renson2019}. A simplified implementation using adaptive filters and circumventing the use of Picard iterations was proposed in~\cite{Abeloos2021}, thereby speeding up the process. Interestingly, the SCBC approach can also be used in a purely numerical setting to solve challenging problems~\cite{Beregi2022}.

            Noise can be handled using local Gaussian process regression~\cite{Renson2019b}. From a NSI perspective, noise effects were shown to be smaller with CBC than with traditional testing methods~\cite{Beregi2021}, and successes using machine learning and/or Gaussian regression processes using SCBC results were demonstrated in~\cite{Beregi2023,Lee2023}. Procedures to determine the open-loop stability were developed based on data-driven linearization~\cite{Barton2017,Dittus2023}. A systematic method to choose the gains of the PD controller to guarantee the absence of folding was also developed for a Duffing oscillator~\cite{Tatzko2023}. Finally, we note that the SCBC found applications beyond structural mechanics, such as aeroelastic systems~\cite{Tartaruga2019}, biochemical processes~\cite{DeCesare2022}, micro-electromechanical systems~\cite{Hayashi2024c} and pedestrian flow control~\cite{Panagiotopoulos2022,Panagiotopoulos2023}. An in-depth dissertation on the SCBC method can be found in~\cite{Abeloos2022Thesis}.
        
    \subsection{Phase-locked loops}

        Resonance conditions are oftentimes complicated to reach in nonlinear systems due to their closeness to unstable regions. This is what motivated Babitsky~\cite{Babitsky1995} to introduce velocity feedback, creating auto-resonant systems that automatically maintain themselves at resonance. This idea was later extended to continuation by Sokolov and Babitsky~\cite{SOKOLOV2001} who proposed the concept of phase-controlled systems: the phase shift between the excitation and the response $\theta$ can be used as a parameter for the NFR. An  experimental implementation of this concept was realized with a delayed acceleration feedback by Mojrzisch, Wallaschek and Bremer~\cite{Mojrzisch2012}. The unfolding of the NFR also comes with feedback stabilization, offering the possibility to obtain its unstable parts. 
        
        Delayed displacement, velocity or acceleration feedback is very simple to realize in practice. It can be used, e.g., for (linear) experimental modal analysis~\cite{Davis2018}. However, it is generally invasive and stabilization can be tricky to achieve. The latter point has recently been discussed in~\cite{Scheel2022} and treated by measuring the structural response at multiple locations. An alternative to this type of approach is to use phase-locked loops (PLLs).

        Invented nearly one hundred years ago, PLLs are extensively used in electronics for signal synchronization purposes and are actually thought to be one of the largest engineering application of feedback control in the world~\cite{Abramovitch2002}. Given a reference signal, PLLs are able to adjust the frequency and phase of another signal so as to lock the phase between these two signals to a prescribed value. Coupling PLLs with (electro-) mechanical systems was achieved even prior to the works of Babitsky for structural control as a solution to modal spillover~\cite{Balas1978,Niezrecki1997} and to track resonance in fatigue testing~\cite{Connally1993}. Phase parametrization of frequency response functions (FRFs) using PLLs was proposed in~\cite{Kern2010,Kern2012}, and the first time PLLs were used in an experimental continuation framework for nonlinear structural systems was the relatively recent work of Peter and Leine~\cite{Peter2017}. The main advantage that PLLs have over delayed feedback is their natural non-invasiveness.

        \begin{figure*}[!ht]
            \centering
        	\includegraphics[scale=.6]{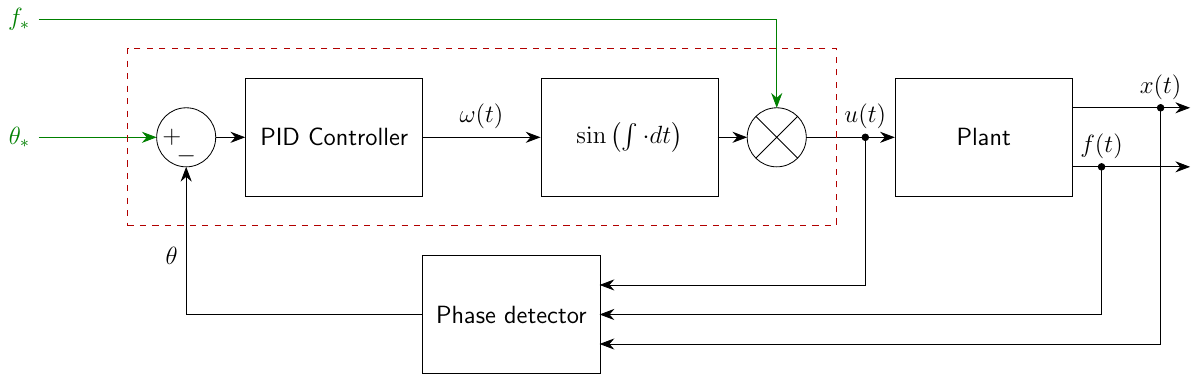}
            \caption{Schematics of the PLL method; the red dashed frame indicates a controller block.}
            \label{fig:PLLSchematics}
        \end{figure*}

        A PLL test set-up consists of a structure under test, a phase detector, a PID controller, and the equivalent of a voltage-controlled oscillator, as depicted in Figure~\ref{fig:PLLSchematics}. The difference between a reference $\theta_*$ and $\theta$, the phase shift between the response $x$ and either the drive $u$ or the force $f$, is fed to the PID controller, which automatically adapts the frequency of excitation $\omega$ to reach the desired phase. Indeed, if $\theta_*$ is constant, since $\omega$ is the output of the PID,
        \begin{equation}
            \begin{array}{rl}
            \displaystyle \omega (t) = & \displaystyle \omega (0) + k_\mathrm{p}(\theta(t) - \theta_*) \\
            & \displaystyle + k_\mathrm{i}\int_0^t (\theta(\tau) - \theta_*) d\tau + k_\mathrm{d}\dot{\theta}(t),
            \end{array}
            \label{eq:PLLPID}
        \end{equation}
        where $k_\mathrm{p}$, $k_\mathrm{i}$ and $k_\mathrm{d}$ are the proportional, integral and derivative gains of the controller, respectively. Differentiating Equation~\eqref{eq:PLLPID} with respect to time,
        \begin{equation}
            \dot{\omega} (t) = k_\mathrm{p}\dot{\theta}(t) + k_\mathrm{i} (\theta(t) - \theta_*) + k_\mathrm{d}\ddot{\theta}(t),
        \end{equation}
        and assuming steady-state conditions ($\dot{\omega} = \dot{\theta} = 0$) shows that, at steady state, $\theta(t) = \theta_*$ as long as $k_\mathrm{i} \neq 0$. Once $\omega$ reaches a steady value, the control is naturally non-invasive and can stabilize unstable responses. By fixing the forcing amplitude $f_*$ and varying the phase $\theta_*$, one can obtain an NFR. Alternatively, by fixing the phase $\theta_*$ (typically to obtain phase resonance) and varying the force $f_*$, one can obtain a constant-phase curve (typically, a phase-resonance backbone curve).

        \subsubsection{PLL applications}
            
\begin{figure*}[!ht]
                \centering
                \includegraphics[width=0.5\textwidth]{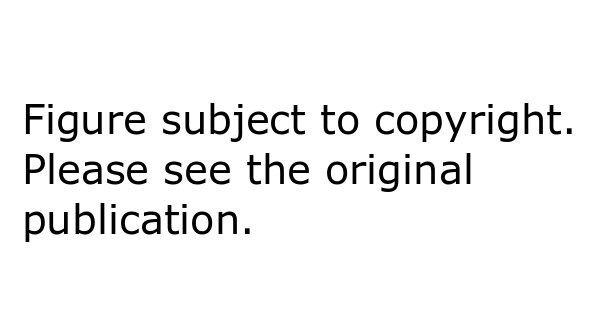}
                \caption{NFRs (colored lines, featuring repeatability) and phase resonance backbone (black dashed line) obtained with a PLL of a turbine blade component with root friction. Reprinted figure with permission from~\cite{Schwarz2020}. Copyright (2020) by the American Society of Mechanical Engineers.}
                \label{fig:Scwharz2020}
            \end{figure*}

            The main thrust for PLL application was (and still is) nonlinear modal analysis, given how easily it can be performed with this control-based method. Early works provided a theoretical study based on normal forms that was developed and validated in~\cite{Denis2018}, and an NFR synthesis method based on the measured phase resonance curve was proposed in~\cite{Peter2018}. Since then, the PLL has been applied to a large variety of structures, including beams with friction~\cite{Scheel2018,Scheel2020b}, impacts~\cite{Peter2019,Woiwode2024}, repelling magnets~\cite{Scheel2020} or featuring extremely large vibration amplitude~\cite{Debeurre2024}, a T-structure~\cite{Nagesh2022}, plates with internal resonances~\cite{Givois2020,Givois2020b} and a turbine blade with friction~\cite{Schwarz2020}. This latter example is an indication of technological maturity, and some of the associated results are plotted in Figure~\ref{fig:Scwharz2020}. In addition, an adaptation for multi-input-multi-output test case was sketched in~\cite{Tang2020} based on a centralized modal-based control scheme, and a framework for phase resonance testing of base-excited systems was proposed in~\cite{Muller2022}.

            Finally, we note that PLL was compared to velocity feedback in~\cite{Scheel2022} and to CBC in~\cite{Abeloos2022}. All methods were shown to be capable of providing relevant results to the experimenter. In-depth dissertations on PLL can be found in~\cite{Peter2018Thesis,Scheel2022Thesis,Abeloos2022Thesis}.
        
    \subsection{Response-controlled stepped-sine testing}

\begin{figure*}[!ht]
            \centering
        	\includegraphics[scale=.6]{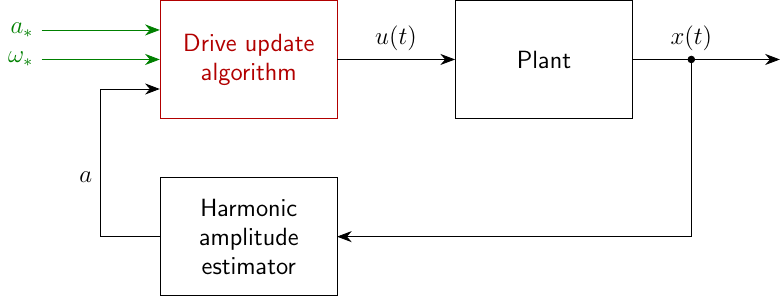}
            \caption{Schematics of the RCT method.}
            \label{fig:RCTSchematics}
        \end{figure*}

        The most recent control-based method is the RCT approach, whose basic principles were first drafted in~\cite{Link2011} and then elaborated in~\cite{Karaagacl2021} together with a theoretical justification. Similarly to the SCBC and PLL, its underlying motivation came from the observation that the frequency-force curves are generally unfolded when the experiment is performed with a constant amplitude response $a_*$ and varied excitation frequency $\omega$. It can thus be seen as a dual of the SCBC method, as illustrated by Table~\ref{tab:continuationClassification}. The general schematics of the RCT method are given in Figure~\ref{fig:RCTSchematics}.

        A key strength of RCT is that, due to numerous cases of modal tests where the structure is base-excited~\cite{Beliveau1986}, conventional linear vibration testing software generally implement a possibility to perform the test at a constant response amplitude. The control principle used in Siemens Testlab (the commercial software used for all RCT demonstrations) is conceptually explained in~\cite{Carrella2020}, although a detailed formulation of the control scheme does not seem to be available. Simply put, it is an automatic trial-and-error procedure that adapts the drive amplitude at a given frequency to match the desired response amplitude. A test at a given target amplitude is started at a given frequency with a trial drive amplitude maintained constant over several periods of excitation, and the harmonic response of the structure is estimated, e.g., with a digital tracking filter. A frequency response function between the response and drive is estimated (or updated) based on the ratio of these amplitudes, and the next drive amplitude is calculated based on this update. The process is repeated until the response amplitude reaches its target (within tolerance bounds) or until a maximum number of corrections have been carried out. The test then carries on to the next frequency. With this method, the control is naturally non-invasive, because the drive is always tonal.

        For a given target amplitude, the test yields a harmonic force curve (HFC). Collecting them for multiple target amplitudes, one obtains a harmonic force surface (HFS), and, similarly to SCBC, NFRs can then be extracted using interpolation of this HFS.
        
        Another interesting feature of the RCT is the individual linear-like look of the frequency-force curves, making it possible to identify them by standard linear system identification tools providing modal characteristics. This quasi-linearity is supported by the single-nonlinear mode theory~\cite{Szemplinska-Stupnicka1979}, which states that if a structure approximately vibrates in a tonal fashion according to one of its mode and if the amplitude of one degree of freedom is maintained constant, the amplitude at the other degrees of freedom will remain constant. This effectively "freezes" the effect of geometrical nonlinearities throughout the structure, giving an intuitive understanding for the observed quasi-linear behavior. This approach also offers an alternative to harmonic force surface interpolation for obtaining NFRs by synthesizing them using identified modal properties.
        
        The RCT method thus provides significant advantages, such as a simplified post-processing procedure leveraging quasi-linearity, and the possibility to use the method directly with commercial vibration measure and control apparatus. When the assumptions for quasi-linearity do not hold, the post-processing might become more complicated but the method still offers a valid control-based framework.

        \subsubsection{RCT applications}
            
            In spite of its relative youth, the RCT has successfully been applied to a variety of examples, including beam setups~\cite{Karaagacl2020,Karaagacl2022,Gurbuz2024}, a bolted missile~\cite{Karaagacl2021}, a stack piezoelectric actuator~\cite{Koyuncu2022} and a highly-damped control fin actuation mechanism~\cite{Karaagacl2024}. We refer to~\cite{Karaagacl2020thesis} for an in-depth dissertation on RCT. The RCT has not been experimentally compared to any other control-based method yet, and such a comparison is carried out in this work.

\section{A derivative-free arclength control-based continuation procedure}
\label{sec:ACBC}

    According to the discussion in Section~\ref{ssec:CBcontinuation}, current control-based methods require a choice from the experimenters. They can use simple, derivative-free methods such as those presented in Table~\ref{tab:continuationClassification}, but the test is then exposed to potential folding issues and the associated possible occurrence of a jump. Some examples will be given later in Section~\ref{ssec:DuffingSuperHNumerical}. Alternatively, they can opt for the more general CBC-FD method, which is however more complex, sensitive to noise and makes the test longer. The main reason for these two downsides is the need for experimental derivatives evaluated through finite differences or Broyden updates.
    
\begin{figure*}[!ht]
        \centering
        \includegraphics[scale=.6]{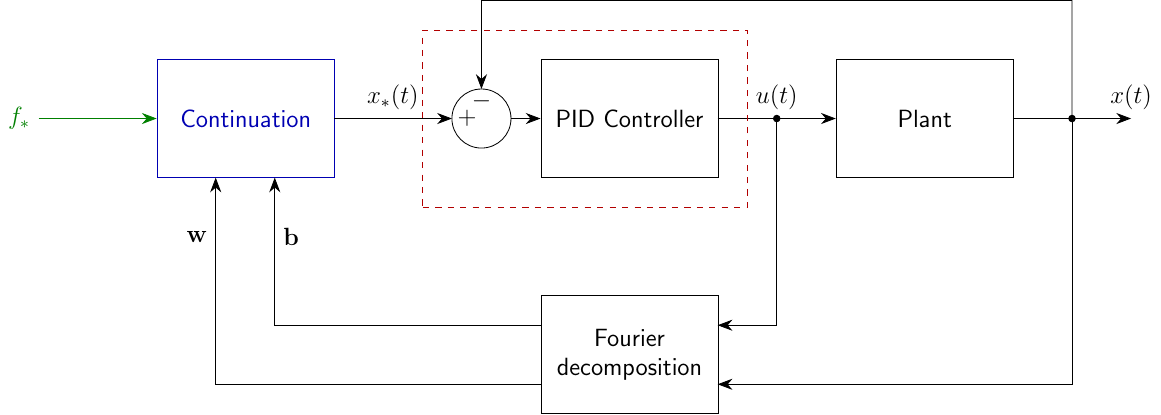}
        \caption{Schematics of the ACBC method; the red dashed frame indicates a controller block.}
        \label{fig:CBCDFSchematics}
    \end{figure*}

    The purpose of this section is to propose a new method which is almost as generic as the CBC-FD method in the sense that it does not assume any parametrization a priori, but which also circumvents the need for experimental derivatives. This procedure is called arclength CBC (ACBC) herein, was originally presented in~\cite{Abeloos2022Thesis}, is schematically depicted in Figure~\ref{fig:CBCDFSchematics} and explained hereafter.

    \subsection{Non-invasive control}

        As in the CBC, a PID-type controller with a reference signal is used to stabilize the structure under test. This reference signal is determined based on two criteria. First, its non-fundamental harmonics are set by the non-invasiveness requirements similarly to the SCBC method. Second, its fundamental harmonic is used together with the excitation frequency to perform continuation.

        Since the method works along the principles of CBC, non-invasiveness is guaranteed as long as the reference and actual signals have identical non-fundamental harmonics. As in~\cite{Abeloos2021}, the non-invasive character of the control signal is ensured by the use of adaptive filters. These filters are able to automatically perform a truncated Fourier decomposition of a signal $x$ with $h$ harmonics, i.e.,
        \begin{equation}
            x(t) \approx \mathbf{h}(t)\mathbf{w},
            \label{eq:FourierDecomposition}
        \end{equation}
        where $\mathbf{w}$ approximate the $2h+1$ Fourier coefficients of $x$.
        
        The  Widrow-Hoff least mean squares (LMS) algorithm can be used to update the Fourier coefficients $\mathbf{w}$ of a signal $x$ through the simple discrete-time law
        \begin{equation}
            \begin{array}{l}
            \displaystyle \mathbf{w}((n+1)t_s) - \mathbf{w}(nt_s) \\
            \displaystyle = \mu t_s (x(nt_s) - \mathbf{h}(nt_s)\mathbf{w}(nt_s))\mathbf{h}^T(nt_s),
            \end{array}
        \end{equation}
        where $\mu$ is the filter gain and $t_s$ is the sampling time.

        Assuming the Fourier coefficients converge, the control can thus be made non-invasive by defining the target signal as
        \begin{equation}
            \begin{array}{rrl}
            \displaystyle x_*(t) &  \displaystyle =  &  \displaystyle a_* \sin(\omega t) + b_* \cos(\omega t)\\
            & & \displaystyle + \mathbf{h}(t)\left(\mathbf{I} - \mathbf{e}_2\mathbf{e}_2^T - \mathbf{e}_{h+2}\mathbf{e}_{h+2}^T\right)\mathbf{w} \\
            &  \displaystyle = & \displaystyle a_* \sin(\omega t) + b_* \cos(\omega t) + \mathbf{h}(t)\mathbf{P}_{\mathrm{nf}}\mathbf{w}.
            \end{array}
            \label{eq:CBCDFXstar}
        \end{equation}
        $\mathbf{P}_{\mathrm{nf}}$ can be interpreted as a matrix whose action retains only the non-fundamental harmonics so that, if Equations~\eqref{eq:FourierDecomposition} and~\eqref{eq:CBCDFXstar} hold, the signal $x_* - x$ is mono-harmonic.
        
    \subsection{Simplified continuation problem}

        Having secured the non-invasiveness of the method, three parameters remain to be set by the method, namely $a_*$, $b_*$ and $\omega$. Without loss of generality, one of the Fourier coefficients can be set to zero, say $b_* = 0$. The two remaining unknowns are set by enforcing the fundamental forcing amplitude $f$ to be equal to the target one $f_*$. 
        
\begin{figure}[!ht]
            \centering
            \begin{subfigure}{.49\textwidth}
                \centering
                \includegraphics[width=1.1\textwidth]{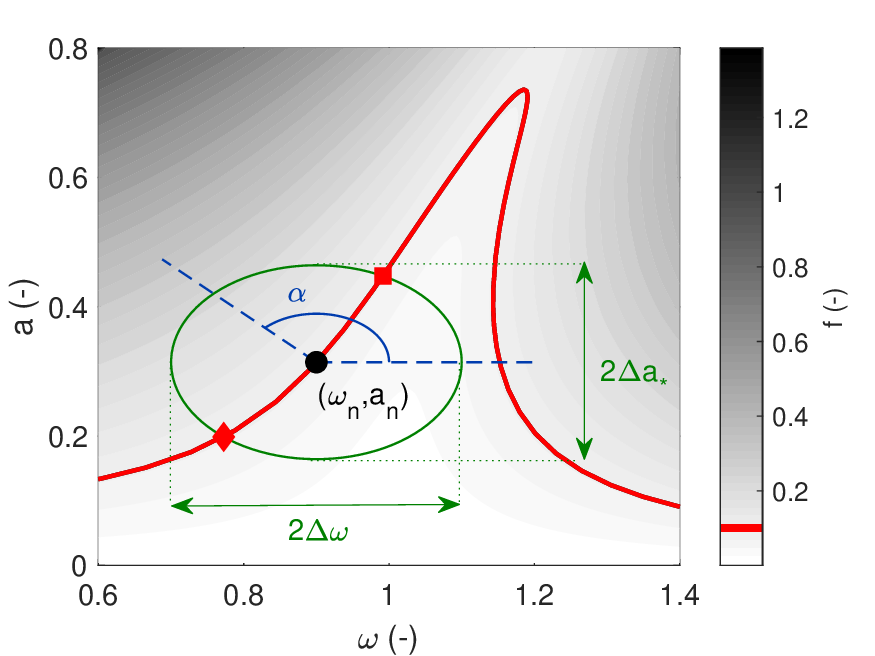}
                \caption{}
                \label{sfig:Ellipse}
            \end{subfigure}
            \begin{subfigure}{.49\textwidth}
                \centering
                \includegraphics[width=1.1\textwidth]{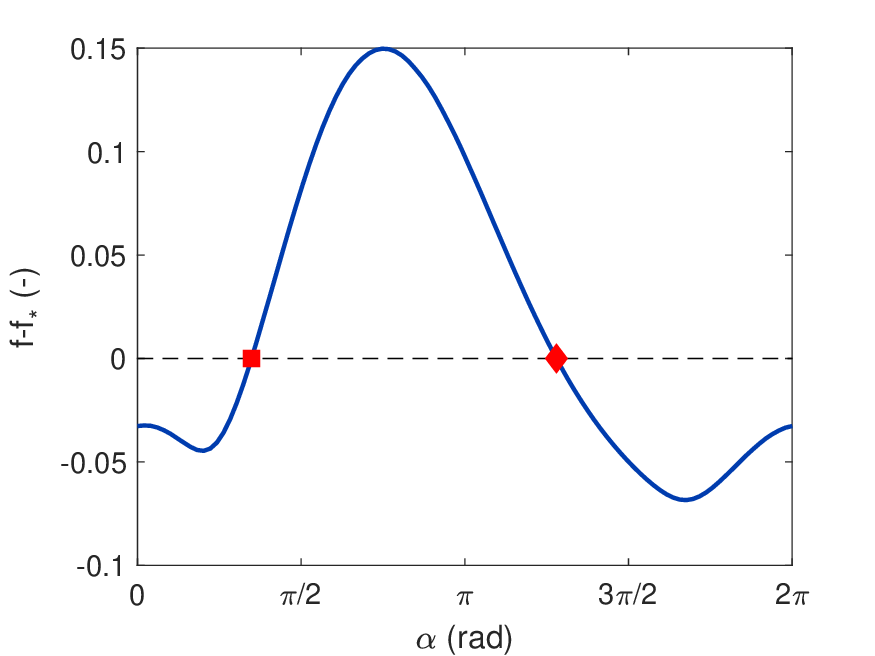}
                \caption{}
                \label{sfig:ForceEllipse}
            \end{subfigure}
            \caption{\subref{sfig:Ellipse}: HFS (levels of gray) and definition of the ellipse (\textcolor{MatlabDarkGreen}{\rule[.2em]{1em}{.2em}}) built around a given point ($\bullet$) used for the ACBC method to obtain points (\textcolor{red}{$\blacklozenge$}, \textcolor{red}{$\blacksquare$}) on the NFR (\textcolor{red}{\rule[.2em]{1em}{.2em}}) at $\bar{f}=\bar{f}_*=0.0875$. \subref{sfig:ForceEllipse}: Profile of the function $\bar{f}-\bar{f}_*$ along the ellipse, whose roots (\textcolor{red}{$\blacklozenge$}, \textcolor{red}{$\blacksquare$}) correspond to the intersection of the ellipse with the NFR.}
            \label{fig:CBCDF}
        \end{figure}

        The method can probably be best understood by looking at the force imposed on the structure as a function of the frequency and target amplitude. On this representation of the HFS depicted in Figure~\ref{sfig:Ellipse}, the NFR can be seen as a level set of this surface, where the forcing amplitude $f$ corresponds to the target $f_*$. In this example, we consider for simplicity infinitely large control gains, such that $a = a_*$. However, we note that the foregoing discussion also holds for finite gains provided they are large enough.

        The new approach for experimental continuation leverages this, and is inspired from arclength continuation~\cite{Riks1979}. Because the NFR is one-dimensional in the $(\omega,a_*)$ space, a sufficiently small ellipse with semi-axes $\Delta \omega$ and $\Delta a_*$ and centered on a point of the branch intersects it twice, as depicted in Figure~\ref{sfig:Ellipse}. The intersections are characterized by $f-f_* = 0$, whereas every other point of the ellipse does not (in general) satisfy this equality. The continuation procedure can thus be formulated as follows: given a current point $(\omega_n,a_n)$ and an ellipse centered around it in the $(\omega,a_*)$ space, the next point of the NFR $(\omega_{n+1},a_{n+1})$ is found as the intersection of the ellipse and the level set where $f = f_*$, as seen in Figure~\ref{sfig:ForceEllipse}. To do so, the procedure sweeps along the ellipse by automatically varying the eccentric anomaly $\alpha$ of a point on the ellipse. The procedure can then be repeated using ellipses centered around each next point, thereby obtaining a set of points on the NFR in the same fashion as Crisfield's arclength continuation~\cite{CRISFIELD1981}.

        A key advantage of this approach is that the harmonic forcing amplitude can easily be measured online, and there are several strategies to find the eccentric anomaly that solves the zero problem $f-f_* = 0$ without need for experimental derivatives: 
        \begin{enumerate}
            \item A first simple possibility consists in sweeping the angle $\alpha (t)$ at a constant rate $r_{\alpha}$ until $f=f_*$. 
            \item A continuous integral control law $\dot{\alpha} = -k_{\alpha} (f-f_*)$ can be used. An advantage of this law is that it automatically makes the algorithm converge toward a point where $f=f_*$, and orients the curve so that the method does not turn back. Indeed, this law makes the point where $f=f_*$ and $\partial f/\partial \alpha >0$ a stable equilibrium, whereas the other solution is unstable. As a consequence, this law converges to intersections which pierce the NFR outward as $\alpha$ is increased. However, since the sweep rate depends on the error $f-f_*$, the sweeping rate of the ellipse depends on the amplitude of the profile $f-f_*$, and is typically low and high away and close to resonances, respectively.
            \item A discontinuous integral control law $\dot{\alpha}(t) = -r_{\alpha} \text{sign} (f-f_*)$ combines the advantages of the sweeping and integral control methods. However, the discontinuity of this control law can make the convergence toward a solution where $f-f_*=0$ difficult.
        \end{enumerate}
        The sweeping law can be chosen freely by the user depending on their preference, and is denoted by $\dot{\alpha} = s(f-f_*)$ in the sequel. The sweeping rate should be high for a fast continuation, but is nevertheless limited by the system's transients, as well as the adaptive filters speed to guarantee a non-invasive control.
        
    \subsection{The algorithm}

        Algorithm~\ref{alg:ACBC} summarizes the proposed derivative-free ACBC procedure. In this implementation, the prediction is made with the previous converged value of $\alpha$, which amounts to implementing a secant predictor. $\rho$ represents a relative tolerance, and $\sigma < 1$ is a factor that makes the algorithm stop sweeping in more strict conditions than the convergence check in order to anticipate transient effects. If this preemptive procedure is not enough to meet the convergence conditions, one can adapt the sweeping law by making it slower.

        One can also decide to measure the NFR in one direction or another (with globally increasing or globally decreasing forcing frequencies) by selecting a positive or negative value for $k_\alpha$ or $r_\alpha$.

        \begin{algorithm} 
        \caption{Derivative-free arclength continuation for CBC experiment.}
        \label{alg:ACBC} 
        \algrenewcommand\algorithmicindent{0.7em}%
        \begin{algorithmic}[1]
            \State \textbf{Input:} $\omega_0,a_0,\Delta \omega, \Delta a_*$ selected by the user.
            \State $n\gets 0$            
            \If{$k_\alpha > 0$ \textbf{or} $r_\alpha > 0$}
                \State $\alpha \gets 0$
            \Else
                \State $\alpha \gets \pi$
            \EndIf
            \Loop
                \State \underline{Prediction}
                \State $(\omega,a_*) \gets (\omega_{n} + \Delta \omega \cos (\alpha) , a_n + \Delta a_* \sin (\alpha) )$
                \State Wait $t_{\mathrm{steady}}$ seconds for steady-state
                \State \underline{Correction} 
                \While{$|f-f_*| > \sigma \rho f_*$}
                    \State $\dot{\alpha}=s(f-f_*)$
                    \State $(\omega,a_*)  \gets (\omega_{n} + \Delta \omega \cos (\alpha) , a_n + \Delta a_* \sin (\alpha) )$
                \EndWhile
                \State Wait $t_{\mathrm{steady}}$ seconds for steady-state
                \State \underline{Convergence check} 
                \If{$|f-f_*| \leq \rho f_*$}
                    \State $(\omega_{n+1},a_{n+1})\gets (\omega,a_*)$
                    \State $n \gets n+1$
                \Else
                    \State Adapt $s$ (optional)
                \EndIf
            \EndLoop
        \end{algorithmic}
        \end{algorithm}

        The proposed algorithm is equivalent to the CBC-FD in terms of generality of the continuation approach, as it does not make any assumption on the unfolding property of a parameter. A key advantage of this approach is its derivative-free nature, brought by the combined use of the adaptive filters and the arclength continuation procedure in a two-dimensional space.

    \subsection{Illustration with a Helmholtz-Duffing oscillator}

        \begin{figure}[!ht]
            \centering
                \begin{subfigure}{.49\textwidth}
                    \centering
                    \includegraphics[width=1.1\textwidth]{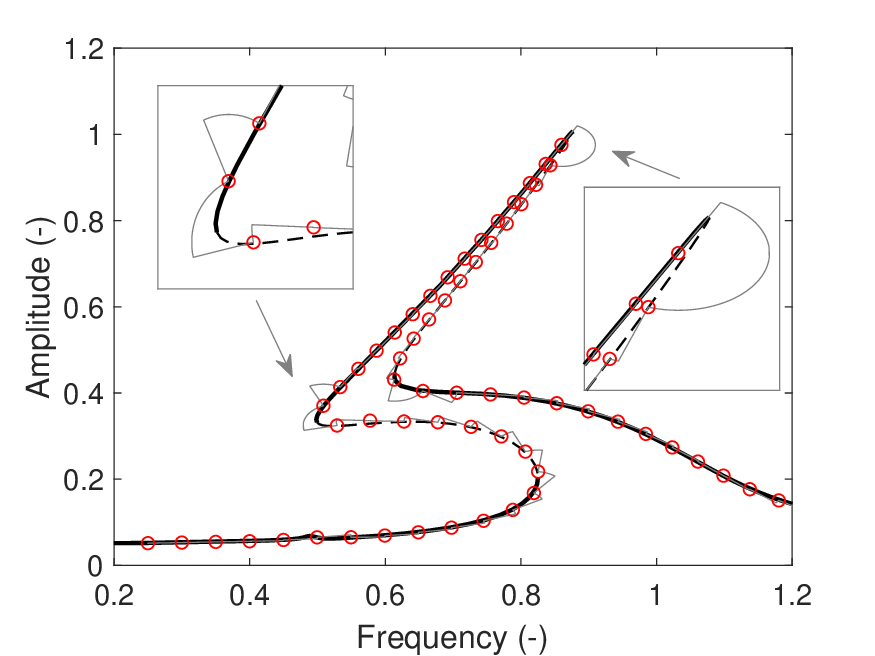}
                    \caption{}
                    \label{sfig:HD_Simulink_1}
                \end{subfigure}
                \begin{subfigure}{.49\textwidth}
                    \centering
                    \includegraphics[width=1.1\textwidth]{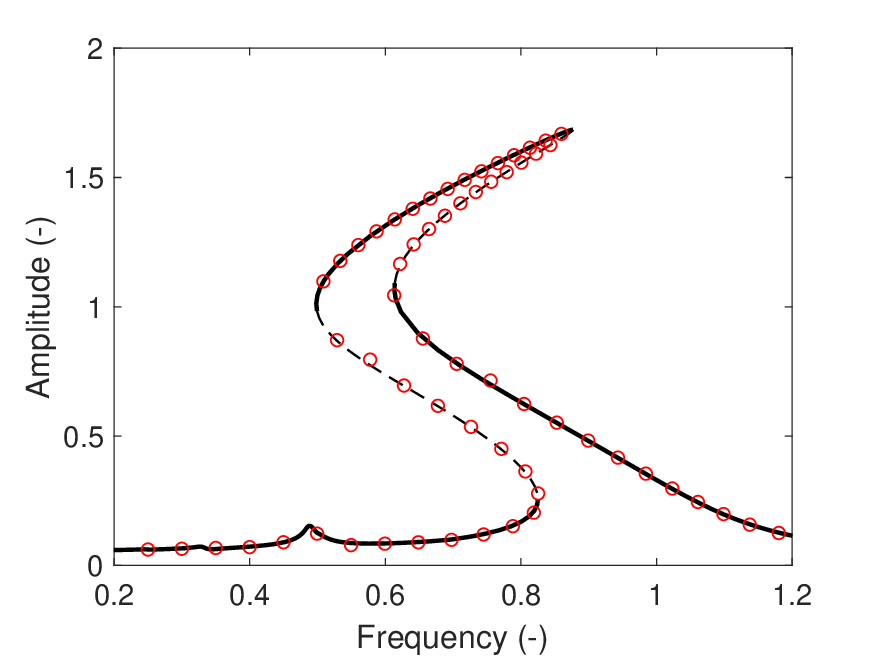}
                    \caption{}
                    \label{sfig:HD_Simulink_2}
                \end{subfigure}
            \caption{Bifurcation diagram of the Helmholtz-Duffing oscillator with $\zeta_0 = 0.05$, $\beta_2 = 1.9$ and a dimensionless forcing amplitude of 0.05 obtained with the ACBC: input parameter space~\subref{sfig:HD_Simulink_1} and NFR~\subref{sfig:HD_Simulink_2}. The full $(\omega,a_*)$ curve is shown (\textcolor{Gray1}{\rule[.2em]{1em}{.2em}}) together with the results at convergence (\textcolor{red}{$\circ$}). The results from numerical continuation are also given for reference (\rule[.2em]{1em}{.2em}), where full and dashed lines represent open-loop stable and unstable solutions, respectively.}
            \label{fig:HelmholtzDuffingACBC}
        \end{figure}

        The Helmholtz-Duffing example in Figure~\ref{fig:HelmholtzDuffingSweep} is revisited with a numerical simulation of the ACBC algorithm in Simulink with $\bar{k}_\mathrm{d}=1$, $k_\mathrm{\alpha}=1$, $h=15$, $\mu=0.0025$, $\Delta \bar{a}_* = \Delta \bar{\omega} = 0.05$ and $\bar{t}_{\mathrm{steady}}=50$. Clearly, the ACBC algorithm is able to capture the full NFR in this example, including its open-loop unstable parts.

\section{Illustration with an electronic Duffing oscillator}
\label{sec:ElectronicDuffing}

    Having reviewed the different control-based approaches for experimental continuation, we illustrate them with two experimental examples.

    A first demonstration is made with an analog electronic circuit that replicates the behavior of a Duffing oscillator~\cite{Jones2001,Srinivasan2009}. This circuit is described in more details in the supplementary materials of this article and in~\cite{Raze2024}. The main advantages of using such an experimental demonstrator are its near perfect repeatability and actuation (i.e., the drive and the force are almost exactly equal), which allows us to test the inherent features of control-based methods. Yet, this set-up is subject to practical experimental challenges, such as noise, unknown non-ideal behavior, and real-time constraints.

    \subsection{Experimental setup}

\begin{figure}[!ht]
        \centering
        \begin{subfigure}{.49\textwidth}
            \centering
        	\includegraphics[scale=.8]{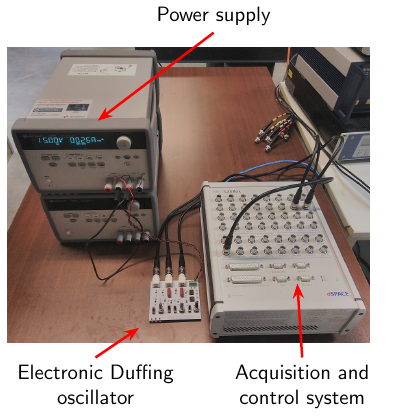}
            \caption{}
            \label{sfig:DuffingSetupPicture}
        \end{subfigure}
        \begin{subfigure}{.49\textwidth}
            \centering
        	\includegraphics[scale=.4]{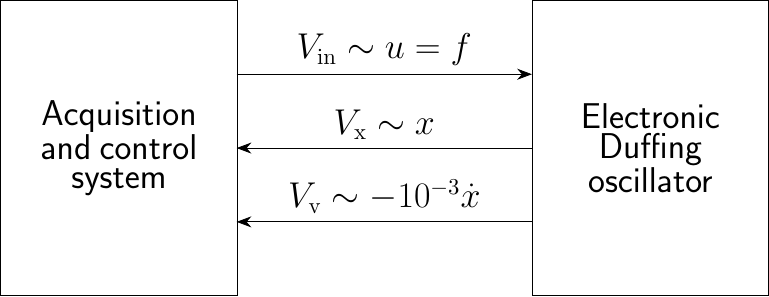}
            \caption{}
            \label{sfig:DuffingSetupSchematics}
        \end{subfigure}
        \caption{Picture~\subref{sfig:DuffingSetupPicture} and schematics~\subref{sfig:DuffingSetupSchematics} of the experimental setup with an electronic Duffing oscillator.}
        \label{fig:DuffingSetupSchematics}
    \end{figure}

    A picture of the experimental apparatus, comprising the electronic Duffing oscillator and an acquisition and control system, is given in Figure~\ref{sfig:DuffingSetupPicture} and schematized in Figure~\ref{sfig:DuffingSetupSchematics}. The oscillator takes a forcing signal $V_\mathrm{in}$ as input, and outputs a displacement signal $V_\mathrm{x}$ and a velocity-proportional signal $V_\mathrm{v}$, with $V_\mathrm{v} \approx -10^{-3} \dot{V}_\mathrm{x}$. Compared to Figure~\ref{fig:controlBasedSchematics}, the acquisition and control system encompasses the continuation, controller and signal analyzer blocks. For all methods except the RCT, this system was a MicroLabBox from dSPACE (as in Figure~\ref{sfig:DuffingSetupPicture}), whereas the RCT was performed with a Scadas Mobile from Siemens.

    The parameters $c$, $k$ and $k_3$ (cf. Equation~\eqref{eq:DuffingODE}) of the electronic Duffing oscillator can be tuned with knobs (potentiometers), and their theoretical values can be computed based on the electrical parameters of the board (see supplementary materials). They are reported in Table~\ref{tab:DuffingParameters}. They were compared to an identification with the nonlinear frequency subspace identification (FNSI) method~\cite{Noel2013} using data from a swept-up sine test. A close agreement between the theoretical and identified values can be observed, with a discrepancy of up to 8\% for the damping coefficient. These parameters correspond to a natural frequency $\omega_0=131.5$ rad/s (20.9 Hz) and a damping ratio $\zeta_0=0.005$. This damping ratio can be considered rather small although it is still representative for a number of metallic structures. It was initially chosen to result in large amplitudes for moderate forcing, and strong nonlinearity, but it also revealed some limitations of various methods that seemed interesting to report herein.

    \begin{table}[!ht]
        \centering
        \caption{Theoretical and identified parameters of the electronic Duffing oscillator.}
        \label{tab:DuffingParameters}
        \begin{tabular}{ccc}
             \hline
             Parameter & Theoretical value & Identified value \\
             \hline
             $m$ ($s^2$) & $10^{-4}$ & $1.0461\times 10^{-4}$\\
             $c$ ($s$) & 1.3$\times 10^{-4}$ & 1.4231$\times 10^{-4}$\\
             $k$ (-) & 1.923 & 1.8099\\
             $k_3$ ($V^{-2}$) & 0.9887 & 1.0077\\
             \hline
        \end{tabular}
    \end{table}

    \subsection{Implementation details}
    
    All methods (except for the RCT) were programmed in a fully online way with dSPACE RTI for Simulink. Almost all methods require a Fourier decomposition. This decomposition was performed online by adaptive filters, which were shown to be quick, effective and accurate~\cite{Abeloos2021,Abeloos2022Thesis,Woiwode2024}. For a tonal signal at frequency $\omega$, an expression for an optimal gain (minimizing the settling time of the Fourier decomposition) was derived in~\cite{Abeloos2022Thesis}\footnote{Note that the symbol $\mu$ used herein corresponds to $\mu t_s$ in~\cite{Abeloos2022Thesis}.}. Here, based on a linearization of this expression with respect to the (small) sampling time, a frequency-dependent gain was used
    \begin{equation}
        \mu = \omega \bar{\mu},
    \end{equation}
    where the optimal dimensionless gain based on this linearization is $\bar{\mu}=2$. For a signal with multiple harmonics, this value was observed to lead to a slightly too aggressive filter, and a smaller value was usually selected. Another procedure to tune $\mu$ was recently proposed in~\cite{Hippold2024}.
        
    Some iterative methods propose specific procedures to detect if the system is at steady-state. In this work, we simply waited a certain time $t_{\mathrm{steady}}$ before subjecting the system to another set of parameters. 

    All CBC methods were implemented with the same differential controller. It should be noted that the velocity signal was used instead of the displacement one, thereby avoiding the use of an approximate derivation with a high-pass filter (as can commonly be done in practical PID implementations).

    Stepped methods (CBC-FD, SCBC, ACBC and RCT) require the definition of steps $\Delta \omega$ and $\Delta a_*$, which are based on the desired resolution of the NFR. For CBC methods, since the feedback is made with a velocity-type signal, the steps were defined in terms of this signal, $\Delta |V_\mathrm{v}|_1$ for simplicity, where $|\cdot|_k$ denotes the $k^{th}$ harmonic amplitude of a signal. For the RCT, a (nearly) equivalent step in terms of displacement $\Delta |V_\mathrm{x}|_1$ was defined instead. The two are related by $\Delta |V_\mathrm{v}|_1 = 10^{-3}\omega \Delta |V_\mathrm{x}|_1$.

    Most of the parameters were tuned on a trial and error basis. A fine-tuning toward optimal performance was not sought given the very large number of parameters and their variety. The results given therein are thus not to be interpreted as the best each method can yield, but rather as an order of idea of how they compare.

    In the sequel, most of the results will be presented in terms of the total amplitude of the response over a period, being represented by the operator $|\cdot|$ (without subscript), i.e.,
    \begin{equation}
        |x| = \max_{t\in [0,2\pi/\omega]} \text{abs}(x(t)),
    \end{equation}
    where $\text{abs}$ is the absolute value operator.

    \subsection{Primary resonance}

        The primary resonance of the oscillator was first tested. The time for the two methods where parameters are swept, namely SWS and PLL, were chosen based on the fastest stepped method (in this case, ACBC). All parameters are given in the supplementary materials. We focus on measuring an NFR at a specific excitation $f=0.03$ V between 50 and 200 rad/s.
        
        \subsubsection{Comparison of the NFR obtained with different measurement approaches}
        \label{ssec:DuffingNFRComparison}

\begin{figure}[!ht]
                \centering
                \begin{subfigure}{.49\textwidth}
                    \centering
                    \includegraphics[width=1.1\textwidth]{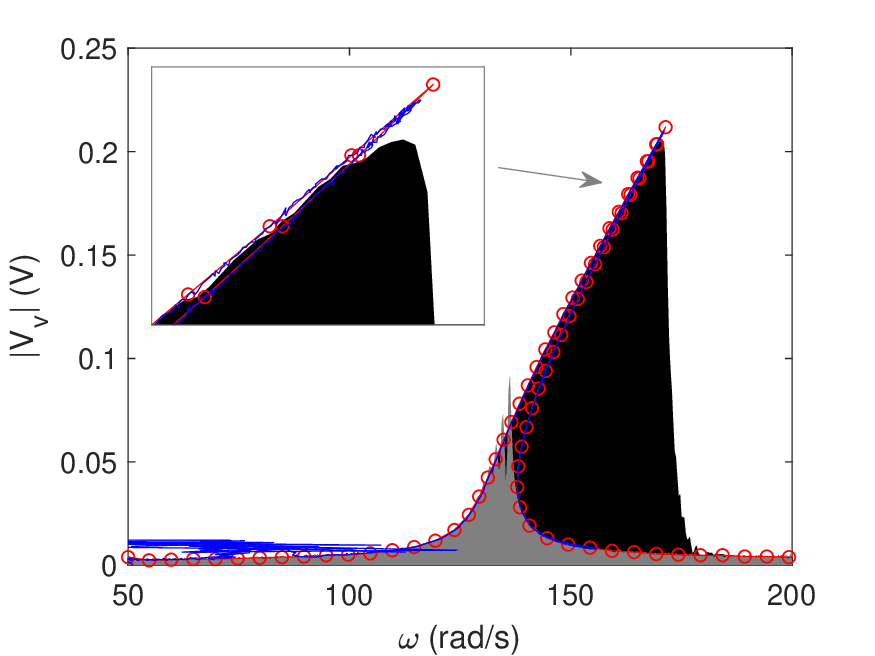}
                    \caption{}
                    \label{sfig:NFR_1}
                \end{subfigure}
                \begin{subfigure}{.49\textwidth}
                    \centering
                    \includegraphics[width=1.1\textwidth]{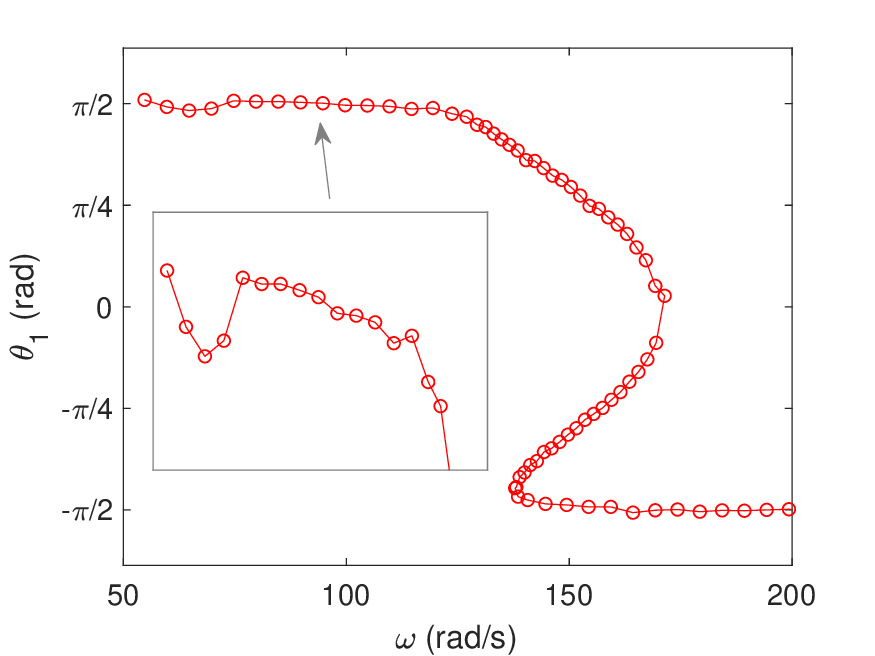}
                    \caption{}
                    \label{sfig:NFR_1_Phase}
                \end{subfigure}
                \caption{Velocity NFR of the electronic Duffing oscillator obtained with the SWS-up (\rule[.2em]{1em}{.2em}) SWS-down (\textcolor{gray}{\rule[.2em]{1em}{.2em}}), PLL (\textcolor{blue}{\rule[.2em]{1em}{.2em}}) and ACBC (\textcolor{red}{\rule[.2em]{1em}{.2em}}) methods (\subref{sfig:NFR_1}: amplitude, \subref{sfig:NFR_1_Phase}: phase).}
                \label{fig:NFR_1}
            \end{figure}

             To start this analysis off, Figure~\ref{fig:NFR_1} shows the velocity signal NFR directly obtained from methods that ran without issues, namely, the SWS, PLL and ACBC approaches. The agreement between the two latter is excellent, even at the resonance peak. Some discrepancies are observable with the SWS data: as usual, transient effects are the cause for a slightly lower peak amplitude. We also note that the SWS data includes higher harmonics as well. The PLL test experiences large noise at low frequency, which is due to the 2:1 superharmonic resonance of the oscillator. While it is too small to be noticeable in the amplitude of the NFR, Figure~\ref{sfig:NFR_1_Phase} shows that the phase of the response evolves in a non-monotonous way at low frequencies, due to this superharmonic resonance. Once the phase becomes unique again (in particular near the primary resonance), this is no longer an issue.
            
            
        \subsubsection{Challenges for the pseudo-arclength CBC method}
        \label{ssec:CBCFDIssues}

\begin{figure}[!ht]
                \centering
                \includegraphics[width=0.539\textwidth]{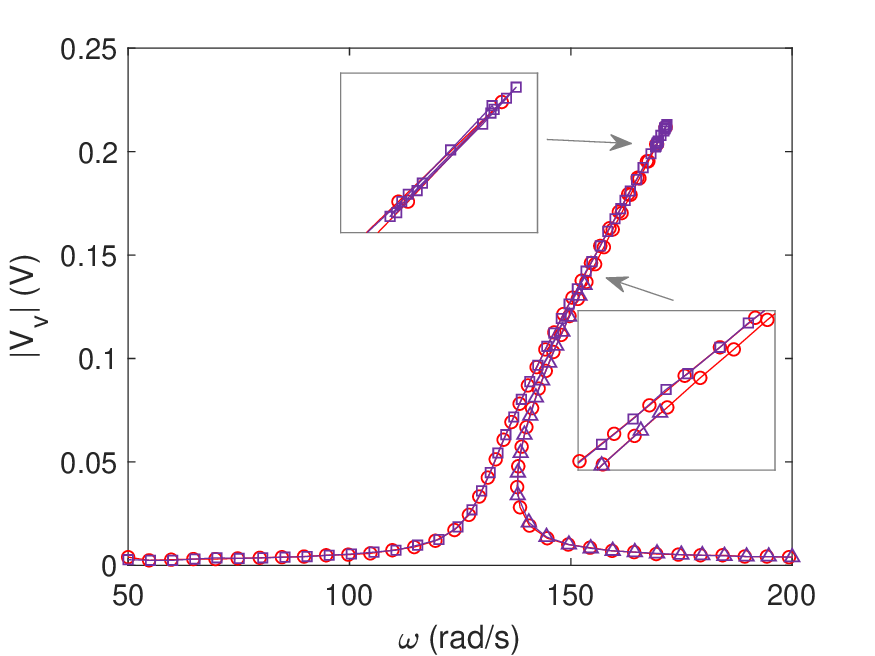}
                \caption{Velocity NFR of the electronic Duffing oscillator obtained with the ACBC (\textcolor{red}{\rule[.2em]{1em}{.2em}}) and CBC-FD (\textcolor{MatlabPurple}{\rule[.2em]{1em}{.2em}}, $\square$: starting from 50 rad/s, $\triangle$: starting from 200 rad/s) methods.}
                \label{fig:NFR_2}
            \end{figure}

            The results obtained with the CBC-FD method are compared to those of the ACBC method in Figure~\ref{fig:NFR_2}. Convergence issues were encountered with the former method in the vicinity of the resonance peak, requiring to restart the method from the high-frequency limit (200 rad/s) to approach the peak from both sides.
            
            The difficulties encountered with CBC-FD are most likely due to two main causes. First, owing to the small damping of the oscillator, the NFR curve experiences a sharp turn near the fold bifurcation. As a consequence, the pseudo-arclength method can struggle to converge to the curve during a correction phase. This issue could be solved by using an arclength continuation scheme, or taking smaller steps. The second issue is linked to the use of finite differences, and again, is mostly problematic in the vicinity of the resonance peak. Indeed, close to this peak, the residual function varies extremely rapidly and has a large curvature (as will also be illustrated by Figure~\ref{sfig:Interpolation_1}). Similar issues were encountered in~\cite{Barton2011,Shen2021,Beregi2023p}. Smaller finite difference steps could have been taken, but this would shift the issues to zones away from resonance, where the residual function is much less sensitive to variations in the input, and where large enough finite difference steps need to be taken to overcome noise. This issue could be solved using an adaptive finite difference step by, e.g., prescribing an ideal change in the residual function (suitably larger than the noise level).

        \subsubsection{Resolution requirements for the methods using interpolation}

\begin{figure}[!ht]
                \centering
                \begin{subfigure}{.49\textwidth}
                    \centering
                    \includegraphics[width=1.1\textwidth]{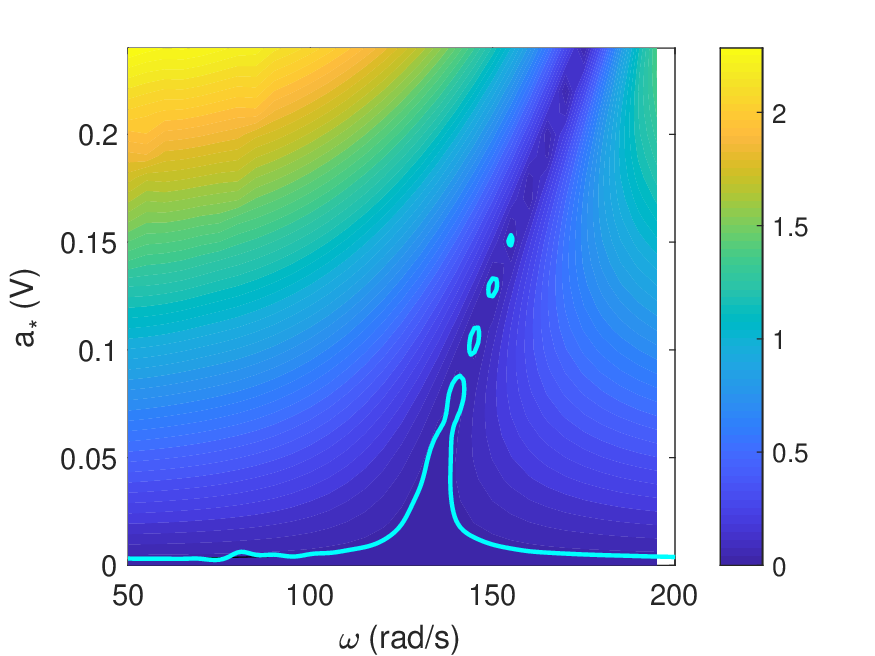}
                    \caption{}
                    \label{sfig:Interpolation_1}
                \end{subfigure}
                \begin{subfigure}{.49\textwidth}
                    \centering
                    \includegraphics[width=1.1\textwidth]{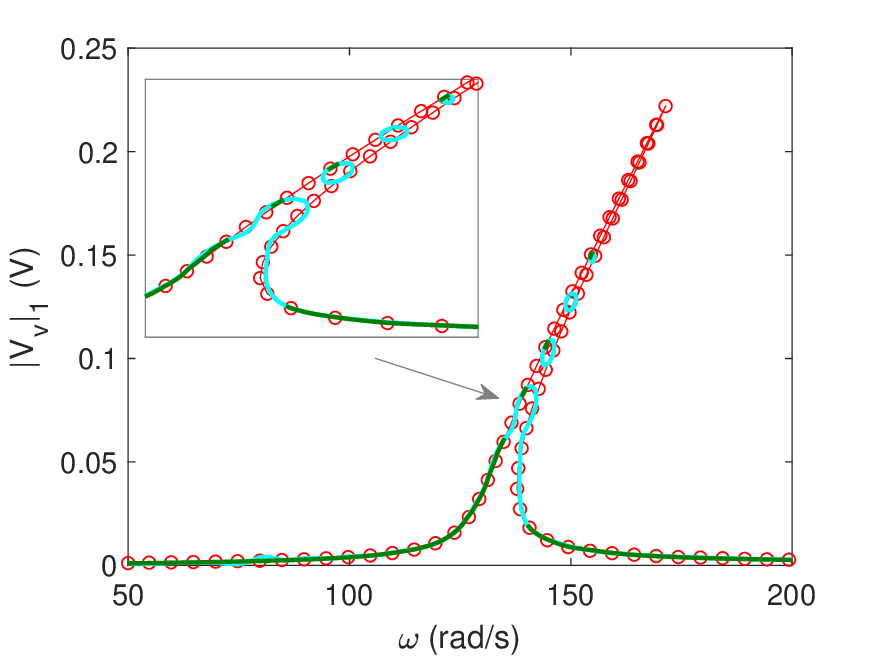}
                    \caption{}
                    \label{sfig:Interpolation_2}
                \end{subfigure}
                \caption{Contour plot of the measured force for the SCBC test with the contour $\hat{f}=0.03$ V highlighted (\textcolor{MatlabCyan}{\rule[.2em]{1em}{.2em}})~\subref{sfig:Interpolation_1}, and velocity NFR of the electronic Duffing oscillator obtained with the ACBC method (\textcolor{red}{\rule[.2em]{1em}{.2em}}) and estimated with the SCBC (\textcolor{MatlabCyan}{\rule[.2em]{1em}{.2em}}) and RCT (\textcolor{MatlabDarkGreen}{\rule[.2em]{1em}{.2em}}) methods~\subref{sfig:Interpolation_2}.}
                \label{fig:Interpolation}
            \end{figure}

            As explained in the respective works where the methods were proposed, both SCBC and RCT obtain a collection of response curves that can be combined into a response surface. NFR curves can be extracted from this surface by slicing it with a plane characterized by $|f|=f_*$, which essentially amounts to performing interpolation on experimental data. An off-the-shelf solution used in this work was to use the \texttt{griddedInterpolant} function in Matlab to fit the measured force and amplitude as 2D spline functions of the frequency and target amplitude, as $\hat{f}(\omega,a_*)$ and $\hat{a}(\omega,a_*)$, respectively (where a hat denotes an interpolated function). The curves satisfying $\hat{f}(\omega,a_*)-f_* = 0$ were then found using the \texttt{fimplicit} function. Figure~\ref{sfig:Interpolation_1} shows such curves obtained for the SCBC. Points of these curves were finally fed to the fitted function $\hat{a}(\omega,a_*)$ to obtain the NFR curve.

\begin{figure}[!ht]
                \centering
                \includegraphics[width=0.539\textwidth]{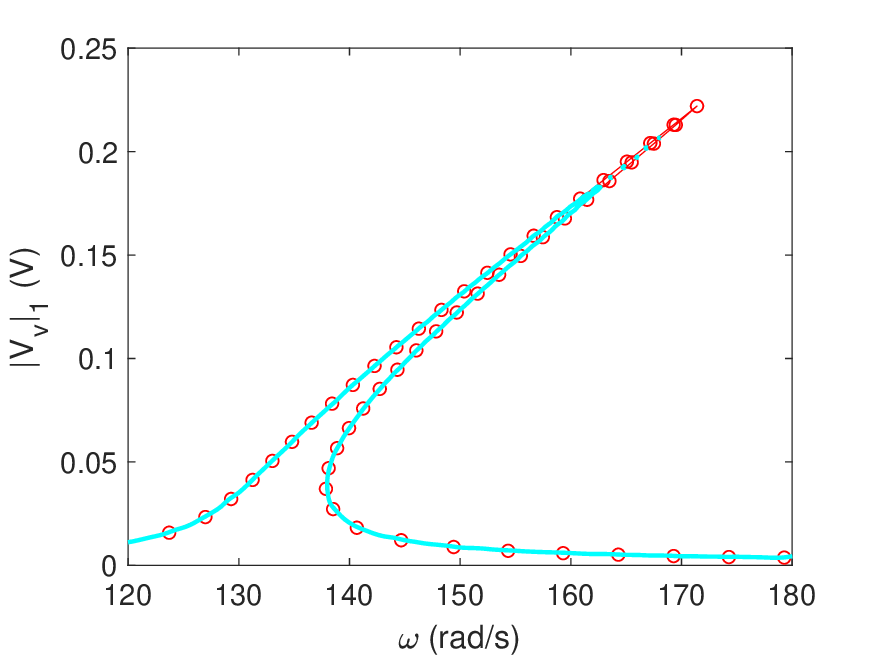}
                \caption{Velocity NFR of the electronic Duffing oscillator obtained with the ACBC method (\textcolor{red}{\rule[.2em]{1em}{.2em}}) and estimated with the SCBC method with higher frequency and amplitude resolution (\textcolor{MatlabCyan}{\rule[.2em]{1em}{.2em}}).}
                \label{fig:Interpolation_3}
            \end{figure}

            The results of the SCBC and RCT obtained with this method are compared to those of the ACBC in Figure~\ref{sfig:Interpolation_2}. In this case, we analyze the fundamental harmonic amplitude rather than the total (multi-harmonic) one, because the latter is more complicated to extract with the results from TestLab with the RCT. These two amplitudes do not differ much anyway. While the interpolated NFRs are close to the one from ACBC away from the resonance peak, there is a clear discrepancy between the two types of results in the vicinity of the peak. This issue is mainly due to the interpolation, and to an insufficient resolution of the ($\omega, a_*$) test grid: near resonance, the function $\hat{f}(\omega,a_*)$ varies rapidly and requires a very fine resolution for the interpolation to be correct. Similar issues were encountered in~\cite{Hayashi2024} for the detection of fold bifurcations from the response surface. This can be solved using either a finer grid (but this requires a longer test time), or using an adaptive sampling strategy instead of a regular grid. Figure~\ref{fig:Interpolation_3} supports this hypothesis by showing the results obtained by narrowing the test frequencies from 120 to 180 rad/s with steps of 2 rad/s, and amplitude steps of 0.05 V, effectively doubling the resolution in both directions. Clearly, the agreement between both methods is improved, at the expense of a longer test time (increasing by a factor 4). The interpolation with the SCBC still inadequately represents a range of about 10 rad/s near the resonance peak, though.
            
            In addition to this issue, we noted that the SCBC method could encounter convergence struggles, as will be discussed in Section~\ref{sec:CBCAIssues}. As for the RCT, another issue related to stability was encountered, which is why the curve in Figure~\ref{sfig:Interpolation_2} appears incomplete. This issue is discussed in Section~\ref{sec:RCTIssues}.

        \subsubsection{Invasiveness of the SCBC in the vicinity of superharmonic resonances}
        \label{sec:CBCAIssues}
        
            The SCBC method was implemented with both the original scheme based on Picard iterations~\cite{Barton2013}, as well as a scheme that exploits the output of an adaptive filter~\cite{Abeloos2021} to make the control non-invasive. The former was observed to be much slower than the latter at some measurement points, forcing a raise of the tolerance to 5\% for the Picard scheme (instead of 1\% for the adaptive scheme) to allow for convergence. As shall be shown hereafter, this occurs in the vicinity of strong superharmonic resonances. Away from them, both methods performed similarly. The difference between the two methods in terms of results presented in Figure~\ref{sfig:Interpolation_2} were negligible, which is why only the results of the method based on adaptive filters were shown.

\begin{figure}[!ht]
                \centering
                \begin{subfigure}{.49\textwidth}
                    \centering
                    \includegraphics[width=1.1\textwidth]{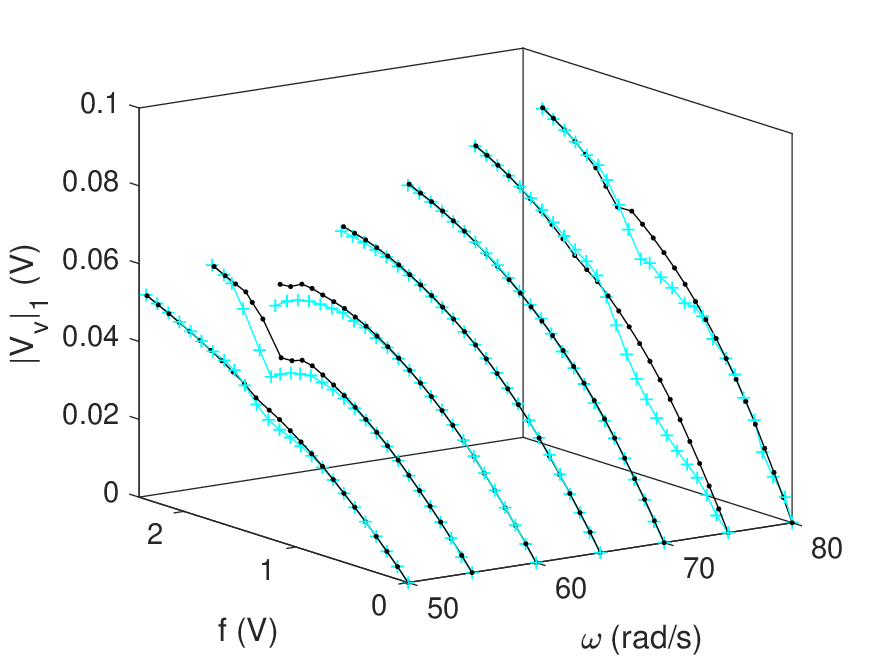}
                    \caption{}
                    \label{sfig:CBCA_1}
                \end{subfigure}
                \begin{subfigure}{.49\textwidth}
                    \centering
                    \includegraphics[width=1.1\textwidth]{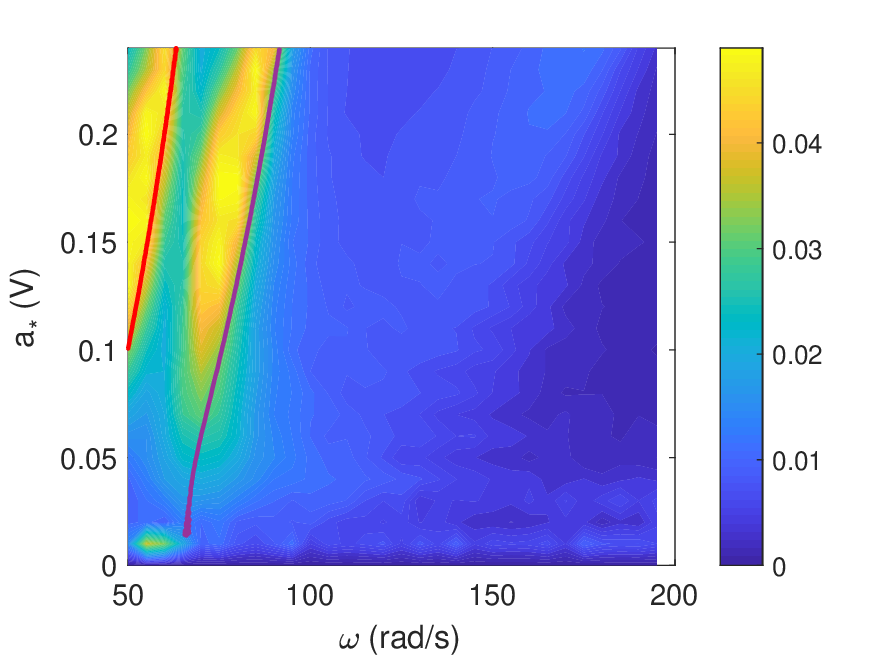}
                    \caption{}
                    \label{sfig:CBCA_2}
                \end{subfigure}
                \caption{Comparison of the response surfaces obtained with the SCBC method with a Picard scheme (\rule[.2em]{1em}{.2em}) and adaptive filters (\textcolor{MatlabCyan}{\rule[.2em]{1em}{.2em}})~\subref{sfig:CBCA_1}, and contour plot of the invasiveness of the SCBC method with a Picard scheme with the locus of 2:1 (\textcolor{MatlabPurple2}{\rule[.2em]{1em}{.2em}}) and 3:1 (\textcolor{red}{\rule[.2em]{1em}{.2em}}) superharmonic resonances~\subref{sfig:CBCA_2}.}
                \label{fig:CBCA}
            \end{figure}

            Figure~\ref{sfig:CBCA_1} compares the results of the two methods in terms of the response surface. Differences are observable at low frequencies, while the agreement is almost perfect for the rest of the surface (not shown here). Figure~\ref{sfig:CBCA_2} shows the relative error, given by the invasiveness (the norm of the left-hand side of Equation~\eqref{eq:CBCANoninvasiveness} divided by $a_*$) made by the Picard scheme, featuring the largest error in the zones of largest discrepancies. The distribution of these zones clearly follows a trend, making superharmonic resonances a potential suspect for this effect. The invasiveness of the SCBC with adaptive filters was also checked a posteriori, and it was also observed that the control shows the largest degree of invasiveness in the same regions, in spite of the expected action of the adaptive filter. This was not detected during the test, because the check for non-invasiveness was performed based on the error made by the signal synthesized by the adaptive filter, $x(t)-\mathbf{h}(t)\mathbf{w}(t)$.

            To ascertain that superharmonic resonances can be an issue, the 2:1 and 3:1 resonances were tracked with the PLL, by locking the phase of the second and third harmonic to $-\pi$ and $-\pi/2$, respectively\footnote{The resonant phase lag for the 2:1 resonance differs from that observed in a perfect Duffing oscillator due to the inevitable presence of asymmetries, making it more similar to a Helmholtz-Duffing oscillator in this regard.}, and varying the excitation amplitude (implementation details are given in the supplementary materials). Overlaying the obtained results\footnote{The forcing amplitude was converted to an equivalent imposed amplitude for the reference signal via $a_* = f/k_\mathrm{d} + a$.} with those of the SCBC test with Picard iterations in Figure~\ref{sfig:CBCA_2} clearly highlights that the issue indeed occurs close to the superharmonic resonances. The non-fundamental harmonic coefficients are set with a Picard iteration scheme, and a failure of this scheme in the vicinity of a resonance is coherent with numerical experience~\cite{Jain2019}.
            
        \subsubsection{Instabilities with the RCT approach}
        \label{sec:RCTIssues}
        
\begin{figure}[!ht]
                \centering
                \begin{subfigure}{.49\textwidth}
                    \centering
                    \includegraphics[width=1.1\textwidth]{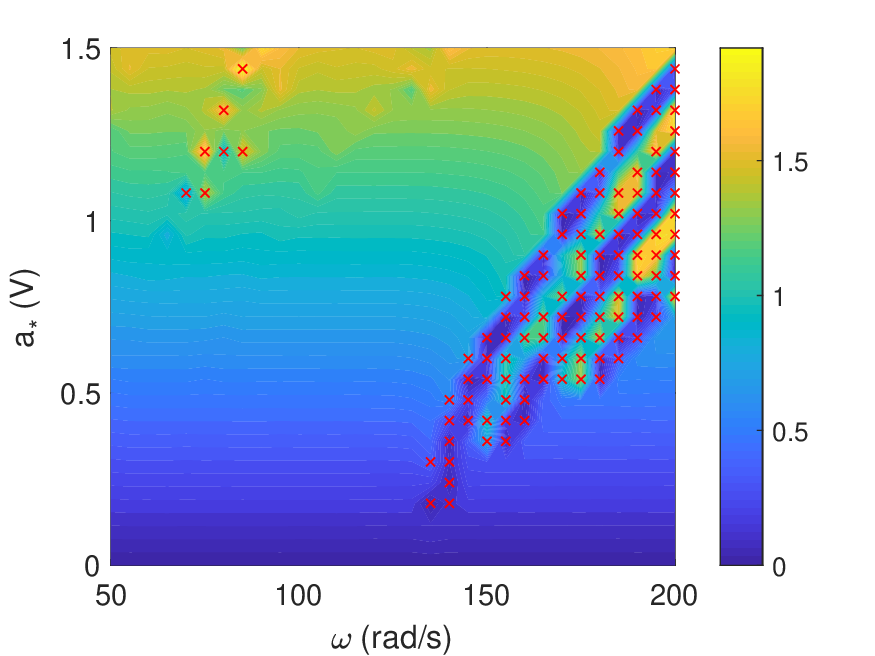}
                    \caption{}
                    \label{sfig:RCT_1}
                \end{subfigure}
                \begin{subfigure}{.49\textwidth}
                    \centering
                    \includegraphics[width=1.1\textwidth]{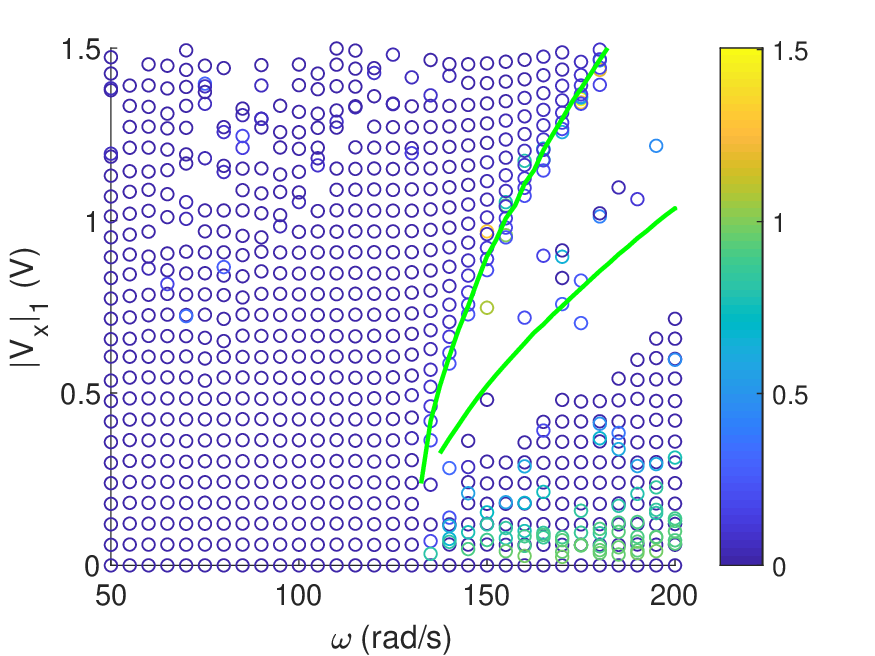}
                    \caption{}
                    \label{sfig:RCT_2}
                \end{subfigure}
                \caption{Contour plot of the actual harmonic (displacement) amplitude of the response of the RCT method (where \textcolor{red}{$\times$} indicate points with a relative error larger than 20\%)~\subref{sfig:RCT_1} and scatter plot of the frequency vs actual (displacement) amplitude (color scales represent the relative error between the target and actual amplitude) with the locus of fold bifurcations (\textcolor{MatlabGreen}{\rule[.2em]{1em}{.2em}})~\subref{sfig:RCT_2}.}
                \label{fig:RCT}
            \end{figure}

            We now turn to the incomplete aspect of the NFR curve obtained with the RCT, which is due to stability issues during the test. Figure~\ref{sfig:RCT_1} plots the actual amplitude of the displacement as a function of the target amplitude and forcing frequency. Zones with large errors are clearly apparent. The points with a relative error larger than 20\% were rejected when performing the interpolation, explaining the incomplete aspect of the NFR curve obtained with the RCT in Figure~\ref{sfig:Interpolation_2}. These zones also appear to organize into a particular pattern reminiscent of the zone of instability of the open-loop oscillator, delineated by the loci of fold bifurcations.

            To verify this hypothesis, the experimental fold bifurcation tracking procedure proposed by~\cite{Renson2017} was implemented (implementation details are given in the supplementary materials), and the resulting loci of fold bifurcations are overlaid in Figure~\ref{sfig:RCT_2}. The near complete depletion of points between the two fold bifurcation curves confirms that the issue is very likely to come from the instability of the oscillator. It thus appears that the RCT was not able to stabilize unstable solutions in this case, and the oscillatory look of the amplitude observed in the large-error region in Figure~\ref{sfig:RCT_1} corresponds to alternating convergence of the response to the upper and lower stable branches. The presence of many unstable points makes the controller struggle and systematically reach its maximum number of iterations, explaining why the RCT took about an order of magnitude longer than the other methods (cf. Table~\ref{tab:primaryResonanceTimes}).

\begin{figure}[!ht]
                \centering
                \includegraphics[width=0.539\textwidth]{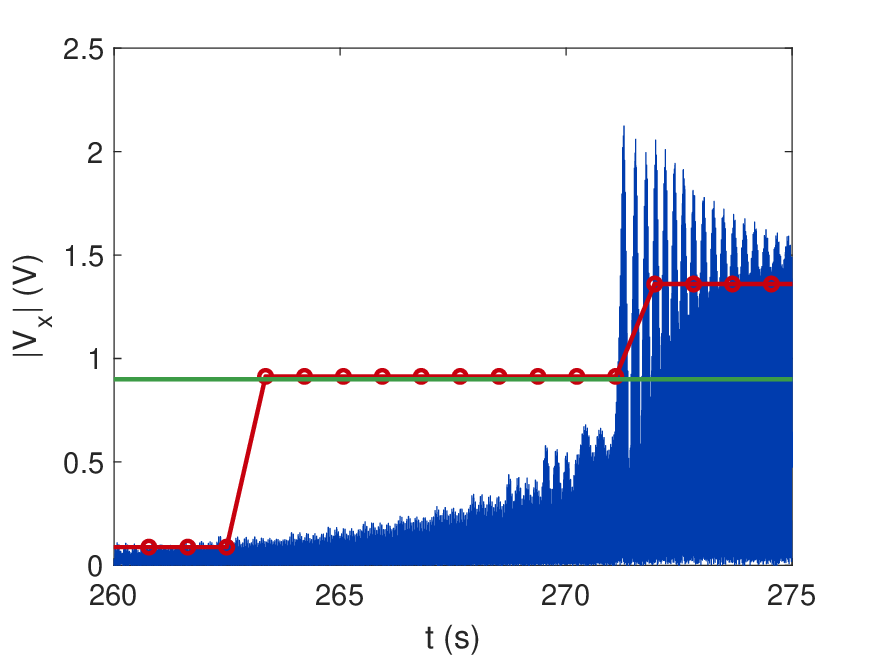}
                \caption{Amplitude of the displacement (\textcolor{MatlabCol1}{\rule[.2em]{1em}{.2em}}), estimated harmonic amplitude (\textcolor{MatlabCol2}{\rule[.2em]{1em}{.2em}}) and target amplitude (\textcolor{MatlabCol3}{\rule[.2em]{1em}{.2em}}) during the RCT at $a_* = 0.9$ V and around $\omega=170$ rad/s.}
                \label{fig:RCT_3}
            \end{figure}

             Some points lay inside the unstable zone delineated in Figure~\ref{sfig:RCT_2}, and an analysis of the time series revealed that they are associated with transient phenomena. The resulting amplitude thus cannot be considered as a steady-state quantity. For illustration, Figure~\ref{fig:RCT_3} depicts the time series associated with the point $(\omega,a_*)=(170,0.9)$ for which the amplitude indicated by the TestLab software results is $|V_\mathrm{x}|_1 = 0.9135$, hence seemingly correct enough. The time series of the amplitude estimated by TestLab was unfortunately unavailable and was retrieved through time-frequency analysis (using the short-time Fourier transform from the \texttt{spectrum} function in Matlab). Its estimation is overlaid together with the target amplitude. It appears that the amplitude of the considered point, indicated by the middle step in Figure~\ref{fig:RCT_3}, is calculated during transients (right before the end of that step) and is thus not representative of the steady-state response of the oscillator. Assessing this in a systematic way is unfortunately difficult given the little amount of information that can be extracted from TestLab results. 

            The accumulation of points at the left of the leftmost branch, and the lack of points at the right of the rightmost branch in Figure~\ref{sfig:RCT_2}, which are both stable regions, was first thought to be due to the stepped-up frequency during the test. However, a stepped-down frequency surprisingly showed the same trend. The reason for this puzzling trend was not investigated further.
            
            Finally, we note that there are points in the upper left corner of Figure~\ref{sfig:RCT_1} which also feature a large error. Similarly to the case exposed in Figure~\ref{sfig:CBCA_2}, this is most likely due to superharmonic resonances, where folding is also expected to occur at high amplitudes.
            
            Because the RCT works with a black-box commercial software, there is no potential implementation issue. Different parameters were tried out without success, as a zone of instability always seemed to appear in the same range (more details are given in the supplementary materials). The authors believe that the main difference between this work and the original~\cite{Karaagacl2020,Karaagacl2021} is the presence of an electrodynamic shaker. This is known to stabilize unstable responses in some cases~\cite{Robbins2023}, which might be the case therein. This issue would deserve further investigation but was deemed beyond the scope of this work. Nevertheless, if this explanation is correct, the extraction of the folded part of the NFR of the structure alone is more related to continuation aspects than stabilization ones, since the RCT itself performs continuation on curves that do not fold but does not appear to inherently change the stability of the system under test (composed of the structure and the shaker). 

        \subsubsection{Testing times}

        We conclude this first series of tests by providing the times needed to perform each experiment in Table~\ref{tab:primaryResonanceTimes}. In terms of testing times, the fact that we focused on one NFR puts direct NFR measurement techniques, namely SWS, PLL, CBC-FD and ACBC, at an advantage compared to the other methods. It should thus be kept in mind when looking at the total time taken by each methods that SCBC and RCT allow for the extraction of multiple NFRs in one measurement run, and a more meaningful measure of the time efficiency of a method can be given by the average time taken per measurement point.

        \begin{table*}[!ht]
        \centering
        \caption{Times taken by each test for the primary resonance.}
        \label{tab:primaryResonanceTimes}
        \begin{tabular}{ccc}
             \hline
             Method & Total time (s) & Time per point (s) \\
             \hline
             SWS (up and down)& 240 & / \\
             CBC-FD & 396 & 4.96 \\
             SCBC (Picard) & 453 & 0.62 \\
             SCBC (adaptive) & 397 & 0.55 \\
             ACBC & 240 & 3.19 \\
             PLL & 240 & / \\
             RCT & 7 847 & 10.13 \\
             \hline
        \end{tabular}
        \end{table*}

    \subsection{Superharmonic resonance}
    \label{sec:DuffingSuperH}

        We now turn to the 3:1 superharmonic resonance of the Duffing oscillator. This type of resonance is seldom studied in control-based methods, but actually challenges all of them.

        \subsubsection{Numerical illustration}
        \label{ssec:DuffingSuperHNumerical}

\begin{figure*}[!ht]
                \centering
                \begin{subfigure}{.49\textwidth}
                    \centering
                    \includegraphics[width=1.1\textwidth]{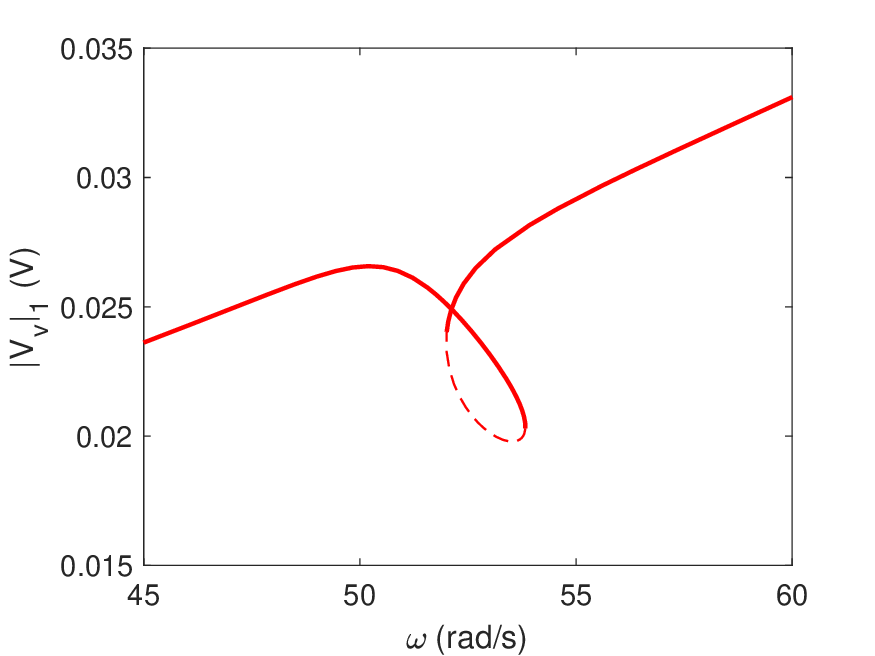}
                    \caption{}
                    \label{sfig:Duffing_3_1_F}
                \end{subfigure}
                \begin{subfigure}{.49\textwidth}
                    \centering
                    \includegraphics[width=1.1\textwidth]{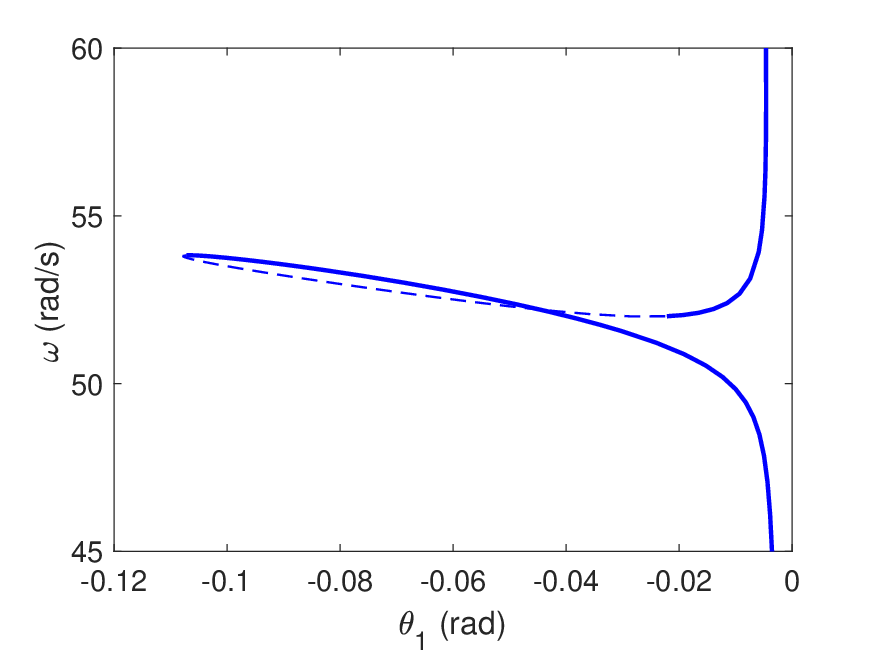}
                    \caption{}
                    \label{sfig:Duffing_3_1_T}
                \end{subfigure}
                
                \begin{subfigure}{.49\textwidth}
                    \centering
                    \includegraphics[width=1.1\textwidth]{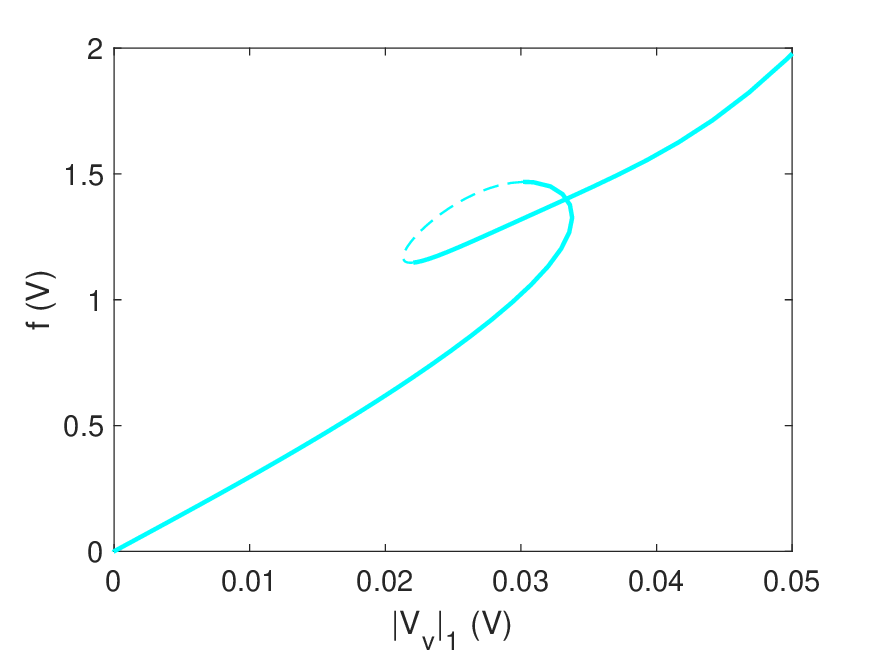}
                    \caption{}
                    \label{sfig:Duffing_3_1_W}
                \end{subfigure}
                \begin{subfigure}{.49\textwidth}
                    \centering
                    \includegraphics[width=1.1\textwidth]{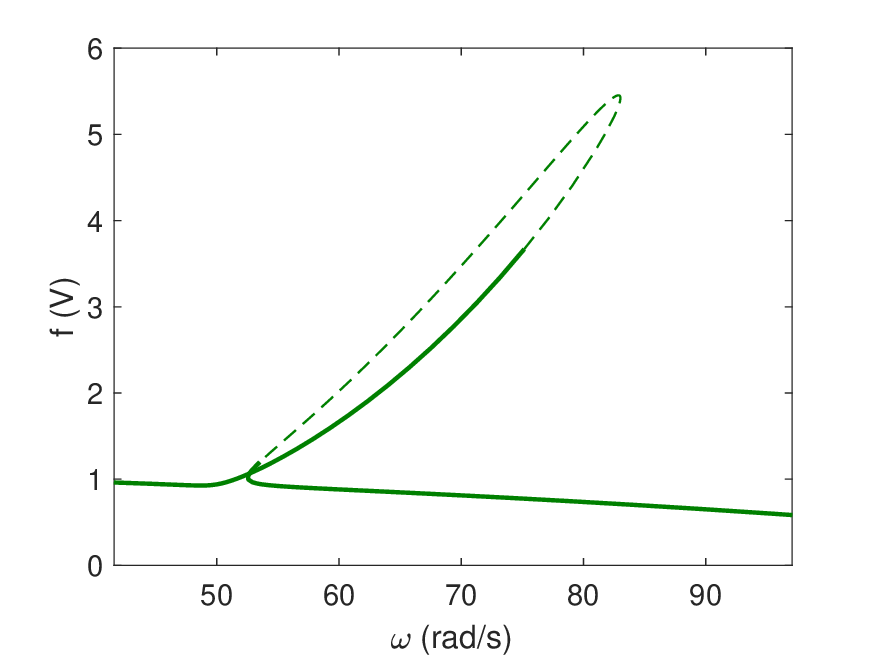}
                    \caption{}
                    \label{sfig:Duffing_3_1_A}
                \end{subfigure}
                \caption{Numerical simulation of the response of the first harmonic of the Duffing oscillator with $f=1$ V (\subref{sfig:Duffing_3_1_F}: frequency parametrization, \subref{sfig:Duffing_3_1_T}: phase parametrization), with $\omega = 55.5$ rad/s \subref{sfig:Duffing_3_1_W}, and $|V_\mathrm{x}|_1 = 0.5$ V \subref{sfig:Duffing_3_1_A}.}
                \label{fig:Duffing_3_1_Num}
            \end{figure*}

        Before moving forward with the experiments, a numerical study has been made around the 3:1 superharmonic resonance of the Duffing oscillator. Variants of a harmonic balance method coupled with numerical continuation~\cite{Detroux2015} have been used to compute constant forcing amplitude, forcing frequency and response amplitude curves with the theoretical parameters in Table~\ref{tab:DuffingParameters}. They are displayed in Figure~\ref{fig:Duffing_3_1_Num}. 

        Clearly, all methods based on natural parameter continuation, namely SWS, SCBC, PLL and RCT will fail, because folding is present in all cases. Similar observations were made in~\cite{Abeloos2022Thesis}. As for the ACBC approach, success is a priori uncertain given that the first harmonic amplitude curve features a loop that does not form a contour level.
    
\begin{figure}[!ht]
                \centering
                \begin{subfigure}{.49\textwidth}
                    \centering
                    \includegraphics[width=1.1\textwidth]{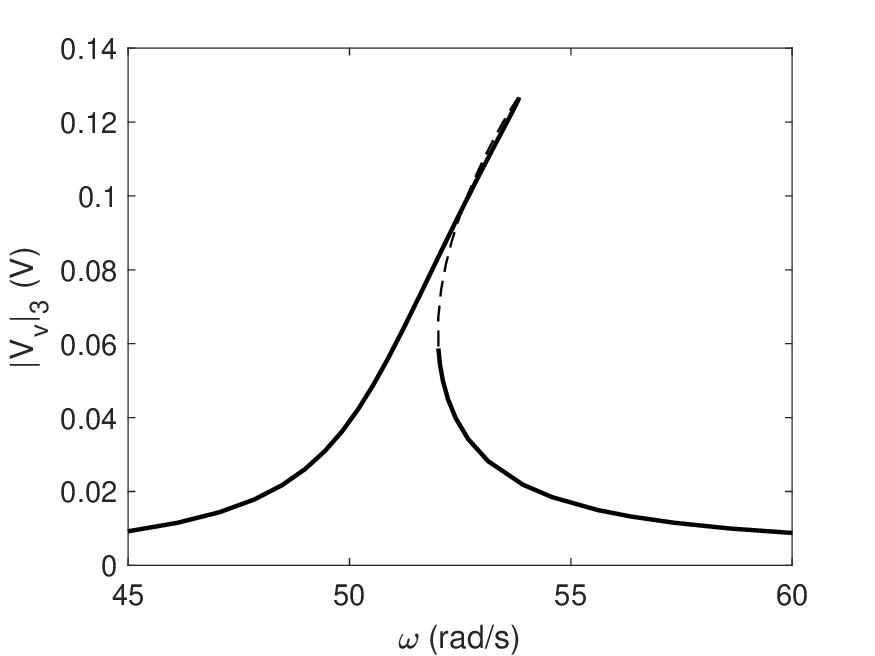}
                    \caption{}
                    \label{sfig:Duffing_3_1_F3}
                \end{subfigure}
                \begin{subfigure}{.49\textwidth}
                    \centering
                    \includegraphics[width=1.1\textwidth]{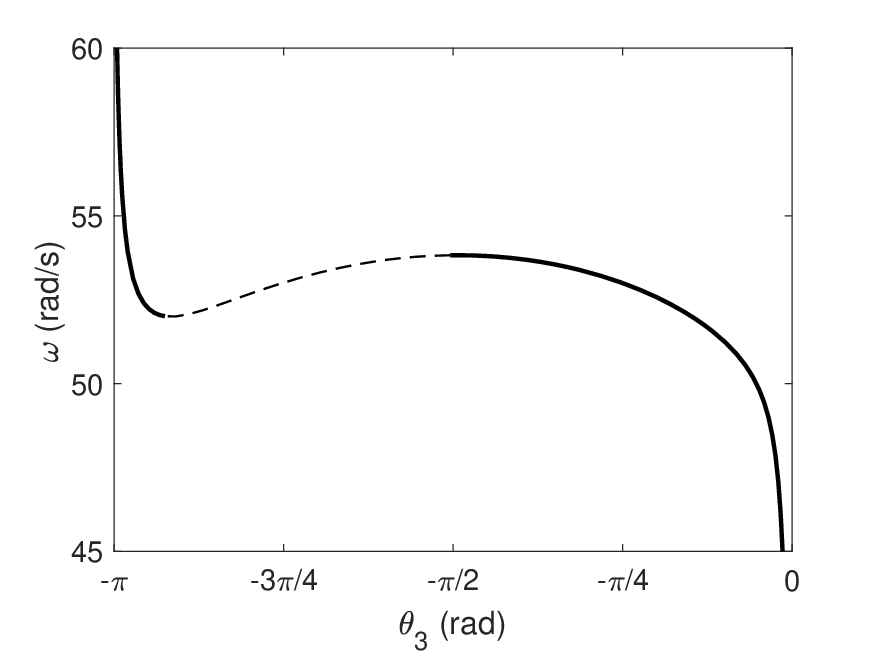}
                    \caption{}
                    \label{sfig:Duffing_3_1_T3}
                \end{subfigure}
                \caption{Numerical simulation of the response of the third harmonic of the Duffing oscillator with $f = 1$ V (\subref{sfig:Duffing_3_1_F3}: amplitude, \subref{sfig:Duffing_3_1_T3}: phase).}
                \label{fig:Duffing_3_1_Num3}
            \end{figure}

        A solution nevertheless exists for the PLL, as exploited in~\cite{zhou2024identification}, which consists in replacing the phase lag of the first harmonic by that of the third harmonic. As Figure~\ref{sfig:Duffing_3_1_T3} evidences, such a parametrization does not exhibit any fold in the vicinity of the 3:1 resonance, making it a viable approach. 

\begin{figure}[!ht]
            \centering
                \includegraphics[width=0.539\textwidth]{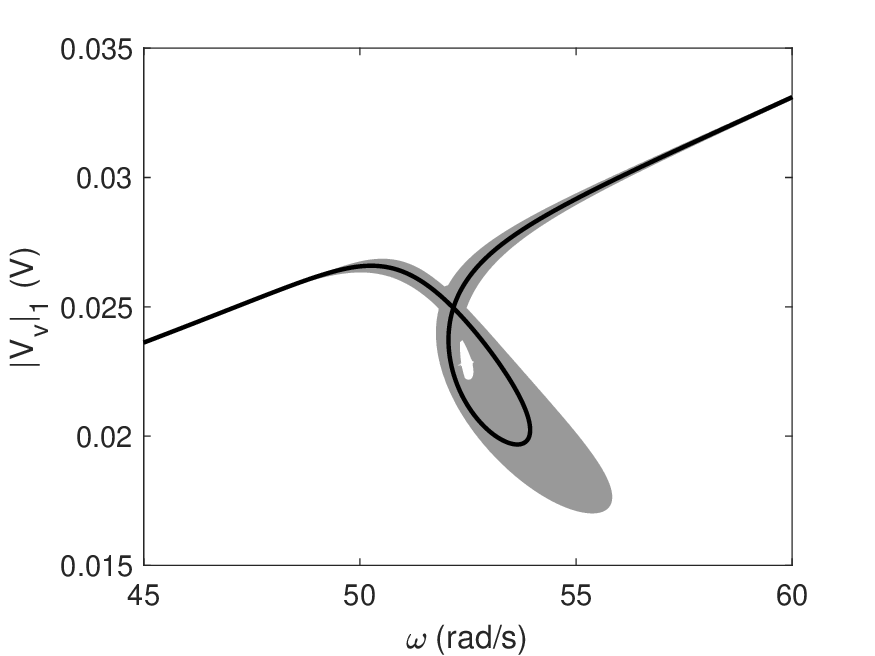}
            \caption{Numerical simulation of the response of the first harmonic of the Duffing oscillator with nominal mono-harmonic forcing with $f=1$ V (\rule[.2em]{1em}{.2em}) and locus of responses with a 1\% additive forcing amplitude and arbitrary phase on the third harmonic (\textcolor{Gray1}{\rule[.2em]{1em}{.2em}}).}
            \label{fig:DuffingInvasiveness_3_1}
        \end{figure}

        Finally, we noted experimental difficulties with the CBC-FD method as well. In the vicinity of the 3:1 resonance, non-invasiveness was found to be crucial for the accuracy of the results and for the success of the method itself. To investigate the impact of invasiveness, a series of numerical simulations were performed by imposing a non-tonal forcing to the oscillator. To the nominal fundamental harmonic forcing was added a forcing on the third harmonic of 1\% relative amplitude and fixed phase. Performing the continuation for multiple values of the phase, we obtained an idea of locus of responses if a 1\% relative tolerance is allowed on the residual, as depicted in Figure~\ref{fig:DuffingInvasiveness_3_1}. In spite of the small tolerance, the loop associated with the resonance can vary by up to about 50\% in size, stressing the importance of a non-invasive control. 

        Numerical simulations of the Duffing oscillator tested with the CBC-FD uncovered a bottleneck linked to this issue. If too large a relative tolerance was chosen (such as 1\%), the continuation procedure was progressing in an erratic fashion inside a locus similar to that depicted in Figure~\ref{fig:DuffingInvasiveness_3_1}. The resulting irregular curve posed challenges for the continuation by making inappropriate predictions which could not be corrected. Alternatively, choosing a small tolerance made the curve less erratic. However, inaccuracies on the residual (coming from the Fourier decomposition) were observed to increase the vicinity of the resonance, eventually exceeding the tolerance threshold, rendering the convergence of the method impossible.

        \subsubsection{Experimental testing}

            In view of the issues outlined in Sections~\ref{sec:CBCAIssues},~\ref{sec:RCTIssues} and~\ref{ssec:DuffingSuperHNumerical}, the SCBC and RCT methods were not deemed suitable for the 3:1 superharmonic resonance. A SWS was performed in spite of its deficiencies to obtain a reference method, together with the CBC-FD, ACBC and PLL (with the phase of the third harmonic), with a forcing amplitude of 1 V. All parameters are given in the supplementary materials.
            
            \begin{figure}[!ht]
                \centering
                \begin{subfigure}{.49\textwidth}
                    \centering
                    \includegraphics[width=1.1\textwidth]{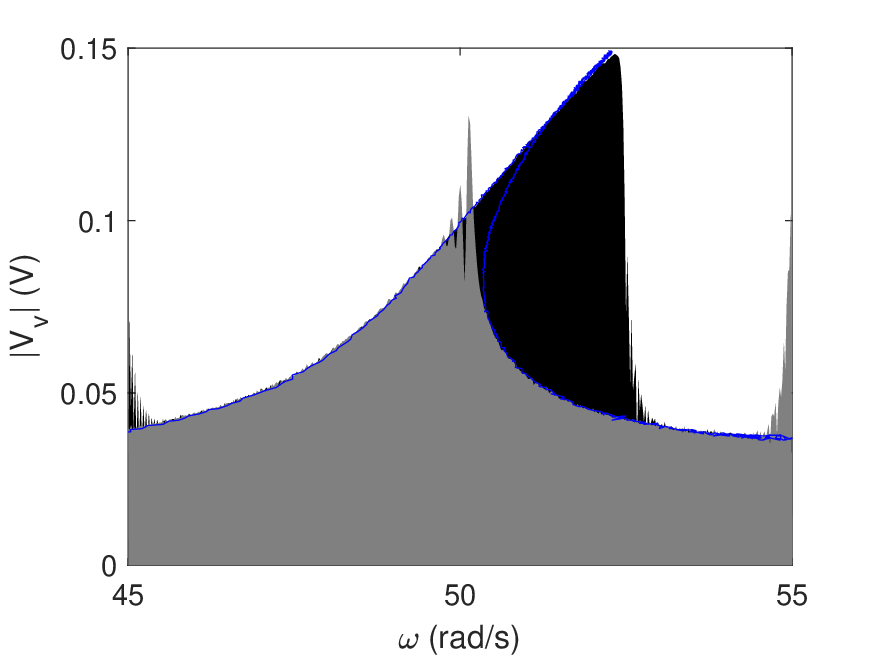}
                    \caption{}
                    \label{sfig:NFR_3_1}
                \end{subfigure}
                \begin{subfigure}{.49\textwidth}
                    \centering
                    \includegraphics[width=1.1\textwidth]{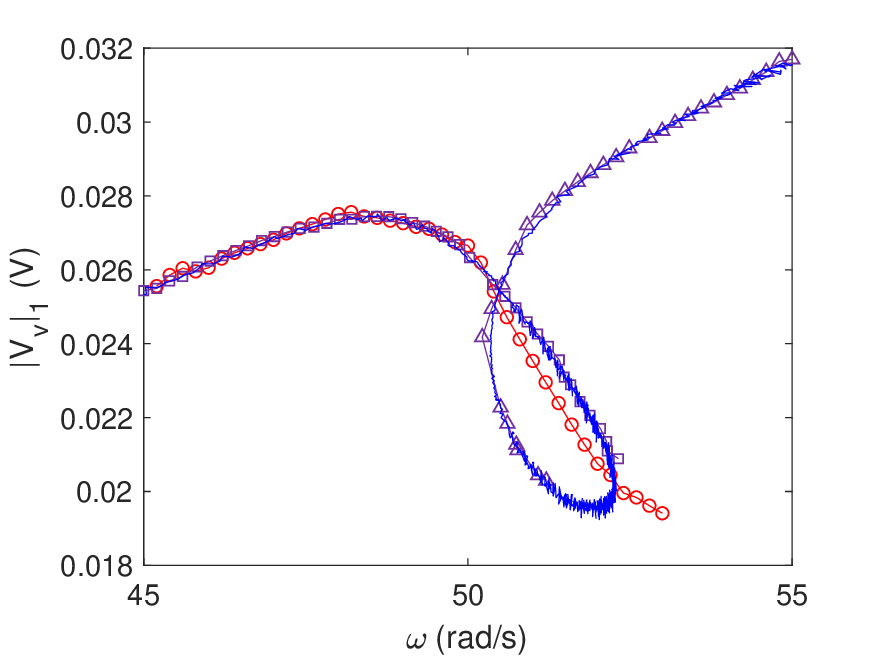}
                    \caption{}
                    \label{sfig:NFR_3_1_2}
                \end{subfigure}
                \caption{\subref{sfig:NFR_3_1}: Velocity NFR around the 3:1 superharmonic resonance of the electronic Duffing oscillator obtained with the PLL (\textcolor{blue}{\rule[.2em]{1em}{.2em}}), SWS-up (\rule[.2em]{1em}{.2em}) and SWS-down (\textcolor{gray}{\rule[.2em]{1em}{.2em}}). \subref{sfig:NFR_3_1_2}: fundamental harmonic of the velocity signal around the 3:1 superharmonic resonance of the electronic Duffing oscillator obtained with the PLL (\textcolor{blue}{\rule[.2em]{1em}{.2em}}), ACBC (\textcolor{red}{\rule[.2em]{1em}{.2em}}) and CBC-FD (\textcolor{MatlabPurple}{\rule[.2em]{1em}{.2em}}, $\square$: starting from 45 rad/s, $\triangle$: starting from 55 rad/s) methods.}
                \label{fig:NFR_3_1}
            \end{figure}    

            Figure~\ref{sfig:NFR_3_1} shows the results of the most successful method, namely the PLL, for the 3:1 resonance, together with swept-sine results for reference. The PLL is able to robustly go through the unstable branches and trace the complete NFR, and it was the only method able to do so thanks to the use of the phase of the third harmonic.

            The NFR follows a loop, which is easier to see in Figure~\ref{sfig:NFR_3_1_2} where only the first harmonic amplitude is depicted. The results of the CBC-FD and ACBC are also shown therein. For the former, convergence of the Newton-Raphson procedure did not occur at some point. This lack of convergence could be due to the increased noise level observed around the 3:1 resonance, making the finite-difference estimation of the Jacobian inaccurate. As for the ACBC, the control was seen to be invasive, signalling that the adaptive filters were failing to fulfill their purpose, as in Section~\ref{sec:CBCAIssues}. The estimated NFR consequently significantly diverges from what is estimated to be its true shape, similarly to the issue raised with Figure~\ref{fig:DuffingInvasiveness_3_1}. In addition, the fact that the NFR loops around the 3:1 resonance (Figure~\ref{sfig:Duffing_3_1_F}) indicates that the HFS cannot be parametrized solely by $\omega$ and $a_*$ in this case, making the proposed ellipse sweeping procedure uncertain.

            The success of the PLL based on the phase lag of the third harmonic could lead to think that the SCBC and RCT would potentially benefit from a similar unfolding provided they work with the amplitude of the third harmonic. Unfortunately, as Figure~\ref{sfig:Duffing_3_1_F3} reveals, such an unfolding does not necessarily occur for this harmonic amplitude, making it uncertain whether these methods would be able to handle superharmonic resonances.

            This series of tests on a superharmonic resonance showed that the assumptions of all methods based on natural parameter continuation can fail even with a single-degree-of-freedom system. Such resonances are nevertheless usually less important and/or relevant than primary resonances, that remain the main interest in nonlinear vibration testing. 

\section{Illustration with a clamped thin plate}
\label{sec:Plate}

    We now move to a real structure with challenging nonlinear dynamics, namely multiple nonlinear modes with a distributed geometrical nonlinearity. This set-up does not feature the same exceptional repeatability as the electronic Duffing oscillator and therefore gives a more realistic idea of what can be expected from a test with a control-based method.

    \subsection{Experimental setup}

        A thin rectangular steel plate of dimensions 500 mm $\times$ 300 mm $\times$ 0.5 mm and clamped on its shortest edges was considered as the next experimental setup. The plate was excited with an electrodynamic shaker (TIRA TV 51075) in current mode, and the response of the plate was monitored with an impedance head (DYTRAN 5860B) as well as a laser vibrometer (Polytec NLV-2500-5). In addition to these measurements from the structure, the voltage across the electrodes of the shaker was also monitored.
        
        \begin{figure*}[!ht]
            \centering
            \begin{subfigure}{.65\textwidth}
                \centering
                \includegraphics[scale=1.]{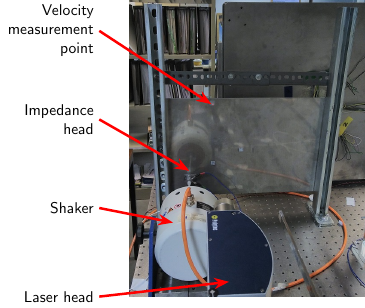}
                \caption{}
                \label{sfig:plate}
            \end{subfigure}
            \begin{subfigure}{.25\textwidth}
                \centering
                \includegraphics[height=5cm]{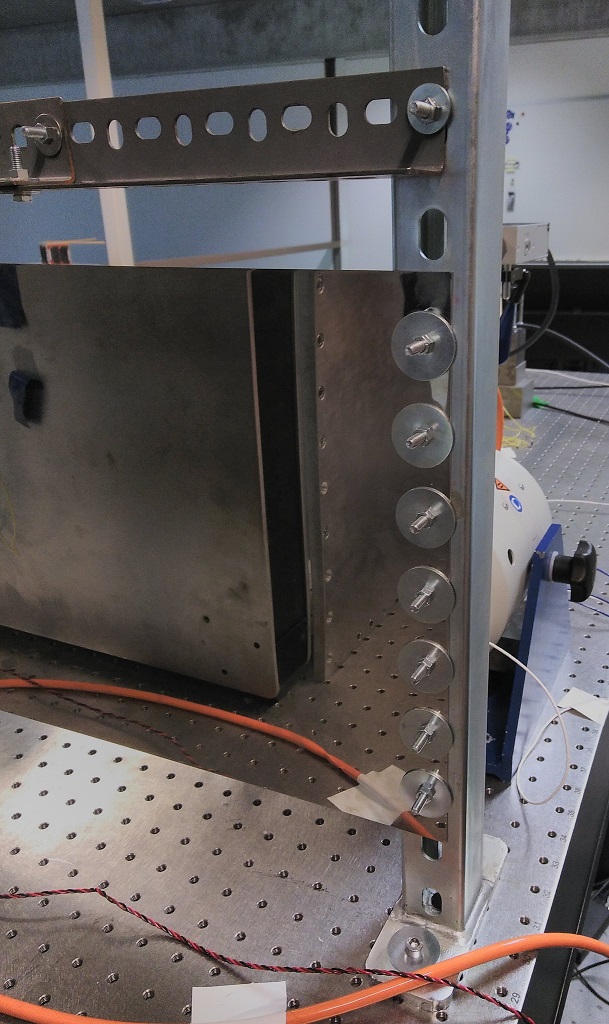}
                \caption{}
                \label{sfig:backside}
            \end{subfigure}
            \caption{Pictures of the experimental setup: \subref{sfig:plate}: overview of the clamped plate and \subref{sfig:backside}: view of the backside of the clamping device.}
            \label{fig:plate}
        \end{figure*} 

        Figure~\ref{fig:plate} shows pictures of the experimental clamped plate. The clamping device was realized with two cantilever arms fixed on an optical table. The plate was then fixed to the arms with bolts, nuts and wide washers to distribute the load as evenly as possible on the clamped side. The cantilever arms were observed to be rigid enough to provide an approximate clamping condition for the bending dynamics of the plate. However, they were not rigid enough with regards to the plate axial stiffness, and a bar was added on top of the plate to stiffen the assembly, allowing for the observation of clear nonlinear midplane stretching effects. Details about the dimensions of the setup are given in the supplementary materials.

        Compared with the electronic Duffing oscillator, several new challenges arise. The most problematic one is the occurrence of shaker-structure interactions. These interactions have two main consequences when a harmonic drive signal is fed to the amplifier of the shaker: (\textit{i}) the fundamental forcing harmonic amplitude is frequency- and amplitude-dependent and (\textit{ii}) harmonics appear in the force. These issues were overlooked in this work by considering the signal fed to the amplifier, $V_\mathrm{in}$ (proportional to the current in the shaker), as the input to the system (as opposed to the force measured by the impedance head). The price to pay for this convenience is that shaker-structure interactions can then substantially dampen the structural modes if $V_\mathrm{in}$ is considered as the input. The shaker was thus placed following a trade-off between conflicting requirements: on the one hand, it should be placed on a responsive location of the targeted modes to excite them well, but on the other hand it should be placed close to the clamping to minimize shaker-structure interactions.

    \subsection{Linear modal analysis}

        A low-level pseudorandom excitation signal was applied to the plate to compute the FRFs when the system behaves linearly. The PolyMAX modal identification method~\cite{Peeters2004} was then used to extract the linear modal characteristics of the first nine modes (below 100 Hz). The characteristics of the first five modes, which will be tested in nonlinear regimes of motion, are gathered in Table~\ref{tab:modalCharacteristics}.

        \begin{table}[!ht]
        \centering
        \caption{Modal characteristics of the first five modes of the plate.}
        \label{tab:modalCharacteristics}
        \begin{tabular}{ccc}
             \hline
             Mode & Frequency (Hz) & Damping ratio (\%) \\
             \hline
             1 & 15.51 & 0.21\\
             2 & 25.86 & 0.89 \\
             3 & 41.11 & 0.91 \\
             4 & 49.02 & 0.33 \\
             5 & 52.90 & 0.33 \\
             \hline
        \end{tabular}
        \end{table}        

        \begin{figure*}[!ht]
            \centering
            \includegraphics[width=\textwidth]{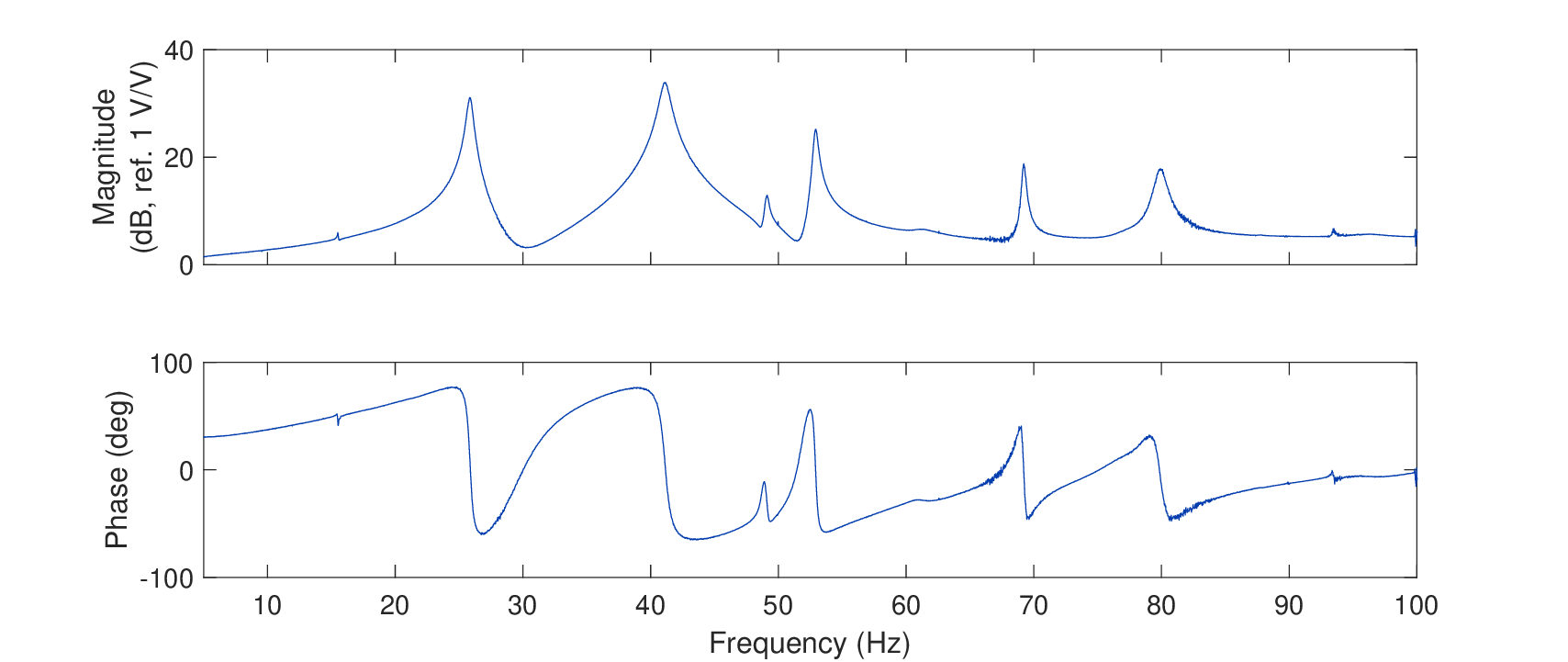}
            \caption{Bode plot of the shaker voltage FRF.}
            \label{fig:PlateFRF_Voltage}
        \end{figure*}

        \begin{figure*}[!ht]
            \centering
            \includegraphics[width=\textwidth]{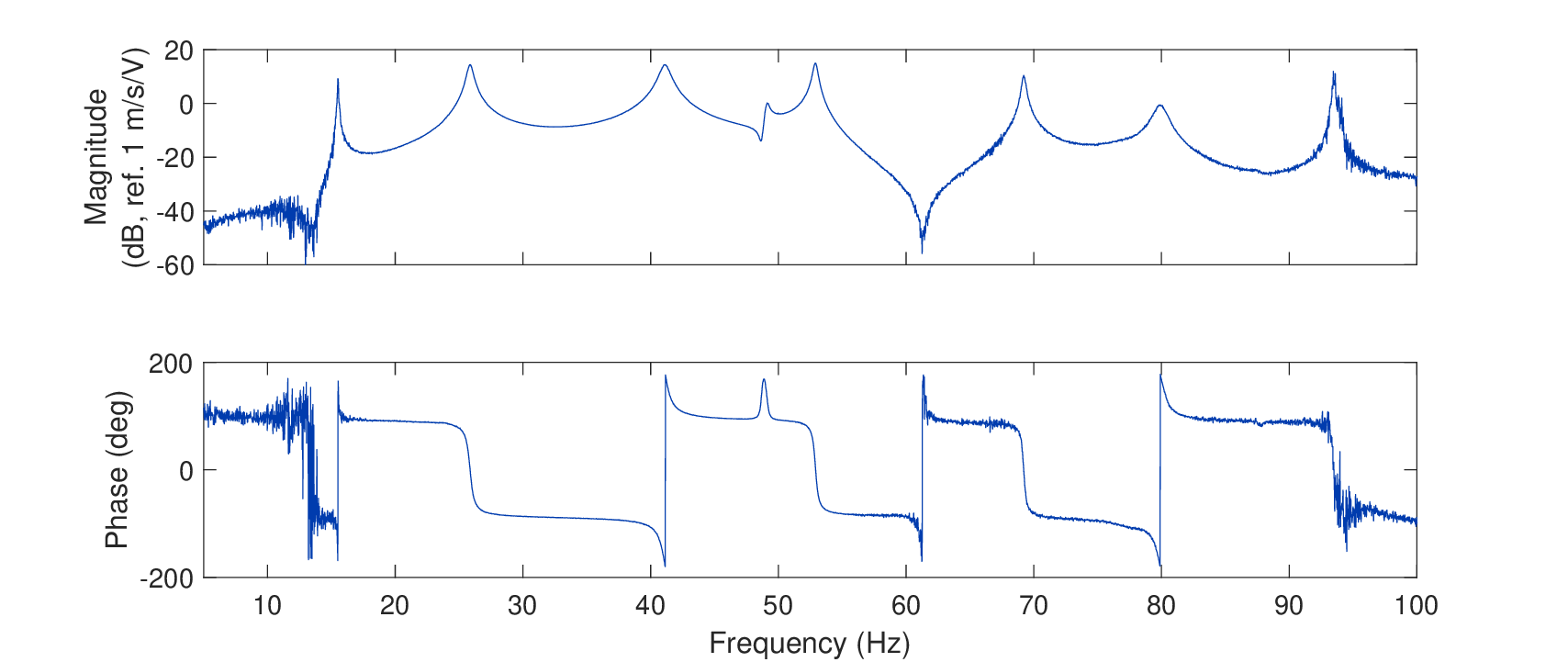}
            \caption{Bode plot of the laser velocity FRF.}
            \label{fig:PlateFRF_Velocity}
        \end{figure*}

        Figures~\ref{fig:PlateFRF_Voltage} and~\ref{fig:PlateFRF_Velocity} show the two FRFs related to the NFRs of interest in the sequel, namely the shaker voltage and laser velocity.

        \begin{figure*}[!ht]
            \centering
            \includegraphics[width=\textwidth]{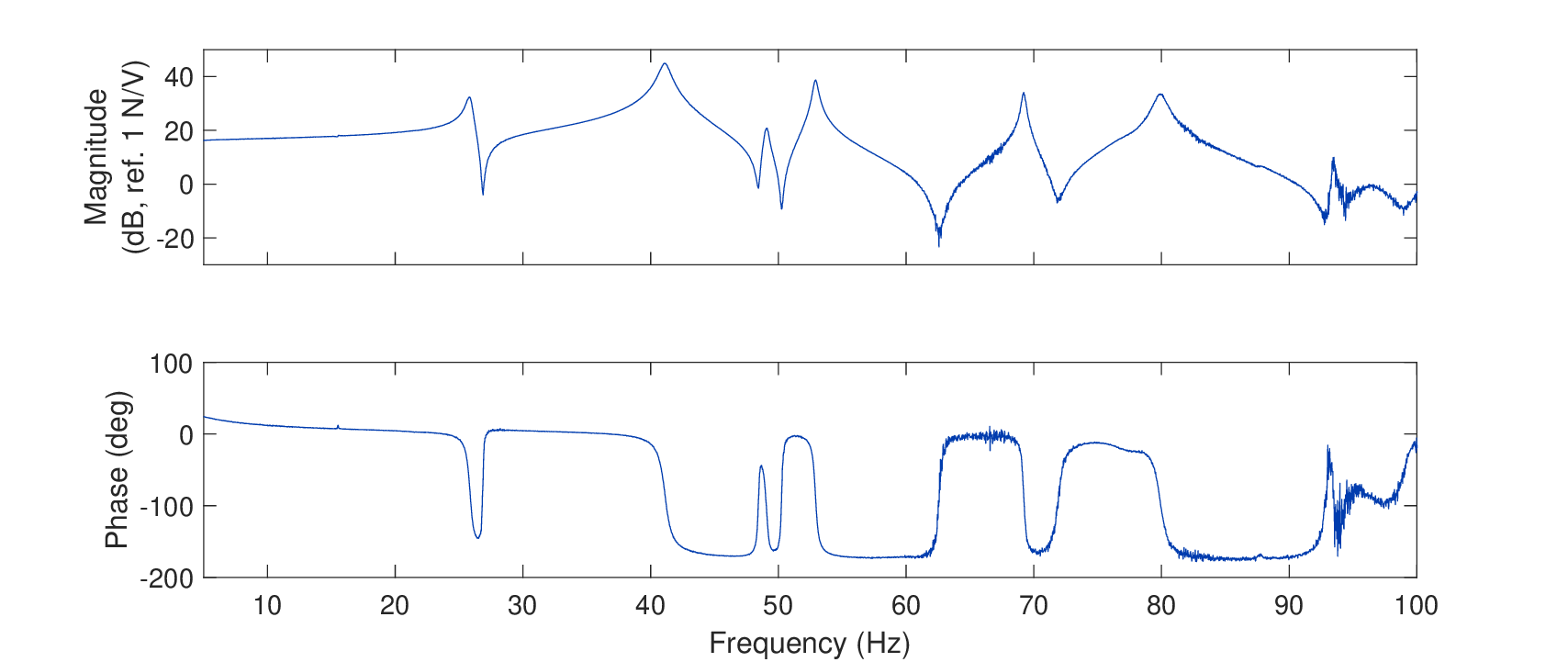}
            \caption{Bode plot of the impedance head force FRF.}
            \label{fig:PlateFRF_Force}
        \end{figure*}

        Figure~\ref{fig:PlateFRF_Force} shows the FRF of the force applied to the structure, as measured by the impedance head. Clear shaker-structure interactions are visible, as the FRF is far from being constant. This is due to the high flexibility of the plate combined with the relatively high mass of the impedance head~\cite{Ewins2009}. Nevertheless, the shaker does not appear to dampen the modes of interest excessively. Controlling the force would prove at best extremely challenging anyway, because this FRF has a nonminimum phase zero (a pair of complex conjugate zeros with positive real part) around 137 Hz, which is very likely to destabilize any attempt of linear feedback control~\cite{Preumont2011}.

        The dynamics of the plate in higher-frequency regions feature a high modal density. While this might not be a limitation for several methods, preliminary SWS tests showed that, in a nonlinear regime of motion, extremely complex phenomena took place therein, with the appearance of internal and combination resonances. Both these phenomena are generally associated with the appearance of quasiperiodic oscillations, that none of the state-of-the-art control-based method can handle. This high frequency range was thus not investigated further in this work.

    \subsection{Implementation details}

        The CBC requires feedback control similar to direct velocity feedback. As is well known in the field of linear control of flexible structures, collocation of the applied forcing and the feedback signal is almost mandatory to hope for a stable closed-loop system~\cite{Preumont2011}. Loosely speaking, collocation allows for the feedback controller to work with power conjugate quantities, thereby controlling the power flow out of the structure. In particular, the control law can be designed to avoid power input to the structure which could lead to closed-loop instabilities. In most works about the CBC, one-degree-of-freedom-like structures were considered, thereby dismissing this question. For the plate, flexibility effects bear a strong importance. To allow for collocated control, one possibility is to measure the motion of the moving part of the shaker~\cite{Abeloos2021,Quaegebeur2023}. Alternatively, the strategy used in this work consists in sensing the voltage across the electrodes of the shaker (which is measured by the shaker amplifier). Since the current and voltage are power conjugates, a truly collocated control is achieved. More specifically, this is confirmed in Figure~\ref{fig:PlateFRF_Voltage} , where the phase bounded between $\pm 90^\circ$ ensures an interlacing pole-zero pattern, typical of a collocated sensor-actuator pair~\cite{Preumont2011}. A first-order, unit-gain low-pass filter with a cut-off frequency at 500 Hz was nevertheless used on the voltage signal to avoid the excitation of high-frequency dynamics and the undesirable effects associated with digital control.

        For the PLL and RCT, such a collocation is not necessary. Since the phase evolution around the first considered mode shows better monotonicity for the velocity signal (cf. Figure~\ref{fig:PlateFRF_Velocity}), the latter was chosen as the feedback signal for PLL.

        Finally, the CBC-FD and RCT methods were unsuccessful. Concerning the former, this is probably due to the stronger noise in this setup. This lead the method to diverge quite far away from the NFR at the start of the test at low amplitudes. Since this behavior was worrisome near resonance, the CBC-FD was not tested in this case. Better results can be hoped for if noise is handled properly~\cite{Schilder2015}, but this is beyond the scope of the present work. The RCT method also had stability issues similar to those encountered for the electronic Duffing. However, since these systematically lead to an overload, the RCT could not be performed fully to guarantee the structural integrity of the plate.

        The tested methods, namely SCBC, PLL and ACBC, were compared once again with a SWS for reference. All parameters are given in the supplementary materials.

        \subsection{Broadband NFR}

        We start off the nonlinear analysis with the aim to obtain an NFR encompassing multiple resonant modes by investigating the frequency range between 10 and 60 Hz. In this case, several of the control-based methods are impractical or even inapplicable. The PLL is inherently hampered by the fact that the phase is non-monotonous (cf. Figures~\ref{fig:PlateFRF_Voltage} and~\ref{fig:PlateFRF_Velocity}). In addition, the phase sensitivity away from resonance is rather small, indicating that phase parametrization is inadequate, and phase noise may cause large frequency noise. As for SCBC and RCT, their main drawback is that they require a certain response amplitude at all frequencies. Near antiresonances, reaching such an amplitude would require unreasonably large forcing levels. Furthermore, since the frequency band is quite large, the structure would have to be excited at maximum response amplitude for a long time, which is obviously undesirable. Hence, while these three methods are suitable to measure the response of a structure around its resonance, they are unsuited for broadband NFR testing. This could nevertheless be different for the two latter if an adapted sampling procedure was chosen instead of a regular grid in the ($\omega, a_*$) space.

        \begin{figure*}[!ht]
            \centering
            \includegraphics[width=\textwidth]{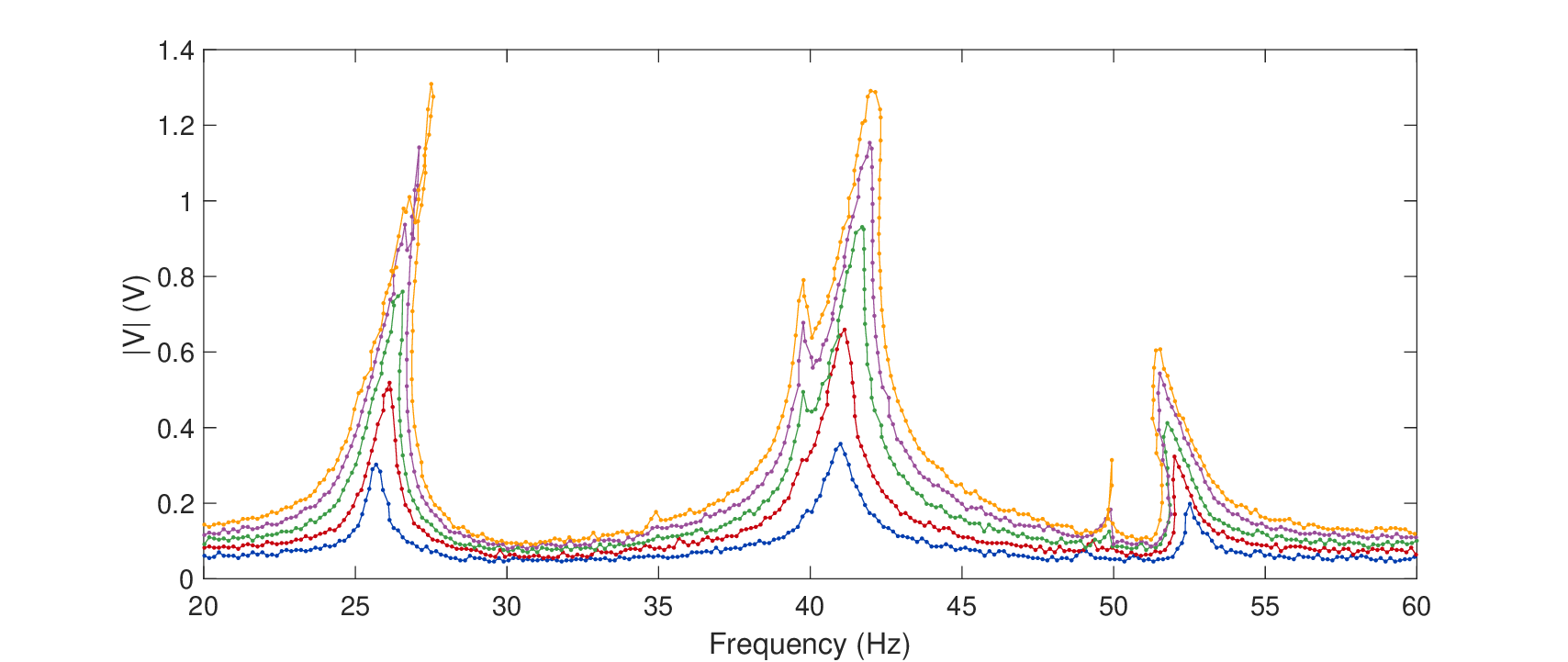}
            \caption{Shaker voltage NFR of the clamped plate for different drive signal amplitudes: $0.01$ V (\textcolor{MatlabCol1}{\rule[.2em]{1em}{.2em}}), $0.02$ V (\textcolor{MatlabCol2}{\rule[.2em]{1em}{.2em}}), $0.03$ V (\textcolor{MatlabCol3}{\rule[.2em]{1em}{.2em}}), $0.04$ V (\textcolor{MatlabCol4}{\rule[.2em]{1em}{.2em}}) and $0.05$ V (\textcolor{MatlabCol5}{\rule[.2em]{1em}{.2em}}).}
            \label{fig:PlateNFRs_Voltage}
        \end{figure*}

        \begin{figure*}[!ht]
            \centering
            \includegraphics[width=\textwidth]{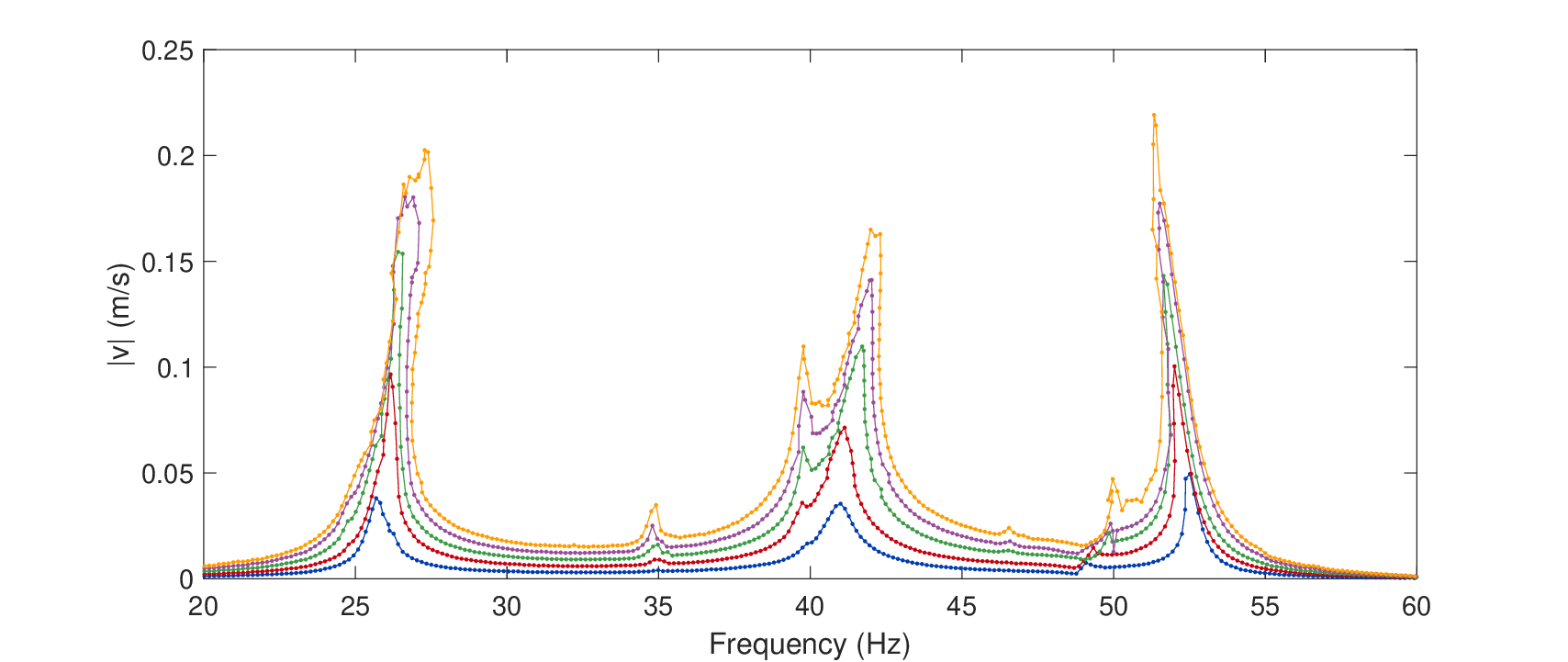}
            \caption{Laser velocity NFR of the clamped plate for different drive signal amplitudes: $0.01$ V (\textcolor{MatlabCol1}{\rule[.2em]{1em}{.2em}}), $0.02$ V (\textcolor{MatlabCol2}{\rule[.2em]{1em}{.2em}}), $0.03$ V (\textcolor{MatlabCol3}{\rule[.2em]{1em}{.2em}}), $0.04$ V (\textcolor{MatlabCol4}{\rule[.2em]{1em}{.2em}}) and $0.05$ V (\textcolor{MatlabCol5}{\rule[.2em]{1em}{.2em}}).}
            \label{fig:PlateNFRs_Velocity}
        \end{figure*}

        Using the ACBC for different levels of excitation, one can obtain a global picture of the plate dynamics. Figures~\ref{fig:PlateNFRs_Voltage} and~\ref{fig:PlateNFRs_Velocity} feature the voltage and velocity NFRs, respectively, for five different drive amplitudes. Although the test started at 10 Hz, the frequency range shown therein starts at 20 Hz because the first mode is barely visible and barely excited by the shaker (as seen in Figure~\ref{fig:PlateFRF_Voltage}), thus remaining linear and being rather uninteresting. Each of these NFRs were obtained under 8 minutes of testing.

        \begin{figure*}[!ht]
            \centering
            \includegraphics[width=\textwidth]{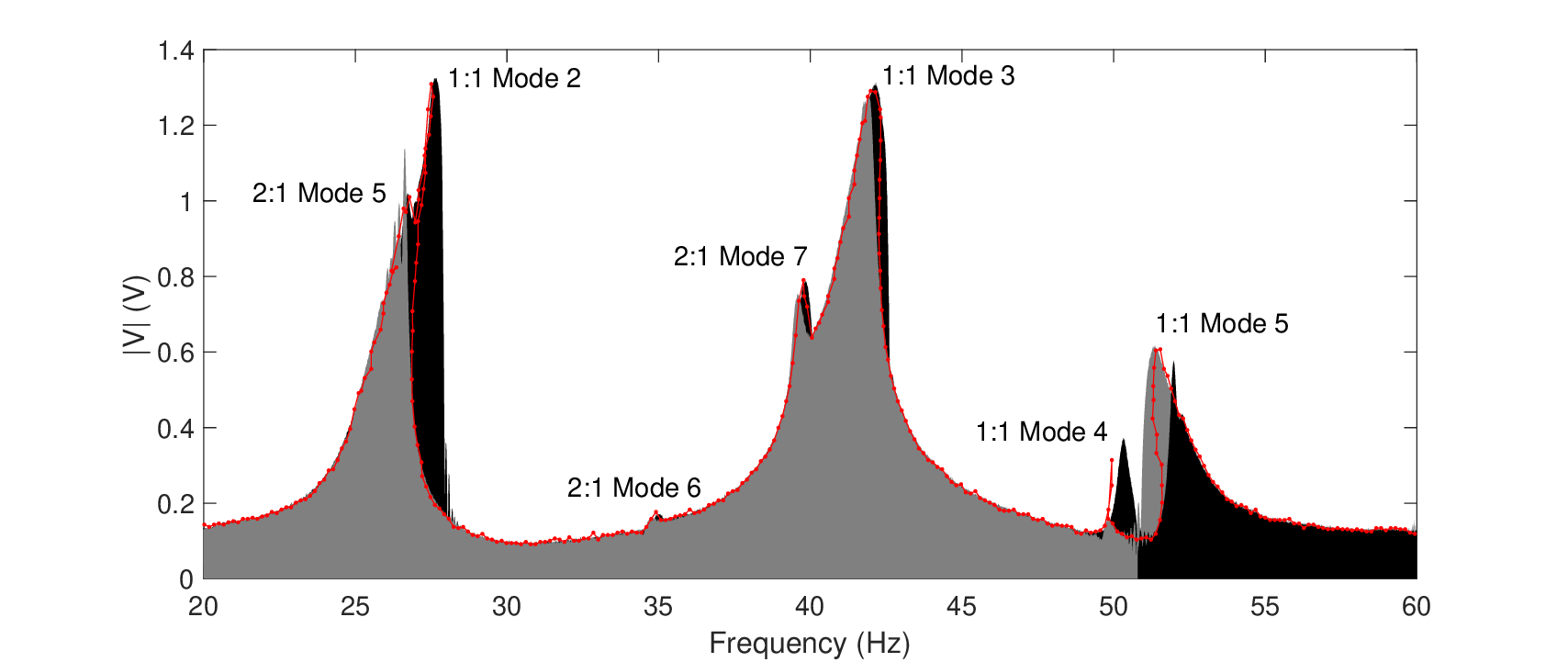}
            \caption{Shaker voltage NFR of the clamped plate for $V_\mathrm{in}=0.05$ V: SWS-up (\rule[.2em]{1em}{.2em}), SWS-down (\textcolor{gray}{\rule[.2em]{1em}{.2em}}) and ACBC (\textcolor{red}{\rule[.2em]{1em}{.2em}}).}
            \label{fig:PlateNFR_Voltage}
        \end{figure*}

        The ACBC method was tested at the highest level of excitation ($V_\mathrm{in} = $ 0.05 V) against SWS for reference. Figure~\ref{fig:PlateNFR_Voltage} presents the voltage NFR. The second and third modes are hardening, while the fifth mode is initially softening, but seems to start its hardening phase. The softening trend may come from a small but inevitable initial curvature of the plate. The results agree almost perfectly with SWS, except for the inevitable hysteresis, small transients near the jump phenomena, and for the fourth mode. The ACBC sometimes had trouble converging in the vicinity of the fourth mode. The issue is thought to come from its poor observability and controllability (this mode being almost absent from the considered driving-point NFR at a linear level, cf. Figure~\ref{fig:PlateFRF_Voltage}). Superharmonic resonances are appearing throughout the test. They are strong enough to be noticed as peaks in the NFR, but not enough to form loops; hence the ACBC method was able to characterize them.

        \begin{figure*}[!ht]
            \centering
            \includegraphics[width=\textwidth]{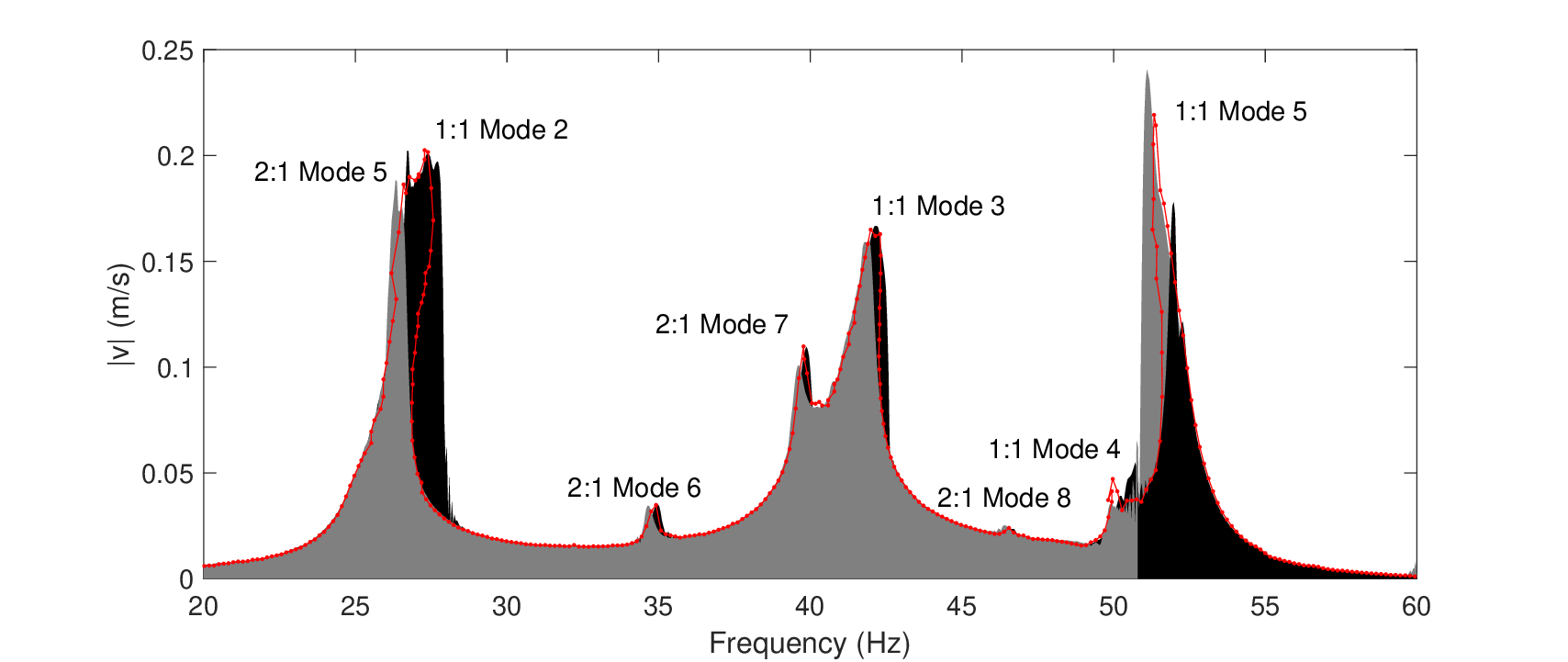}
            \caption{Edge velocity NFR of the clamped plate for $V_\mathrm{in}=0.05$ V: SWS-up (\rule[.2em]{1em}{.2em}), SWS-down (\textcolor{gray}{\rule[.2em]{1em}{.2em}}) and ACBC (\textcolor{red}{\rule[.2em]{1em}{.2em}}). }
            \label{fig:PlateNFR_Velocity}
        \end{figure*}

        Figure~\ref{fig:PlateNFR_Velocity} shows the velocity NFR at the laser measurement point. Mode 5 is more visible therein, as well as the superharmonic resonances.

        All modes feature an interesting and strongly nonlinear behavior at the highest considered forcing amplitude. To close this analysis, we investigate the dynamics of the second mode more in depth. It should be noted that since the frequency of mode 5 is about twice that of mode 2, the structure is in near-internal resonance conditions, resulting in an unconventional shape of the resonance peak.

    \subsection{Second mode}

        \begin{figure}[!ht]
                \centering
                \begin{subfigure}{.49\textwidth}
                    \centering
                    \includegraphics[width=1.1\textwidth]{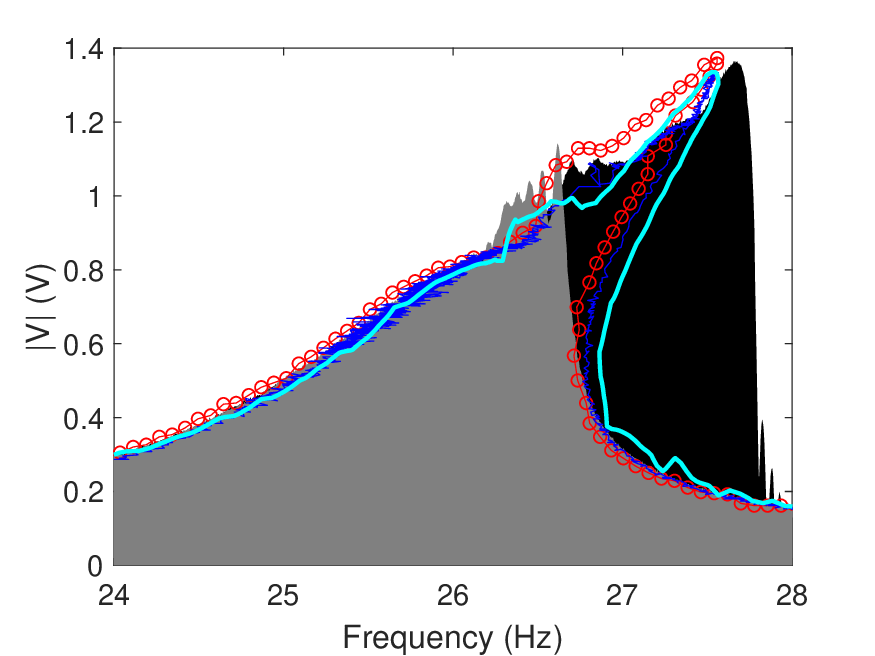}
                    \caption{}
                    \label{sfig:PlateMode2_NFR1}
                \end{subfigure}
                \begin{subfigure}{.49\textwidth}
                    \centering
                    \includegraphics[width=1.1\textwidth]{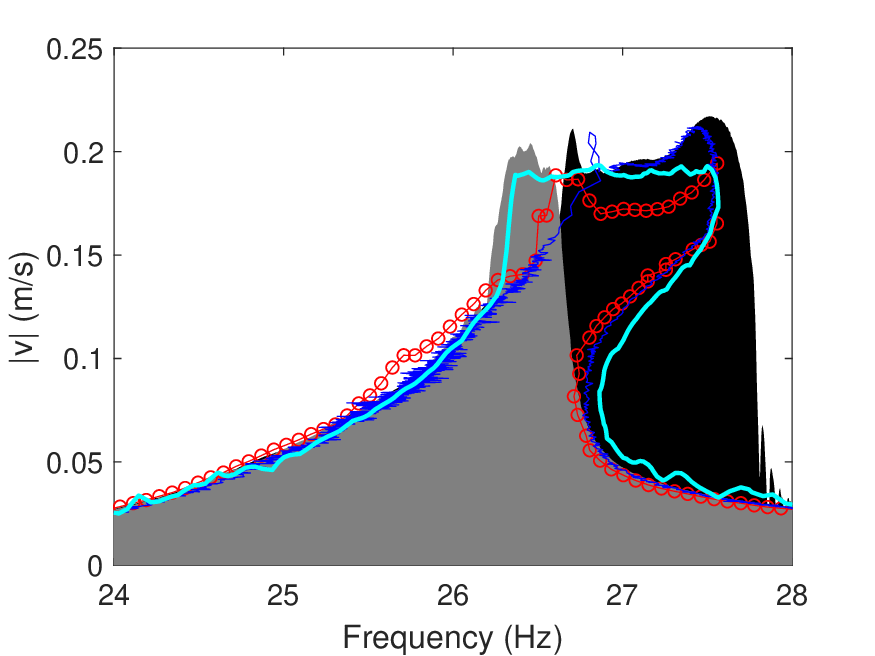}
                    \caption{}
                    \label{sfig:PlateMode2_NFR2}
                \end{subfigure}
                \caption{Shaker voltage~\subref{sfig:PlateMode2_NFR1} and laser velocity~\subref{sfig:PlateMode2_NFR2} NFRs around the second mode: SWS-up (\rule[.2em]{1em}{.2em}), SWS-down (\textcolor{gray}{\rule[.2em]{1em}{.2em}}), ACBC (\textcolor{red}{\rule[.2em]{1em}{.2em}}), PLL (\textcolor{blue}{\rule[.2em]{1em}{.2em}}) and SCBC (\textcolor{MatlabCyan}{\rule[.2em]{1em}{.2em}}).}
                \label{fig:PlateMode2_NFR}
            \end{figure}  

        We now focus on a single resonance, the second flexible mode of the assembly. Tests that were performed for the Duffing oscillator are repeated for the considered methods, namely SWS, SCBC, ACBC and PLL. All tests were performed on the same day within a time span of 3 h 30 min, and Table~\ref{tab:plateTimes} provides the time taken by each method (NB: as in the Duffing case, the times for SWS and PLL were chosen based on the fastest method). Figure~\ref{sfig:PlateMode2_NFR1} presents the results for the voltage NFR. All methods agree relatively well away from resonance, and in particular along the unstable branch (except of course for SWS), but small discrepancies are visible near resonance. These discrepancies are more marked for the velocity NFR, as illustrated by Figure~\ref{sfig:PlateMode2_NFR2}, where two peaks are generally visible.

        \begin{table}[!ht]
        \centering
        \caption{Times taken by each test for the primary resonance of the second mode of the plate.}
        \label{tab:plateTimes}
        \begin{tabular}{ccc}
             \hline
             Method & Total time (s) & Time per point (s) \\
             \hline
             SWS (up and down)& 190 & / \\
             SCBC (adaptive) & 397 & 0.58 \\
             ACBC & 190 & 1.13 \\
             PLL & 190 & / \\
             \hline
        \end{tabular}
        \end{table}

        \begin{figure}[!ht]
            \centering
                \includegraphics[width=0.539\textwidth]{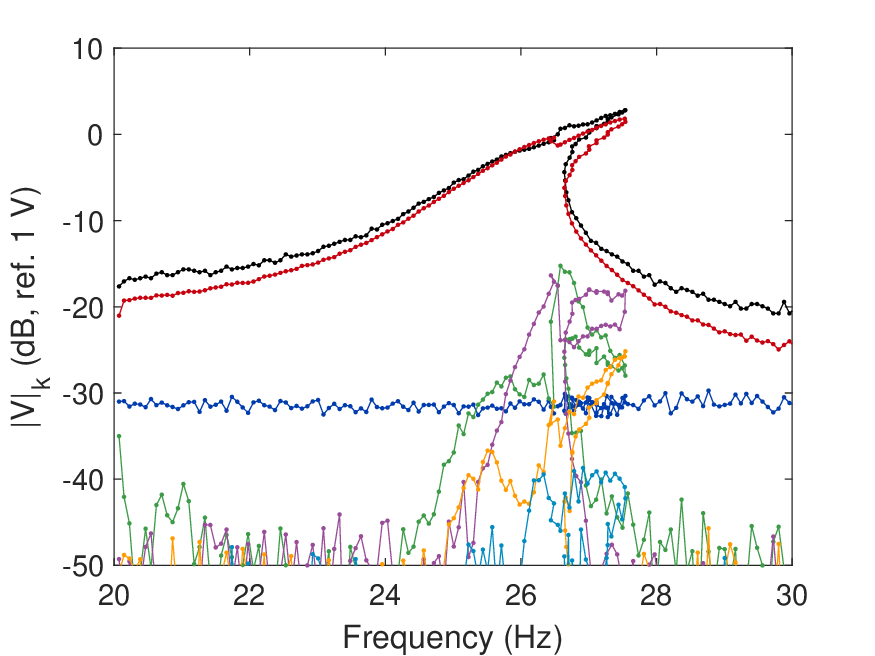}
            \caption{Harmonic amplitudes of the voltage signal: zeroth (\textcolor{MatlabCol1}{\rule[.2em]{1em}{.2em}}), first (\textcolor{MatlabCol2}{\rule[.2em]{1em}{.2em}}), second (\textcolor{MatlabCol3}{\rule[.2em]{1em}{.2em}}), third (\textcolor{MatlabCol4}{\rule[.2em]{1em}{.2em}}), fourth (\textcolor{MatlabCol5}{\rule[.2em]{1em}{.2em}}) and fifth (\textcolor{MatlabCol6}{\rule[.2em]{1em}{.2em}}) harmonic amplitudes, and total amplitude (\rule[.2em]{1em}{.2em}).}
            \label{fig:PlateMode2_Harmonics}
        \end{figure}

        Adaptive filters provide means to analyze the reason for the presence of a double peak in the NFR. As Figure~\ref{fig:PlateMode2_Harmonics} reveals, the second harmonic becomes more prominent around 26.5 Hz, which is slightly more than half the linear resonance frequency of mode 5. A 2:1 superharmonic resonance of mode 5 thus occurs right before the primary resonance of mode 2, explaining the appearance of the second peak.

        \begin{figure}[!ht]
            \centering
                \includegraphics[width=0.539\textwidth]{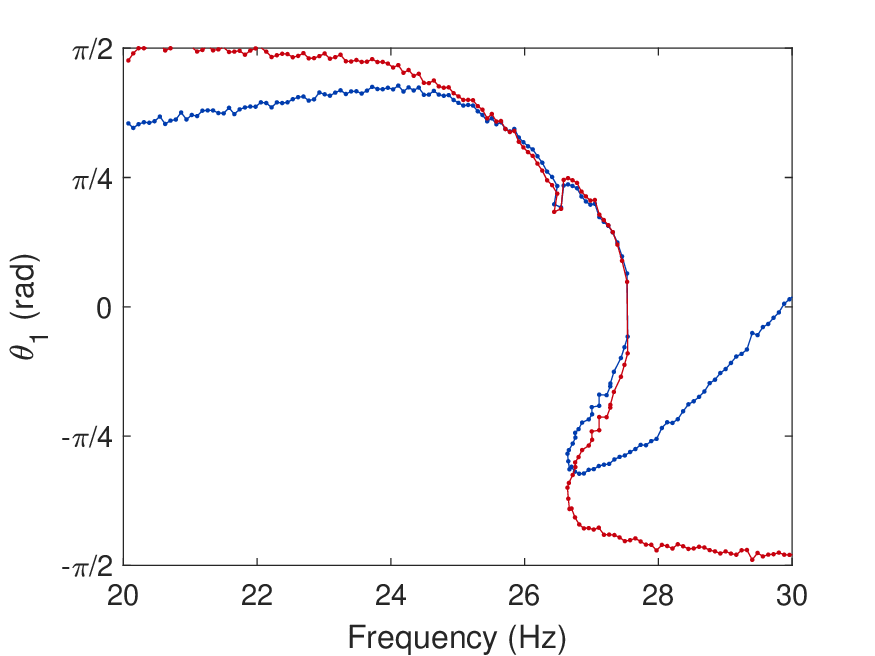}
            \caption{Phase of the first harmonic for the voltage (\textcolor{MatlabCol1}{\rule[.2em]{1em}{.2em}}) and velocity (\textcolor{MatlabCol2}{\rule[.2em]{1em}{.2em}}) NFRs.}
            \label{fig:PlateMode2_Phase}
        \end{figure}
        
        In addition, Figure~\ref{fig:PlateMode2_Phase} shows the evolution of the phase of the first harmonic for the voltage and velocity NFRs obtained with the ACBC. Similarly to the FRFs (Figures~\ref{fig:PlateFRF_Voltage} and~\ref{fig:PlateFRF_Velocity}), the voltage NFR does not feature a monotonous phase decrease around the resonance, whereas the velocity NFR shows this trend overall. Moreover, the 2:1 resonance manifests itself through a small phase trough, which is responsible for a small jump of the PLL method. With attention to detail, this jump can be observed in Figure~\ref{fig:PlateMode2_NFR}, approximately between 26.5 and 27 Hz, wherein the PLL curve has less resolution.

        \begin{figure}[!ht]
                \centering
                \begin{subfigure}{.49\textwidth}
                    \centering
                    \includegraphics[width=1.1\textwidth]{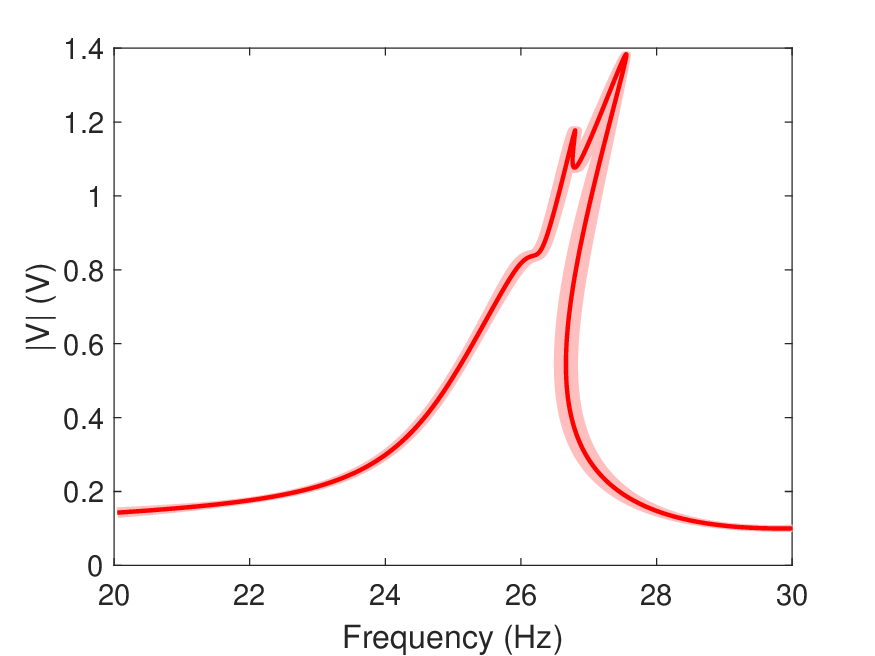}
                    \caption{}
                    \label{sfig:PlateMode2_SpreadACBC}
                \end{subfigure}
                \begin{subfigure}{.49\textwidth}
                    \centering
                    \includegraphics[width=1.1\textwidth]{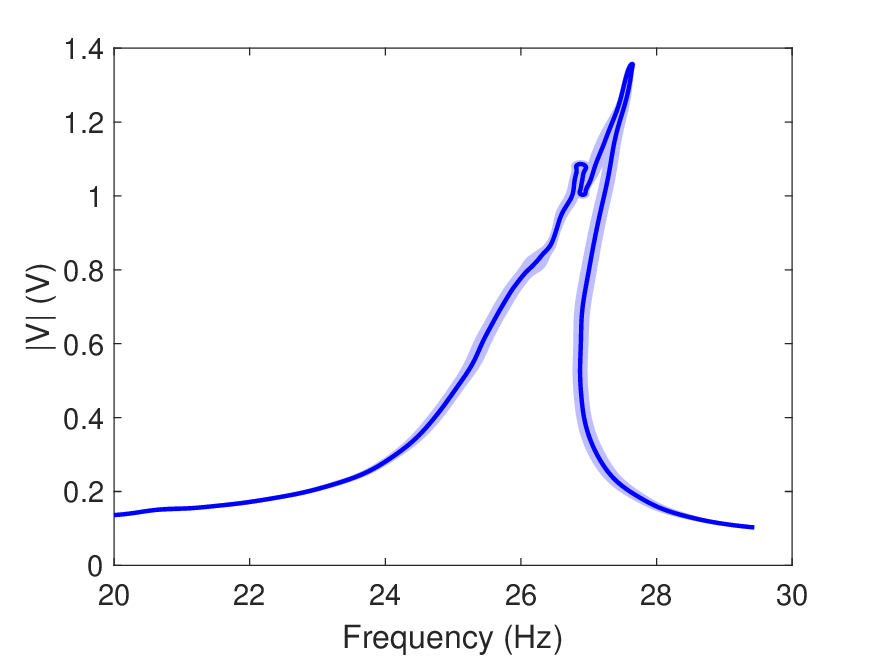}
                    \caption{}
                    \label{sfig:PlateMode2_SpreadPLL}
                \end{subfigure}
                \caption{Average NFR (dark curve), and 95\% confidence interval (light area) obtained from 10 repeatability tests: ACBC~\subref{sfig:PlateMode2_SpreadACBC} and PLL~\subref{sfig:PlateMode2_SpreadPLL}.}
                \label{fig:PlateMode2_Spread}
            \end{figure}  

        To assess the source of discrepancies, ten repeatability tests were performed for the ACBC and PLL. Statistics on the point cloud formed by the measurement were then obtained by fitting it with a spline with eleven control points and Gaussian process regressions~\cite{Goldberg1998}, as explained in the supplementary materials. Figure~\ref{fig:PlateMode2_Spread} presents the results for the voltage NFR. Except for an outlier for the ACBC, both types of test yield relatively repeatable results, although the differences between test partially explain the discrepancies observed in Figure~\ref{fig:PlateMode2_NFR}. The remaining reasons may come from the repeated tests at high amplitude (it can be deduced from Figure~\ref{sfig:PlateMode2_NFR2} that the plate undergoes a displacement amplitude of more than twice its thickness at the laser measurement point at resonance) over a relatively long time span. Indeed, there was at least 30 min of heavy testing in between each curve in Figure~\ref{fig:PlateMode2_NFR}. It is also interesting to note that the averaging process makes the small peak associated with the superharmonic resonance appear more clearly in both cases.

\section{Conclusion}
\label{sec:Conclusion}

With the growing need to exhaustively test nonlinear systems, experimental continuation approaches are emerging as robust and efficient means to obtain bifurcation diagrams experimentally, encompassing the open-loop unstable responses of the system under test. By reviewing their essential ingredients, namely continuation and control, this article presented the main variants in experimental continuation in a unified manner. In particular, a common framework for the derivative-free experimental continuation approaches which are usually presented in an individual way, namely SCBC, PLL and RCT, was given. It also objectively compared these different methods with an electronic Duffing oscillator and a thin clamped plate. The main conclusions coming from this work and the literature are outlined below, summarizing the strengths and weaknesses of the different methods in their current form.

SWS is typically the first type of test that one tries on a system. In spite of its inherent deficiencies, the fact that this method is open-loop makes it vastly simpler than the other methods presented in this work. It is an excellent approach to obtain preliminary results and to have an idea of the orders of magnitude involved with the system under test.

CBC offers the most mature and comprehensive suite for experimental continuation. With its numerous variants, it is possible to trade between simplicity and advanced nonlinear testing capabilities. The controller being  separated from the continuation algorithm, it is possible to have an intuitive understanding of its effects. Yet, CBC can be the most complex of the approaches outlined herein.  

PLL offers a quick and easy testing solution with a simple architecture. It is particularly effective and fast for testing around a specific resonance, obtaining the NFR close to this resonance, or tracking its backbone. Yet, tuning its parameters can sometimes prove difficult for the unstable parts of the NFR, and the fact that the phase folds as soon as there are multiple resonances (including secondary ones) inherently limits this approach for broadband testing. 

The RCT is an attractive method for the possibility to implement it with commercial software, and the simple postprocessing associated with tests made at constant vibration amplitude. There remains yet to be formal proofs of its stabilization properties (or the conditions for them), as examples herein showed that this is not guaranteed.

Finally, we look back at the proposed ACBC procedure. The fact that it is based on arclength continuation and derivative-free makes it a powerful and general tool, allowing the experimenter to obtain complex bifurcation diagrams, as featured in the plate example. Nevertheless, limitations were also highlighted, e.g., when the bifurcation diagram features a loop, or when adaptive filters fail to secure non-invasiveness. Future works will address these issues.

Since the pioneering work of Sieber and Krauskopf~\cite{SIEBER2007}, substantial progress has been made over the last twenty years in experimental continuation, but great challenges also lie ahead. Some of them are listed below:
\begin{enumerate}
    \item \textbf{Gain selection and adaptivity:} In almost all works, the control parameters are selected manually, often on a trial and error basis. A few theoretical studies have been made on a single-degree-of-freedom Duffing oscillator~\cite{Denis2018,Abeloos2022Thesis,Tatzko2023}, and a general but time-consuming parameter space exploration procedure was proposed in~\cite{Bureau2013}. More systematic procedures should be developed to enable the use of control-based methods for non-expert users. We note that recent works proposed to use adaptive procedures to control the system with solid theoretical properties~\cite{Li2021,Li2023,Rezaee2023}. These methods have stringent requirements for an experimental implementation, such as full-state feedback, but are very promising and represent a significant step forward. Alternatively, practical ways of setting the parameters of a PLL with limited information from the plant were proposed in~\cite{Hippold2024}. Finally, the data-driven approach based on parameter sweeps recently proposed in~\cite{Kruse2024} could form another interesting solution.
    \item \textbf{Stability and bifurcations:} When a test is successfully carried out, the closed-loop stability of every equilibria is ascertained. However, since the system of interest is the open-loop one, it is interesting to know whether the measured equilibria are open-loop stable or not. Procedures based on temporarily shutting the control down~\cite{Bureau2014} or on real-time linearization~\cite{Barton2017,Dittus2023} were proposed for the CBC, and similar procedures could be developed for the other control-based approaches. Going beyond that, different bifurcations could also be detected~\cite{Lim2011,Habib2023} and even tracked in an experimental setting, as in the pioneering works of Renson et al~\cite{Renson2017} and Shen et al~\cite{Shen2021b} (in the field of statics).  
    \item \textbf{Isolated solutions:} From a numerical standpoint, isolated solutions can be challenging to find with traditional continuation tools, but may be reachable in experiments. As a few works demonstrated~\cite{Bureau2014,Renson2019,Woiwode2024,zhou2024identification}, control-based methods could offer robust tools to trace out a complete isola, but reaching or even detecting them in a systematic way still represents a challenge.
    \item \textbf{Aperiodic responses:} None of the control-based method is currently able to deal with aperiodic (quasiperiodic or chaotic) responses, even though they may generally appear in nonlinear systems. As stated in Section~\ref{sec:NLControl}, it is possible to stabilize unstable periodic orbits embedded in a chaotic attractor~\cite{Ott1990,Pyragas1992}, but, to the authors' knowledge, little work has been done to stabilize an aperiodic attractor itself.
    \item \textbf{Actuator-system interaction:} In most cases, a structure will not be perfectly actuated. This problem pertains to the breadth of nonlinear vibration testing problems~\cite{Pacini2022} and not just to control-based methods. Although some attempts have been made~\cite{Magnevall2006,Novak2018}, no ideal solution currently exists in the literature.
    \item \textbf{Large-scale systems:} Except a few notable exceptions, control-based methods were almost always tested on systems of relatively low dimensionality, being either single-degree-of-freedom-like, or focusing on a single resonance widely spaced from the others. To achieve their full potential and provide convincing results for a future industrial use, control-based methods should turn toward more challenging, large-scale systems. The use of multiple sensors and/or actuators and the study of underactuated problems also remain largely uncovered.
\end{enumerate}

\backmatter












\begin{appendices}

\end{appendices}


\bibliography{bibliography}


\begin{thebibliography}{163}
\ifx \bisbn   \undefined \def \bisbn  #1{ISBN #1}\fi
\ifx \binits  \undefined \def \binits#1{#1}\fi
\ifx \bauthor  \undefined \def \bauthor#1{#1}\fi
\ifx \batitle  \undefined \def \batitle#1{#1}\fi
\ifx \bjtitle  \undefined \def \bjtitle#1{#1}\fi
\ifx \bvolume  \undefined \def \bvolume#1{\textbf{#1}}\fi
\ifx \byear  \undefined \def \byear#1{#1}\fi
\ifx \bissue  \undefined \def \bissue#1{#1}\fi
\ifx \bfpage  \undefined \def \bfpage#1{#1}\fi
\ifx \blpage  \undefined \def \blpage #1{#1}\fi
\ifx \burl  \undefined \def \burl#1{\textsf{#1}}\fi
\ifx \doiurl  \undefined \def \doiurl#1{\url{https://doi.org/#1}}\fi
\ifx \betal  \undefined \def \betal{\textit{et al.}}\fi
\ifx \binstitute  \undefined \def \binstitute#1{#1}\fi
\ifx \binstitutionaled  \undefined \def \binstitutionaled#1{#1}\fi
\ifx \bctitle  \undefined \def \bctitle#1{#1}\fi
\ifx \beditor  \undefined \def \beditor#1{#1}\fi
\ifx \bpublisher  \undefined \def \bpublisher#1{#1}\fi
\ifx \bbtitle  \undefined \def \bbtitle#1{#1}\fi
\ifx \bedition  \undefined \def \bedition#1{#1}\fi
\ifx \bseriesno  \undefined \def \bseriesno#1{#1}\fi
\ifx \blocation  \undefined \def \blocation#1{#1}\fi
\ifx \bsertitle  \undefined \def \bsertitle#1{#1}\fi
\ifx \bsnm \undefined \def \bsnm#1{#1}\fi
\ifx \bsuffix \undefined \def \bsuffix#1{#1}\fi
\ifx \bparticle \undefined \def \bparticle#1{#1}\fi
\ifx \barticle \undefined \def \barticle#1{#1}\fi
\bibcommenthead
\ifx \bconfdate \undefined \def \bconfdate #1{#1}\fi
\ifx \botherref \undefined \def \botherref #1{#1}\fi
\ifx \url \undefined \def \url#1{\textsf{#1}}\fi
\ifx \bchapter \undefined \def \bchapter#1{#1}\fi
\ifx \bbook \undefined \def \bbook#1{#1}\fi
\ifx \bcomment \undefined \def \bcomment#1{#1}\fi
\ifx \oauthor \undefined \def \oauthor#1{#1}\fi
\ifx \citeauthoryear \undefined \def \citeauthoryear#1{#1}\fi
\ifx \endbibitem  \undefined \def \endbibitem {}\fi
\ifx \bconflocation  \undefined \def \bconflocation#1{#1}\fi
\ifx \arxivurl  \undefined \def \arxivurl#1{\textsf{#1}}\fi
\csname PreBibitemsHook\endcsname

\bibitem[\protect\citeauthoryear{Nayfeh and Mook}{1995}]{Nayfeh1985}
\begin{bbook}
\bauthor{\bsnm{Nayfeh}, \binits{A.H.}},
\bauthor{\bsnm{Mook}, \binits{D.T.}}:
\bbtitle{Nonlinear Oscillations},
p. \bfpage{720}.
\bpublisher{Wiley},
\blocation{Weinheim}
(\byear{1995}).
\doiurl{10.1002/9783527617586} .
\burl{https://onlinelibrary.wiley.com/doi/book/10.1002/9783527617586}
\end{bbook}
\endbibitem

\bibitem[\protect\citeauthoryear{Guckenheimer and
  Holmes}{1983}]{Guckenheimer1984}
\begin{bbook}
\bauthor{\bsnm{Guckenheimer}, \binits{J.}},
\bauthor{\bsnm{Holmes}, \binits{P.}}:
\bbtitle{Nonlinear Oscillations, Dynamical Systems, and Bifurcations of Vector
  Fields}.
\bsertitle{Applied Mathematical Sciences},
vol. \bseriesno{42}.
\bpublisher{Springer},
\blocation{New York, NY}
(\byear{1983}).
\doiurl{10.1007/978-1-4612-1140-2} .
\burl{http://link.springer.com/10.1007/978-1-4612-1140-2}
\end{bbook}
\endbibitem

\bibitem[\protect\citeauthoryear{Kuznetsov}{2004}]{Kuznetsov2004}
\begin{bbook}
\bauthor{\bsnm{Kuznetsov}, \binits{Y.A.}}:
\bbtitle{Elements of Applied Bifurcation Theory}.
\bsertitle{Applied Mathematical Sciences},
vol. \bseriesno{112}.
\bpublisher{Springer},
\blocation{New York, NY}
(\byear{2004}).
\doiurl{10.1007/978-1-4757-3978-7} .
\burl{http://link.springer.com/10.1007/978-1-4757-3978-7}
\end{bbook}
\endbibitem

\bibitem[\protect\citeauthoryear{Wriggers}{2008}]{Wriggers}
\begin{bbook}
\bauthor{\bsnm{Wriggers}, \binits{P.}}:
\bbtitle{Nonlinear Finite Element Methods}.
\bpublisher{Springer},
\blocation{Berlin Heidelberg}
(\byear{2008})
\end{bbook}
\endbibitem

\bibitem[\protect\citeauthoryear{Keller}{1977}]{Keller1977}
\begin{bchapter}
\bauthor{\bsnm{Keller}, \binits{H.B.}}:
\bctitle{Numerical solution of bifurcation and nonlinear eigenvalue problems}.
In: \beditor{\bsnm{Rabinowitz}, \binits{P.H.}} (ed.)
\bbtitle{Applications of Bifurcation Theory: Proceedings of an Advanced
  Seminar},
pp. \bfpage{359}--\blpage{384}.
\bpublisher{Academic Press},
\blocation{New York, NY}
(\byear{1977})
\end{bchapter}
\endbibitem

\bibitem[\protect\citeauthoryear{Allgower and
  Georg}{2003}]{allgower2003introduction}
\begin{bbook}
\bauthor{\bsnm{Allgower}, \binits{E.L.}},
\bauthor{\bsnm{Georg}, \binits{K.}}:
\bbtitle{Introduction to Numerical Continuation Methods}.
\bpublisher{SIAM},
\blocation{Philadelphia}
(\byear{2003})
\end{bbook}
\endbibitem

\bibitem[\protect\citeauthoryear{Dhooge et~al.}{2003}]{Dhooge2003}
\begin{barticle}
\bauthor{\bsnm{Dhooge}, \binits{A.}},
\bauthor{\bsnm{Govaerts}, \binits{W.}},
\bauthor{\bsnm{Kuznetsov}, \binits{Y.A.}}:
\batitle{Matcont: A matlab package for numerical bifurcation analysis of odes}.
\bjtitle{ACM Transactions on Mathematical Software}
\bvolume{29}(\bissue{2}),
\bfpage{141}--\blpage{164}
(\byear{2003})
\doiurl{10.1145/779359.779362}
\end{barticle}
\endbibitem

\bibitem[\protect\citeauthoryear{Govaerts}{2000}]{govaerts2000numerical}
\begin{bbook}
\bauthor{\bsnm{Govaerts}, \binits{W.J.}}:
\bbtitle{Numerical Methods for Bifurcations of Dynamical Equilibria}.
\bpublisher{SIAM},
\blocation{Philadelphia}
(\byear{2000})
\end{bbook}
\endbibitem

\bibitem[\protect\citeauthoryear{Nayfeh and Balachandran}{1995}]{Nayfeh1995}
\begin{bbook}
\bauthor{\bsnm{Nayfeh}, \binits{A.H.}},
\bauthor{\bsnm{Balachandran}, \binits{B.}}:
\bbtitle{Applied Nonlinear Dynamics}.
\bpublisher{Wiley},
\blocation{Weinheim}
(\byear{1995}).
\doiurl{10.1002/9783527617548} .
\burl{https://onlinelibrary.wiley.com/doi/book/10.1002/9783527617548}
\end{bbook}
\endbibitem

\bibitem[\protect\citeauthoryear{Knowles et~al.}{2011}]{Knowles}
\begin{barticle}
\bauthor{\bsnm{Knowles}, \binits{J.}},
\bauthor{\bsnm{Krauskopf}, \binits{B.}},
\bauthor{\bsnm{Lowenberg}, \binits{M.}}:
\batitle{Numerical continuation applied to landing gear mechanism analysis}.
\bjtitle{Journal of Aircraft}
\bvolume{48}(\bissue{4}),
\bfpage{1254}--\blpage{1262}
(\byear{2011})
\doiurl{10.2514/1.C031247}
\end{barticle}
\endbibitem

\bibitem[\protect\citeauthoryear{Detroux et~al.}{2015}]{Detroux2015}
\begin{barticle}
\bauthor{\bsnm{Detroux}, \binits{T.}},
\bauthor{\bsnm{Renson}, \binits{L.}},
\bauthor{\bsnm{Masset}, \binits{L.}},
\bauthor{\bsnm{Kerschen}, \binits{G.}}:
\batitle{The harmonic balance method for bifurcation analysis of large-scale
  nonlinear mechanical systems}.
\bjtitle{Computer Methods in Applied Mechanics and Engineering}
\bvolume{296},
\bfpage{18}--\blpage{38}
(\byear{2015})
\doiurl{10.1016/j.cma.2015.07.017}
{\href{https://arxiv.org/abs/1604.05621}{{arXiv:1604.05621}}}
\end{barticle}
\endbibitem

\bibitem[\protect\citeauthoryear{Dankowicz and Schilder}{2011}]{Dankowicz2011}
\begin{botherref}
\oauthor{\bsnm{Dankowicz}, \binits{H.}},
\oauthor{\bsnm{Schilder}, \binits{F.}}:
An extended continuation problem for bifurcation analysis in the presence of
  constraints.
Journal of Computational and Nonlinear Dynamics
\textbf{6}(3)
(2011)
\doiurl{10.1115/1.4002684}
\end{botherref}
\endbibitem

\bibitem[\protect\citeauthoryear{Freour et~al.}{2020}]{Freour}
\begin{barticle}
\bauthor{\bsnm{Freour}, \binits{V.}},
\bauthor{\bsnm{Guillot}, \binits{L.}},
\bauthor{\bsnm{Masuda}, \binits{H.}},
\bauthor{\bsnm{Usa}, \binits{S.}},
\bauthor{\bsnm{Tominaga}, \binits{E.}},
\bauthor{\bsnm{Tohgi}, \binits{Y.}},
\bauthor{\bsnm{Vergez}, \binits{C.}},
\bauthor{\bsnm{Cochelin}, \binits{B.}}:
\batitle{Numerical continuation of a physical model of brass instruments:
  Application to trumpet comparisons}.
\bjtitle{Journal of the Acoustical SOciety of America}
\bvolume{148}(\bissue{2}),
\bfpage{748}
(\byear{2020})
\doiurl{10.1121/10.0001603}
\end{barticle}
\endbibitem

\bibitem[\protect\citeauthoryear{Haller and Ponsioen}{2016}]{Haller2016}
\begin{barticle}
\bauthor{\bsnm{Haller}, \binits{G.}},
\bauthor{\bsnm{Ponsioen}, \binits{S.}}:
\batitle{{Nonlinear normal modes and spectral submanifolds: existence,
  uniqueness and use in model reduction}}.
\bjtitle{Nonlinear Dynamics}
\bvolume{86}(\bissue{3}),
\bfpage{1493}--\blpage{1534}
(\byear{2016})
\doiurl{10.1007/s11071-016-2974-z}
{\href{https://arxiv.org/abs/1602.00560}{{arXiv:1602.00560}}}
\end{barticle}
\endbibitem

\bibitem[\protect\citeauthoryear{Vizzaccaro et~al.}{2021}]{Vizzaccaro2021}
\begin{botherref}
\oauthor{\bsnm{Vizzaccaro}, \binits{A.}},
\oauthor{\bsnm{Shen}, \binits{Y.}},
\oauthor{\bsnm{Salles}, \binits{L.}},
\oauthor{\bsnm{Blaho{\v{s}}}, \binits{J.}},
\oauthor{\bsnm{Touz{\'{e}}}, \binits{C.}}:
{Direct computation of nonlinear mapping via normal form for reduced-order
  models of finite element nonlinear structures}.
Computer Methods in Applied Mechanics and Engineering
\textbf{384}
(2021)
\doiurl{10.1016/j.cma.2021.113957}
{\href{https://arxiv.org/abs/2009.12145}{{arXiv:2009.12145}}}
\end{botherref}
\endbibitem

\bibitem[\protect\citeauthoryear{Jain and Haller}{2022}]{Jain2022}
\begin{barticle}
\bauthor{\bsnm{Jain}, \binits{S.}},
\bauthor{\bsnm{Haller}, \binits{G.}}:
\batitle{{How to compute invariant manifolds and their reduced dynamics in
  high-dimensional finite element models}}.
\bjtitle{Nonlinear Dynamics}
\bvolume{107}(\bissue{2}),
\bfpage{1417}--\blpage{1450}
(\byear{2022})
\doiurl{10.1007/s11071-021-06957-4}
{\href{https://arxiv.org/abs/2103.10264}{{arXiv:2103.10264}}}
\end{barticle}
\endbibitem

\bibitem[\protect\citeauthoryear{Touz{\'{e}} et~al.}{2021}]{Touze2021}
\begin{barticle}
\bauthor{\bsnm{Touz{\'{e}}}, \binits{C.}},
\bauthor{\bsnm{Vizzaccaro}, \binits{A.}},
\bauthor{\bsnm{Thomas}, \binits{O.}}:
\batitle{{Model order reduction methods for geometrically nonlinear structures:
  a review of nonlinear techniques}}.
\bjtitle{Nonlinear Dynamics}
\bvolume{105}(\bissue{2}),
\bfpage{1141}--\blpage{1190}
(\byear{2021})
\doiurl{10.1007/s11071-021-06693-9}
{\href{https://arxiv.org/abs/2107.05077}{{arXiv:2107.05077}}}
\end{barticle}
\endbibitem

\bibitem[\protect\citeauthoryear{Cenedese et~al.}{2022}]{Cenedese2022}
\begin{barticle}
\bauthor{\bsnm{Cenedese}, \binits{M.}},
\bauthor{\bsnm{Ax{\aa}s}, \binits{J.}},
\bauthor{\bsnm{B{\"{a}}uerlein}, \binits{B.}},
\bauthor{\bsnm{Avila}, \binits{K.}},
\bauthor{\bsnm{Haller}, \binits{G.}}:
\batitle{{Data-driven modeling and prediction of non-linearizable dynamics via
  spectral submanifolds}}.
\bjtitle{Nature Communications}
\bvolume{13}(\bissue{1}),
\bfpage{872}
(\byear{2022})
\doiurl{10.1038/s41467-022-28518-y}
\end{barticle}
\endbibitem

\bibitem[\protect\citeauthoryear{Virgin}{2000}]{Virgin}
\begin{bbook}
\bauthor{\bsnm{Virgin}, \binits{L.}}:
\bbtitle{Introduction to Experimental Nonlinear Dynamics}.
\bpublisher{Cambridge University Press},
\blocation{Cambridge}
(\byear{2000})
\end{bbook}
\endbibitem

\bibitem[\protect\citeauthoryear{Kerschen et~al.}{2006}]{Kerschen2006}
\begin{barticle}
\bauthor{\bsnm{Kerschen}, \binits{G.}},
\bauthor{\bsnm{Worden}, \binits{K.}},
\bauthor{\bsnm{Vakakis}, \binits{A.F.}},
\bauthor{\bsnm{Golinval}, \binits{J.-C.}}:
\batitle{Past, present and future of nonlinear system identification in
  structural dynamics}.
\bjtitle{Mechanical Systems and Signal Processing}
\bvolume{20}(\bissue{3}),
\bfpage{505}--\blpage{592}
(\byear{2006})
\doiurl{10.1016/j.ymssp.2005.04.008}
\end{barticle}
\endbibitem

\bibitem[\protect\citeauthoryear{No{\"{e}}l and Kerschen}{2017}]{Noel2017}
\begin{barticle}
\bauthor{\bsnm{No{\"{e}}l}, \binits{J.P.}},
\bauthor{\bsnm{Kerschen}, \binits{G.}}:
\batitle{Nonlinear system identification in structural dynamics: 10 more years
  of progress}.
\bjtitle{Mechanical Systems and Signal Processing}
\bvolume{83},
\bfpage{2}--\blpage{35}
(\byear{2017})
\doiurl{10.1016/j.ymssp.2016.07.020}
\end{barticle}
\endbibitem

\bibitem[\protect\citeauthoryear{Anastasio and
  Marchesiello}{2023}]{Anastasio2023}
\begin{barticle}
\bauthor{\bsnm{Anastasio}, \binits{D.}},
\bauthor{\bsnm{Marchesiello}, \binits{S.}}:
\batitle{{Nonlinear frequency response curves estimation and stability analysis
  of randomly excited systems in the subspace framework}}.
\bjtitle{Nonlinear Dynamics}
\bvolume{111}(\bissue{9}),
\bfpage{8115}--\blpage{8133}
(\byear{2023})
\doiurl{10.1007/s11071-023-08280-6}
\end{barticle}
\endbibitem

\bibitem[\protect\citeauthoryear{Sieber and Krauskopf}{2008}]{Sieber2008}
\begin{barticle}
\bauthor{\bsnm{Sieber}, \binits{J.}},
\bauthor{\bsnm{Krauskopf}, \binits{B.}}:
\batitle{Control based bifurcation analysis for experiments}.
\bjtitle{Nonlinear Dynamics}
\bvolume{51}(\bissue{3}),
\bfpage{365}--\blpage{377}
(\byear{2008})
\doiurl{10.1007/s11071-007-9217-2}
\end{barticle}
\endbibitem

\bibitem[\protect\citeauthoryear{Sieber et~al.}{2008}]{Sieber2008b}
\begin{barticle}
\bauthor{\bsnm{Sieber}, \binits{J.}},
\bauthor{\bsnm{Gonzalez-Buelga}, \binits{A.}},
\bauthor{\bsnm{Neild}, \binits{S.A.}},
\bauthor{\bsnm{Wagg}, \binits{D.J.}},
\bauthor{\bsnm{Krauskopf}, \binits{B.}}:
\batitle{Experimental continuation of periodic orbits through a fold}.
\bjtitle{Physical Review Letters}
\bvolume{100}(\bissue{24}),
\bfpage{1}--\blpage{4}
(\byear{2008})
\doiurl{10.1103/PhysRevLett.100.244101}
{\href{https://arxiv.org/abs/0804.0320}{{arXiv:0804.0320}}}
\end{barticle}
\endbibitem

\bibitem[\protect\citeauthoryear{Barton and Sieber}{2013}]{Barton2013}
\begin{barticle}
\bauthor{\bsnm{Barton}, \binits{D.A.W.}},
\bauthor{\bsnm{Sieber}, \binits{J.}}:
\batitle{Systematic experimental exploration of bifurcations with noninvasive
  control}.
\bjtitle{Physical Review E}
\bvolume{87}(\bissue{5}),
\bfpage{052916}
(\byear{2013})
\doiurl{10.1103/PhysRevE.87.052916}
\end{barticle}
\endbibitem

\bibitem[\protect\citeauthoryear{Peter and Leine}{2017}]{Peter2017}
\begin{barticle}
\bauthor{\bsnm{Peter}, \binits{S.}},
\bauthor{\bsnm{Leine}, \binits{R.I.}}:
\batitle{Excitation power quantities in phase resonance testing of nonlinear
  systems with phase-locked-loop excitation}.
\bjtitle{Mechanical Systems and Signal Processing}
\bvolume{96},
\bfpage{139}--\blpage{158}
(\byear{2017})
\doiurl{10.1016/j.ymssp.2017.04.011}
\end{barticle}
\endbibitem

\bibitem[\protect\citeauthoryear{Karaağa{\c{c}}lı and
  {\"{O}}zg{\"{u}}ven}{2020}]{Karaagacl2020}
\begin{barticle}
\bauthor{\bsnm{Karaağa{\c{c}}lı}, \binits{T.}},
\bauthor{\bsnm{{\"{O}}zg{\"{u}}ven}, \binits{H.N.}}:
\batitle{Experimental identification of backbone curves of strongly nonlinear
  systems by using response-controlled stepped-sine testing (rct)}.
\bjtitle{Vibration}
\bvolume{3}(\bissue{3}),
\bfpage{266}--\blpage{280}
(\byear{2020})
\doiurl{10.3390/vibration3030019}
\end{barticle}
\endbibitem

\bibitem[\protect\citeauthoryear{Abeloos}{2022}]{Abeloos2022Thesis}
\begin{botherref}
\oauthor{\bsnm{Abeloos}, \binits{G.}}:
Control-based methods for the identification of nonlinear structures.
PhD thesis,
University of Li{\`{e}}ge
(2022).
\url{https://hdl.handle.net/2268/295414}
\end{botherref}
\endbibitem

\bibitem[\protect\citeauthoryear{Kovacic and Brennan}{2011}]{Kovacic2011}
\begin{bbook}
\bauthor{\bsnm{Kovacic}, \binits{I.}},
\bauthor{\bsnm{Brennan}, \binits{M.J.}}:
\bbtitle{The {D}uffing Equation}.
\bpublisher{John Wiley {\&} Sons, Ltd},
\blocation{Chichester, UK}
(\byear{2011}).
\doiurl{10.1002/9780470977859} .
\burl{http://doi.wiley.com/10.1002/9780470977859}
\end{bbook}
\endbibitem

\bibitem[\protect\citeauthoryear{Seydel}{2009}]{seydel2009practical}
\begin{bbook}
\bauthor{\bsnm{Seydel}, \binits{R.}}:
\bbtitle{Practical Bifurcation and Stability Analysis}
vol. \bseriesno{5}.
\bpublisher{Springer},
\blocation{New York Dordrecht Heidelberg London}
(\byear{2009})
\end{bbook}
\endbibitem

\bibitem[\protect\citeauthoryear{Krack and Gross}{2019}]{Krack2019}
\begin{bbook}
\bauthor{\bsnm{Krack}, \binits{M.}},
\bauthor{\bsnm{Gross}, \binits{J.}}:
\bbtitle{Harmonic Balance for Nonlinear Vibration Problems}.
\bsertitle{Mathematical Engineering},
p. \bfpage{159}.
\bpublisher{Springer},
\blocation{Cham}
(\byear{2019}).
\doiurl{10.1007/978-3-030-14023-6} .
\burl{http://link.springer.com/10.1007/978-3-030-14023-6}
\end{bbook}
\endbibitem

\bibitem[\protect\citeauthoryear{Dankowicz and
  Schilder}{2013}]{dankowicz2013recipes}
\begin{bbook}
\bauthor{\bsnm{Dankowicz}, \binits{H.}},
\bauthor{\bsnm{Schilder}, \binits{F.}}:
\bbtitle{Recipes for Continuation}.
\bpublisher{SIAM},
\blocation{Philadelphia}
(\byear{2013})
\end{bbook}
\endbibitem

\bibitem[\protect\citeauthoryear{Doedel et~al.}{1991}]{DOEDEL1991}
\begin{barticle}
\bauthor{\bsnm{Doedel}, \binits{E.}},
\bauthor{\bsnm{Keller}, \binits{H.B.}},
\bauthor{\bsnm{Kernevez}, \binits{J.P.}}:
\batitle{Numerical analysis and control of bifurcation problems (i):
  Bifurcation in finite dimensions}.
\bjtitle{International Journal of Bifurcation and Chaos}
\bvolume{01}(\bissue{03}),
\bfpage{493}--\blpage{520}
(\byear{1991})
\doiurl{10.1142/S0218127491000397}
\end{barticle}
\endbibitem

\bibitem[\protect\citeauthoryear{Cochelin and Vergez}{2009}]{Cochelin2009}
\begin{barticle}
\bauthor{\bsnm{Cochelin}, \binits{B.}},
\bauthor{\bsnm{Vergez}, \binits{C.}}:
\batitle{A high order purely frequency-based harmonic balance formulation for
  continuation of periodic solutions}.
\bjtitle{Journal of Sound and Vibration}
\bvolume{324}(\bissue{1-2}),
\bfpage{243}--\blpage{262}
(\byear{2009})
\doiurl{10.1016/j.jsv.2009.01.054}
\end{barticle}
\endbibitem

\bibitem[\protect\citeauthoryear{Krauskopf
  et~al.}{2007}]{krauskopf2007numerical}
\begin{bbook}
\bauthor{\bsnm{Krauskopf}, \binits{B.}},
\bauthor{\bsnm{Osinga}, \binits{H.M.}},
\bauthor{\bsnm{Gal{\'a}n-Vioque}, \binits{J.}}:
\bbtitle{Numerical Continuation Methods for Dynamical Systems}
vol. \bseriesno{2}.
\bpublisher{Springer},
\blocation{Dordrecht}
(\byear{2007})
\end{bbook}
\endbibitem

\bibitem[\protect\citeauthoryear{Haselgrove}{1961}]{Haselgrove1961}
\begin{barticle}
\bauthor{\bsnm{Haselgrove}, \binits{C.B.}}:
\batitle{The solution of non-linear equations and of differential equations
  with two-point boundary conditions}.
\bjtitle{The Computer Journal}
\bvolume{4}(\bissue{3}),
\bfpage{255}--\blpage{259}
(\byear{1961})
\doiurl{10.1093/comjnl/4.3.255}
\end{barticle}
\endbibitem

\bibitem[\protect\citeauthoryear{Riks}{1979}]{Riks1979}
\begin{barticle}
\bauthor{\bsnm{Riks}, \binits{E.}}:
\batitle{An incremental approach to the solution of snapping and buckling
  problems}.
\bjtitle{International Journal of Solids and Structures}
\bvolume{15}(\bissue{7}),
\bfpage{529}--\blpage{551}
(\byear{1979})
\doiurl{10.1016/0020-7683(79)90081-7}
\end{barticle}
\endbibitem

\bibitem[\protect\citeauthoryear{Crisfield}{1981}]{CRISFIELD1981}
\begin{barticle}
\bauthor{\bsnm{Crisfield}, \binits{M.A.}}:
\batitle{{A fast incremental/iterative solution procedure that handles
  “snap-through”}}.
\bjtitle{Computers {\&} Structures}
\bvolume{13}(\bissue{1-3}),
\bfpage{55}--\blpage{62}
(\byear{1981})
\doiurl{10.1016/0045-7949(81)90108-5}
\end{barticle}
\endbibitem

\bibitem[\protect\citeauthoryear{Lewandowski}{1992}]{Lewandowski1992}
\begin{barticle}
\bauthor{\bsnm{Lewandowski}, \binits{R.}}:
\batitle{Non-linear, steady-state vibration of structures by harmonic
  balance/finite element method}.
\bjtitle{Computers {\&} Structures}
\bvolume{44}(\bissue{1-2}),
\bfpage{287}--\blpage{296}
(\byear{1992})
\doiurl{10.1016/0045-7949(92)90248-X}
\end{barticle}
\endbibitem

\bibitem[\protect\citeauthoryear{Sundararajan and
  Noah}{1997}]{Sundararajan1997}
\begin{barticle}
\bauthor{\bsnm{Sundararajan}, \binits{P.}},
\bauthor{\bsnm{Noah}, \binits{S.T.}}:
\batitle{Dynamics of forced nonlinear systems using shooting/are-length
  continuation method-application to rotor systems}.
\bjtitle{Journal of Vibration and Acoustics, Transactions of the ASME}
\bvolume{119}(\bissue{1}),
\bfpage{9}--\blpage{20}
(\byear{1997})
\doiurl{10.1115/1.2889694}
\end{barticle}
\endbibitem

\bibitem[\protect\citeauthoryear{Slotine}{1991}]{slotine1991applied}
\begin{bbook}
\bauthor{\bsnm{Slotine}, \binits{J.-J.E.}}:
\bbtitle{Applied Nonlinear Control}.
\bpublisher{Prentice-Hall},
\blocation{Englewood Cliffs}
(\byear{1991})
\end{bbook}
\endbibitem

\bibitem[\protect\citeauthoryear{Parlitz and Lauterborn}{1985}]{Parlitz1985}
\begin{barticle}
\bauthor{\bsnm{Parlitz}, \binits{U.}},
\bauthor{\bsnm{Lauterborn}, \binits{W.}}:
\batitle{Superstructure in the bifurcation set of the {D}uffing equation}.
\bjtitle{Physics Letters A}
\bvolume{107}(\bissue{8}),
\bfpage{351}--\blpage{355}
(\byear{1985})
\doiurl{10.1016/0375-9601(85)90687-5}
\end{barticle}
\endbibitem

\bibitem[\protect\citeauthoryear{Marchionne et~al.}{2018}]{Marchionne2018}
\begin{barticle}
\bauthor{\bsnm{Marchionne}, \binits{A.}},
\bauthor{\bsnm{Ditlevsen}, \binits{P.}},
\bauthor{\bsnm{Wieczorek}, \binits{S.}}:
\batitle{Synchronisation vs. resonance: Isolated resonances in damped nonlinear
  oscillators}.
\bjtitle{Physica D: Nonlinear Phenomena}
\bvolume{380-381},
\bfpage{8}--\blpage{16}
(\byear{2018})
\doiurl{10.1016/j.physd.2018.05.004}
\end{barticle}
\endbibitem

\bibitem[\protect\citeauthoryear{Stoker}{1950}]{stoker1950nonlinear}
\begin{bbook}
\bauthor{\bsnm{Stoker}, \binits{J.J.}}:
\bbtitle{Nonlinear Vibrations in Mechanical and Electrical Systems}
vol. \bseriesno{2}.
\bpublisher{Interscience publishers},
\blocation{New York}
(\byear{1950})
\end{bbook}
\endbibitem

\bibitem[\protect\citeauthoryear{Pacini et~al.}{2022}]{Pacini2022}
\begin{barticle}
\bauthor{\bsnm{Pacini}, \binits{B.R.}},
\bauthor{\bsnm{Kuether}, \binits{R.J.}},
\bauthor{\bsnm{Roettgen}, \binits{D.R.}}:
\batitle{Shaker-structure interaction modeling and analysis for nonlinear force
  appropriation testing}.
\bjtitle{Mechanical Systems and Signal Processing}
\bvolume{162}(\bissue{May 2021}),
\bfpage{108000}
(\byear{2022})
\doiurl{10.1016/j.ymssp.2021.108000}
\end{barticle}
\endbibitem

\bibitem[\protect\citeauthoryear{Virgin et~al.}{1998}]{Virgin1998}
\begin{barticle}
\bauthor{\bsnm{Virgin}, \binits{L.N.}},
\bauthor{\bsnm{Todd}, \binits{M.D.}},
\bauthor{\bsnm{Begley}, \binits{C.J.}},
\bauthor{\bsnm{Trickey}, \binits{S.T.}},
\bauthor{\bsnm{Dowell}, \binits{E.H.}}:
\batitle{Basins of attraction in experimental nonlinear oscillators}.
\bjtitle{International Journal of Bifurcation and Chaos}
\bvolume{08}(\bissue{03}),
\bfpage{521}--\blpage{533}
(\byear{1998})
\doiurl{10.1142/S0218127498000334}
\end{barticle}
\endbibitem

\bibitem[\protect\citeauthoryear{Ludeke}{1951}]{Ludeke1951}
\begin{barticle}
\bauthor{\bsnm{Ludeke}, \binits{C.A.}}:
\batitle{Predominantly subharmonic oscillations}.
\bjtitle{Journal of Applied Physics}
\bvolume{22}(\bissue{11}),
\bfpage{1321}--\blpage{1326}
(\byear{1951})
\doiurl{10.1063/1.1699858}
\end{barticle}
\endbibitem

\bibitem[\protect\citeauthoryear{Sieber and Krauskopf}{2007}]{SIEBER2007}
\begin{barticle}
\bauthor{\bsnm{Sieber}, \binits{J.}},
\bauthor{\bsnm{Krauskopf}, \binits{B.}}:
\batitle{Control-based continuation of periodic orbits with a time-delayed
  difference scheme}.
\bjtitle{International Journal of Bifurcation and Chaos}
\bvolume{17}(\bissue{08}),
\bfpage{2579}--\blpage{2593}
(\byear{2007})
\doiurl{10.1142/S0218127407018646}
\end{barticle}
\endbibitem

\bibitem[\protect\citeauthoryear{Mayr}{1970}]{Mayr1970}
\begin{barticle}
\bauthor{\bsnm{Mayr}, \binits{O.}}:
\batitle{The origins of feedback control}.
\bjtitle{Scientific American}
\bvolume{223}(\bissue{4}),
\bfpage{110}--\blpage{118}
(\byear{1970})
\doiurl{10.1038/scientificamerican1070-110}
\end{barticle}
\endbibitem

\bibitem[\protect\citeauthoryear{Fuller}{1963}]{Fuller1963}
\begin{barticle}
\bauthor{\bsnm{Fuller}, \binits{A.T.}}:
\batitle{Directions of research in control}.
\bjtitle{Automatica}
\bvolume{1}(\bissue{4}),
\bfpage{289}--\blpage{296}
(\byear{1963})
\doiurl{10.1016/0005-1098(63)90013-X}
\end{barticle}
\endbibitem

\bibitem[\protect\citeauthoryear{Bennett}{1996}]{Bennett1996}
\begin{barticle}
\bauthor{\bsnm{Bennett}, \binits{S.}}:
\batitle{A brief history of automatic control}.
\bjtitle{IEEE Control Systems}
\bvolume{16}(\bissue{3}),
\bfpage{17}--\blpage{25}
(\byear{1996})
\doiurl{10.1109/37.506394}
\end{barticle}
\endbibitem

\bibitem[\protect\citeauthoryear{Fuller}{1976a}]{Fuller1976a}
\begin{barticle}
\bauthor{\bsnm{Fuller}, \binits{A.T.}}:
\batitle{The early development of control theory}.
\bjtitle{Journal of Dynamic Systems, Measurement, and Control}
\bvolume{98}(\bissue{2}),
\bfpage{109}--\blpage{118}
(\byear{1976})
\doiurl{10.1115/1.3426994}
\end{barticle}
\endbibitem

\bibitem[\protect\citeauthoryear{Fuller}{1976b}]{Fuller1976}
\begin{barticle}
\bauthor{\bsnm{Fuller}, \binits{A.T.}}:
\batitle{The early development of control theory. ii}.
\bjtitle{Journal of Dynamic Systems, Measurement, and Control}
\bvolume{98}(\bissue{3}),
\bfpage{224}--\blpage{235}
(\byear{1976})
\doiurl{10.1115/1.3427026}
\end{barticle}
\endbibitem

\bibitem[\protect\citeauthoryear{Atherton}{1996}]{Atherton1996}
\begin{barticle}
\bauthor{\bsnm{Atherton}, \binits{D.P.}}:
\batitle{Early developments in nonlinear control}.
\bjtitle{IEEE Control Systems}
\bvolume{16}(\bissue{3}),
\bfpage{34}--\blpage{43}
(\byear{1996})
\doiurl{10.1109/37.506396}
\end{barticle}
\endbibitem

\bibitem[\protect\citeauthoryear{Kokotovi{\'{c}} and
  Arcak}{2001}]{Kokotovic2001}
\begin{barticle}
\bauthor{\bsnm{Kokotovi{\'{c}}}, \binits{P.}},
\bauthor{\bsnm{Arcak}, \binits{M.}}:
\batitle{Constructive nonlinear control: a historical perspective}.
\bjtitle{Automatica}
\bvolume{37}(\bissue{5}),
\bfpage{637}--\blpage{662}
(\byear{2001})
\doiurl{10.1016/S0005-1098(01)00002-4}
\end{barticle}
\endbibitem

\bibitem[\protect\citeauthoryear{Iqbal et~al.}{2017}]{Iqbal2017}
\begin{barticle}
\bauthor{\bsnm{Iqbal}, \binits{J.}},
\bauthor{\bsnm{Ullah}, \binits{M.}},
\bauthor{\bsnm{Khan}, \binits{S.G.}},
\bauthor{\bsnm{Khelifa}, \binits{B.}},
\bauthor{\bsnm{{\'{C}}ukovi{\'{c}}}, \binits{S.}}:
\batitle{Nonlinear control systems - a brief overview of historical and recent
  advances}.
\bjtitle{Nonlinear Engineering}
\bvolume{6}(\bissue{4}),
\bfpage{301}--\blpage{312}
(\byear{2017})
\doiurl{10.1515/nleng-2016-0077}
\end{barticle}
\endbibitem

\bibitem[\protect\citeauthoryear{Khalil}{2002}]{khalil2002control}
\begin{bbook}
\bauthor{\bsnm{Khalil}, \binits{H.K.}}:
\bbtitle{Nonlinear Control}.
\bpublisher{Pearson Education},
\blocation{Harlow}
(\byear{2002})
\end{bbook}
\endbibitem

\bibitem[\protect\citeauthoryear{Sepulchre
  et~al.}{2012}]{sepulchre2012constructive}
\begin{bbook}
\bauthor{\bsnm{Sepulchre}, \binits{R.}},
\bauthor{\bsnm{Jankovic}, \binits{M.}},
\bauthor{\bsnm{Kokotovic}, \binits{P.V.}}:
\bbtitle{Constructive Nonlinear Control}.
\bpublisher{Springer},
\blocation{London}
(\byear{2012})
\end{bbook}
\endbibitem

\bibitem[\protect\citeauthoryear{Wagg and Neild}{2015}]{Wagg2015}
\begin{bbook}
\bauthor{\bsnm{Wagg}, \binits{D.}},
\bauthor{\bsnm{Neild}, \binits{S.}}:
\bbtitle{Nonlinear Vibration with Control}.
\bsertitle{Solid Mechanics and Its Applications},
vol. \bseriesno{218}.
\bpublisher{Springer},
\blocation{Cham}
(\byear{2015}).
\doiurl{10.1007/978-3-319-10644-1} .
\burl{https://link.springer.com/10.1007/978-3-319-10644-1}
\end{bbook}
\endbibitem

\bibitem[\protect\citeauthoryear{Ott et~al.}{1990}]{Ott1990}
\begin{barticle}
\bauthor{\bsnm{Ott}, \binits{E.}},
\bauthor{\bsnm{Grebogi}, \binits{C.}},
\bauthor{\bsnm{Yorke}, \binits{J.A.}}:
\batitle{Controlling chaos}.
\bjtitle{Physical Review Letters}
\bvolume{64}(\bissue{11}),
\bfpage{1196}--\blpage{1199}
(\byear{1990})
\doiurl{10.1103/PhysRevLett.64.1196}
\end{barticle}
\endbibitem

\bibitem[\protect\citeauthoryear{Shinbrot et~al.}{1993}]{Shinbrot1993}
\begin{barticle}
\bauthor{\bsnm{Shinbrot}, \binits{T.}},
\bauthor{\bsnm{Grebogi}, \binits{C.}},
\bauthor{\bsnm{Yorke}, \binits{J.A.}},
\bauthor{\bsnm{Ott}, \binits{E.}}:
\batitle{Using small perturbations to control chaos}.
\bjtitle{Nature}
\bvolume{363}(\bissue{6428}),
\bfpage{411}--\blpage{417}
(\byear{1993})
\doiurl{10.1038/363411a0}
\end{barticle}
\endbibitem

\bibitem[\protect\citeauthoryear{Boccaletti et~al.}{2000}]{Boccaletti2000}
\begin{barticle}
\bauthor{\bsnm{Boccaletti}, \binits{S.}},
\bauthor{\bsnm{Grebogi}, \binits{C.}},
\bauthor{\bsnm{Lai}, \binits{Y.C.}},
\bauthor{\bsnm{Mancini}, \binits{H.}},
\bauthor{\bsnm{Maza}, \binits{D.}}:
\batitle{The control of chaos: Theory and applications}.
\bjtitle{Physics Report}
\bvolume{329}(\bissue{3}),
\bfpage{103}--\blpage{197}
(\byear{2000})
\doiurl{10.1016/S0370-1573(99)00096-4}
\end{barticle}
\endbibitem

\bibitem[\protect\citeauthoryear{Fradkov and Evans}{2005}]{Fradkov2005}
\begin{barticle}
\bauthor{\bsnm{Fradkov}, \binits{A.L.}},
\bauthor{\bsnm{Evans}, \binits{R.J.}}:
\batitle{Control of chaos: Methods and applications in engineering}.
\bjtitle{Annual Reviews in Control}
\bvolume{29}(\bissue{1}),
\bfpage{33}--\blpage{56}
(\byear{2005})
\doiurl{10.1016/j.arcontrol.2005.01.001}
\end{barticle}
\endbibitem

\bibitem[\protect\citeauthoryear{Pyragas}{1992}]{Pyragas1992}
\begin{barticle}
\bauthor{\bsnm{Pyragas}, \binits{K.}}:
\batitle{Continuous control of chaos by self-controlling feedback}.
\bjtitle{Physics Letters A}
\bvolume{170}(\bissue{6}),
\bfpage{421}--\blpage{428}
(\byear{1992})
\doiurl{10.1016/0375-9601(92)90745-8}
\end{barticle}
\endbibitem

\bibitem[\protect\citeauthoryear{Carroll et~al.}{1992}]{Carroll1992}
\begin{barticle}
\bauthor{\bsnm{Carroll}, \binits{T.L.}},
\bauthor{\bsnm{Triandaf}, \binits{I.}},
\bauthor{\bsnm{Schwartz}, \binits{I.}},
\bauthor{\bsnm{Pecora}, \binits{L.}}:
\batitle{Tracking unstable orbits in an experiment}.
\bjtitle{Physical Review A}
\bvolume{46}(\bissue{10}),
\bfpage{6189}--\blpage{6192}
(\byear{1992})
\doiurl{10.1103/PhysRevA.46.6189}
\end{barticle}
\endbibitem

\bibitem[\protect\citeauthoryear{Gills et~al.}{1992}]{Gills1992}
\begin{barticle}
\bauthor{\bsnm{Gills}, \binits{Z.}},
\bauthor{\bsnm{Iwata}, \binits{C.}},
\bauthor{\bsnm{Roy}, \binits{R.}},
\bauthor{\bsnm{Schwartz}, \binits{I.B.}},
\bauthor{\bsnm{Triandaf}, \binits{I.}}:
\batitle{Tracking unstable steady states: Extending the stability regime of a
  multimode laser system}.
\bjtitle{Physical Review Letters}
\bvolume{69}(\bissue{22}),
\bfpage{3169}--\blpage{3172}
(\byear{1992})
\doiurl{10.1103/PhysRevLett.69.3169}
\end{barticle}
\endbibitem

\bibitem[\protect\citeauthoryear{Neville et~al.}{2020}]{Neville2020}
\begin{barticle}
\bauthor{\bsnm{Neville}, \binits{R.M.}},
\bauthor{\bsnm{Groh}, \binits{R.M.J.}},
\bauthor{\bsnm{Pirrera}, \binits{A.}},
\bauthor{\bsnm{Schenk}, \binits{M.}}:
\batitle{Beyond the fold: experimentally traversing limit points in nonlinear
  structures}.
\bjtitle{Proceedings of the Royal Society A: Mathematical, Physical and
  Engineering Sciences}
\bvolume{476}(\bissue{2233}),
\bfpage{20190576}
(\byear{2020})
\doiurl{10.1098/rspa.2019.0576}
\end{barticle}
\endbibitem

\bibitem[\protect\citeauthoryear{Tatzko et~al.}{2023}]{Tatzko2023}
\begin{barticle}
\bauthor{\bsnm{Tatzko}, \binits{S.}},
\bauthor{\bsnm{Kleyman}, \binits{G.}},
\bauthor{\bsnm{Wallaschek}, \binits{J.}}:
\batitle{Continuation methods for lab experiments of nonlinear vibrations}.
\bjtitle{GAMM Mitteilungen}
\bvolume{46}(\bissue{2}),
\bfpage{1}--\blpage{13}
(\byear{2023})
\doiurl{10.1002/gamm.202300009}
\end{barticle}
\endbibitem

\bibitem[\protect\citeauthoryear{Hayashi et~al.}{2024}]{Hayashi2024}
\begin{barticle}
\bauthor{\bsnm{Hayashi}, \binits{S.}},
\bauthor{\bsnm{Gutschmidt}, \binits{S.}},
\bauthor{\bsnm{Murray}, \binits{R.}},
\bauthor{\bsnm{Krauskopf}, \binits{B.}}:
\batitle{Experimental bifurcation analysis of a clamped beam with designed
  mechanical nonlinearity}.
\bjtitle{Nonlinear Dynamics}
(\byear{2024})
\doiurl{10.1007/s11071-024-09873-5}
\end{barticle}
\endbibitem

\bibitem[\protect\citeauthoryear{Bureau et~al.}{2013}]{Bureau2013}
\begin{barticle}
\bauthor{\bsnm{Bureau}, \binits{E.}},
\bauthor{\bsnm{Schilder}, \binits{F.}},
\bauthor{\bsnm{{Ferreira Santos}}, \binits{I.}},
\bauthor{\bsnm{{Juel Thomsen}}, \binits{J.}},
\bauthor{\bsnm{Starke}, \binits{J.}}:
\batitle{Experimental bifurcation analysis of an impact oscillator—tuning a
  non-invasive control scheme}.
\bjtitle{Journal of Sound and Vibration}
\bvolume{332}(\bissue{22}),
\bfpage{5883}--\blpage{5897}
(\byear{2013})
\doiurl{10.1016/j.jsv.2013.05.033}
\end{barticle}
\endbibitem

\bibitem[\protect\citeauthoryear{Bureau}{2014}]{Bureau2014Thesis}
\begin{botherref}
\oauthor{\bsnm{Bureau}, \binits{E.}}:
Experimental bifurcation analysis using control-based continuation.
PhD thesis,
Technical University of Denmark
(2014).
\url{https://findit.dtu.dk/en/catalog/53be7cc5f7d9a21b490001de}
\end{botherref}
\endbibitem

\bibitem[\protect\citeauthoryear{Siettos et~al.}{2004}]{SIETTOS2004}
\begin{barticle}
\bauthor{\bsnm{Siettos}, \binits{C.I.}},
\bauthor{\bsnm{Kevrekidis}, \binits{I.G.}},
\bauthor{\bsnm{Maroudas}, \binits{D.}}:
\batitle{Coarse bifurcation diagrams via microscopic simulators: A
  state-feedback control-based approach}.
\bjtitle{International Journal of Bifurcation and Chaos}
\bvolume{14}(\bissue{01}),
\bfpage{207}--\blpage{220}
(\byear{2004})
\doiurl{10.1142/S0218127404009193}
\end{barticle}
\endbibitem

\bibitem[\protect\citeauthoryear{Barton et~al.}{2012}]{Barton2012}
\begin{barticle}
\bauthor{\bsnm{Barton}, \binits{D.A.W.}},
\bauthor{\bsnm{Mann}, \binits{B.P.}},
\bauthor{\bsnm{Burrow}, \binits{S.G.}}:
\batitle{Control-based continuation for investigating nonlinear experiments}.
\bjtitle{Journal of Vibration and Control}
\bvolume{18}(\bissue{4}),
\bfpage{509}--\blpage{520}
(\byear{2012})
\doiurl{10.1177/1077546310384004}
\end{barticle}
\endbibitem

\bibitem[\protect\citeauthoryear{Sieber et~al.}{2011}]{Sieber2011}
\begin{barticle}
\bauthor{\bsnm{Sieber}, \binits{J.}},
\bauthor{\bsnm{Krauskopf}, \binits{B.}},
\bauthor{\bsnm{Wagg}, \binits{D.}},
\bauthor{\bsnm{Neild}, \binits{S.}},
\bauthor{\bsnm{Gonzalez-Buelga}, \binits{A.}}:
\batitle{Control-based continuation of unstable periodic orbits}.
\bjtitle{Journal of Computational and Nonlinear Dynamics}
\bvolume{6}(\bissue{1}),
\bfpage{1}--\blpage{9}
(\byear{2011})
\doiurl{10.1115/1.4002101}
\end{barticle}
\endbibitem

\bibitem[\protect\citeauthoryear{Barton and Burrow}{2011}]{Barton2011}
\begin{barticle}
\bauthor{\bsnm{Barton}, \binits{D.A.W.}},
\bauthor{\bsnm{Burrow}, \binits{S.G.}}:
\batitle{Numerical continuation in a physical experiment: Investigation of a
  nonlinear energy harvester}.
\bjtitle{Journal of Computational and Nonlinear Dynamics}
\bvolume{6}(\bissue{1}),
\bfpage{1}--\blpage{6}
(\byear{2011})
\doiurl{10.1115/1.4002380}
\end{barticle}
\endbibitem

\bibitem[\protect\citeauthoryear{Kleyman et~al.}{2024}]{Kleyman2024}
\begin{bchapter}
\bauthor{\bsnm{Kleyman}, \binits{G.}},
\bauthor{\bsnm{Jahn}, \binits{M.}},
\bauthor{\bsnm{Tatzko}, \binits{S.}},
\bauthor{\bsnm{Scheidt}, \binits{L.P.-v.}}:
\bctitle{A combined numerical-experimental approach for the damping evaluation
  of non-linear dissipative vibration systems}.
In: \bbtitle{Lecture Notes in Applied and Computational Mechanics}
vol. \bseriesno{102},
pp. \bfpage{285}--\blpage{303}.
\bpublisher{Springer},
\blocation{Cham}
(\byear{2024}).
\doiurl{10.1007/978-3-031-36143-2_15} .
\burl{https://link.springer.com/10.1007/978-3-031-36143-2{\_}15}
\end{bchapter}
\endbibitem

\bibitem[\protect\citeauthoryear{van Iderstein and
  Wiebe}{2019}]{VanIderstein2019}
\begin{barticle}
\bauthor{\bsnm{Iderstein}, \binits{T.}},
\bauthor{\bsnm{Wiebe}, \binits{R.}}:
\batitle{Experimental path following of unstable static equilibria for
  snap-through buckling}.
\bjtitle{Conference Proceedings of the Society for Experimental Mechanics
  Series}
\bvolume{1},
\bfpage{17}--\blpage{22}
(\byear{2019})
\doiurl{10.1007/978-3-319-74280-9_2}
\end{barticle}
\endbibitem

\bibitem[\protect\citeauthoryear{Shen et~al.}{2021}]{Shen2021}
\begin{barticle}
\bauthor{\bsnm{Shen}, \binits{J.}},
\bauthor{\bsnm{Groh}, \binits{R.M.J.}},
\bauthor{\bsnm{Schenk}, \binits{M.}},
\bauthor{\bsnm{Pirrera}, \binits{A.}}:
\batitle{Experimental path-following of equilibria using newton's method. part
  i: Theory, modelling, experiments}.
\bjtitle{International Journal of Solids and Structures}
\bvolume{210-211},
\bfpage{203}--\blpage{223}
(\byear{2021})
\doiurl{10.1016/j.ijsolstr.2020.11.037}
\end{barticle}
\endbibitem

\bibitem[\protect\citeauthoryear{Beregi et~al.}{2023}]{Beregi2023p}
\begin{botherref}
\oauthor{\bsnm{Beregi}, \binits{S.}},
\oauthor{\bsnm{Barton}, \binits{D.A.W.}},
\oauthor{\bsnm{Rezgui}, \binits{D.}},
\oauthor{\bsnm{Neild}, \binits{S.A.}}:
Real-time hybrid testing using iterative control for periodic oscillations
(2023)
{\href{https://arxiv.org/abs/2312.06362}{{arXiv:2312.06362}}}
\end{botherref}
\endbibitem

\bibitem[\protect\citeauthoryear{Melville and Suarez}{2017}]{Melville2017}
\begin{barticle}
\bauthor{\bsnm{Melville}, \binits{R.C.}},
\bauthor{\bsnm{Suarez}, \binits{A.}}:
\batitle{Experimental investigation of bifurcation behavior in nonlinear
  microwave circuits}.
\bjtitle{IEEE Transactions on Microwave Theory and Techniques}
\bvolume{65}(\bissue{5}),
\bfpage{1545}--\blpage{1559}
(\byear{2017})
\doiurl{10.1109/TMTT.2016.2640955}
\end{barticle}
\endbibitem

\bibitem[\protect\citeauthoryear{Bureau et~al.}{2014}]{Bureau2014}
\begin{barticle}
\bauthor{\bsnm{Bureau}, \binits{E.}},
\bauthor{\bsnm{Schilder}, \binits{F.}},
\bauthor{\bsnm{Elmeg{\aa}rd}, \binits{M.}},
\bauthor{\bsnm{Santos}, \binits{I.F.}},
\bauthor{\bsnm{Thomsen}, \binits{J.J.}},
\bauthor{\bsnm{Starke}, \binits{J.}}:
\batitle{Experimental bifurcation analysis of an impact
  oscillator—determining stability}.
\bjtitle{Journal of Sound and Vibration}
\bvolume{333}(\bissue{21}),
\bfpage{5464}--\blpage{5474}
(\byear{2014})
\doiurl{10.1016/j.jsv.2014.05.032}
\end{barticle}
\endbibitem

\bibitem[\protect\citeauthoryear{Schilder et~al.}{2015}]{Schilder2015}
\begin{barticle}
\bauthor{\bsnm{Schilder}, \binits{F.}},
\bauthor{\bsnm{Bureau}, \binits{E.}},
\bauthor{\bsnm{Santos}, \binits{I.F.}},
\bauthor{\bsnm{Thomsen}, \binits{J.J.}},
\bauthor{\bsnm{Starke}, \binits{J.}}:
\batitle{Experimental bifurcation analysis—continuation for
  noise-contaminated zero problems}.
\bjtitle{Journal of Sound and Vibration}
\bvolume{358},
\bfpage{251}--\blpage{266}
(\byear{2015})
\doiurl{10.1016/j.jsv.2015.08.008}
\end{barticle}
\endbibitem

\bibitem[\protect\citeauthoryear{Blyth et~al.}{2023}]{Blyth2023}
\begin{barticle}
\bauthor{\bsnm{Blyth}, \binits{M.}},
\bauthor{\bsnm{Tsaneva-Atanasova}, \binits{K.}},
\bauthor{\bsnm{Marucci}, \binits{L.}},
\bauthor{\bsnm{Renson}, \binits{L.}}:
\batitle{Numerical methods for control-based continuation of relaxation
  oscillations}.
\bjtitle{Nonlinear Dynamics}
\bvolume{111}(\bissue{9}),
\bfpage{7975}--\blpage{7992}
(\byear{2023})
\doiurl{10.1007/s11071-023-08288-y}
\end{barticle}
\endbibitem

\bibitem[\protect\citeauthoryear{Schilder and Dankowicz}{}]{COCO}
\begin{botherref}
\oauthor{\bsnm{Schilder}, \binits{F.}},
\oauthor{\bsnm{Dankowicz}, \binits{H.}}:
Continuation Core and Toolboxes (COCO).
\url{https://sourceforge.net/projects/cocotools/}
\end{botherref}
\endbibitem

\bibitem[\protect\citeauthoryear{Cenedese and Haller}{2020}]{Cenedese2020}
\begin{barticle}
\bauthor{\bsnm{Cenedese}, \binits{M.}},
\bauthor{\bsnm{Haller}, \binits{G.}}:
\batitle{How do conservative backbone curves perturb into forced responses? a
  {M}elnikov function analysis}.
\bjtitle{Proceedings of the Royal Society A: Mathematical, Physical and
  Engineering Sciences}
\bvolume{476}(\bissue{2234}),
\bfpage{20190494}
(\byear{2020})
\doiurl{10.1098/rspa.2019.0494}
{\href{https://arxiv.org/abs/1908.00721}{{arXiv:1908.00721}}}
\end{barticle}
\endbibitem

\bibitem[\protect\citeauthoryear{Renson et~al.}{2017}]{Renson2017}
\begin{barticle}
\bauthor{\bsnm{Renson}, \binits{L.}},
\bauthor{\bsnm{Barton}, \binits{D.A.W.}},
\bauthor{\bsnm{Neild}, \binits{S.A.}}:
\batitle{Experimental tracking of limit-point bifurcations and backbone curves
  using control-based continuation}.
\bjtitle{International Journal of Bifurcation and Chaos}
\bvolume{27}(\bissue{1}),
\bfpage{1}--\blpage{19}
(\byear{2017})
\doiurl{10.1142/S0218127417300026}
\end{barticle}
\endbibitem

\bibitem[\protect\citeauthoryear{Woiwode and Krack}{2024}]{Woiwode2024}
\begin{botherref}
\oauthor{\bsnm{Woiwode}, \binits{L.}},
\oauthor{\bsnm{Krack}, \binits{M.}}:
Experimentally uncovering isolas via backbone tracking.
Journal of Structural Dynamics,
122--143
(2024)
\doiurl{10.25518/2684-6500.180}
\end{botherref}
\endbibitem

\bibitem[\protect\citeauthoryear{Kleyman et~al.}{2020}]{Kleyman2020}
\begin{barticle}
\bauthor{\bsnm{Kleyman}, \binits{G.}},
\bauthor{\bsnm{Paehr}, \binits{M.}},
\bauthor{\bsnm{Tatzko}, \binits{S.}}:
\batitle{Application of control-based-continuation for characterization of
  dynamic systems with stiffness and friction nonlinearities}.
\bjtitle{Mechanics Research Communications}
\bvolume{106},
\bfpage{103520}
(\byear{2020})
\doiurl{10.1016/j.mechrescom.2020.103520}
\end{barticle}
\endbibitem

\bibitem[\protect\citeauthoryear{Renson et~al.}{2016}]{Renson2016}
\begin{barticle}
\bauthor{\bsnm{Renson}, \binits{L.}},
\bauthor{\bsnm{Gonzalez-Buelga}, \binits{A.}},
\bauthor{\bsnm{Barton}, \binits{D.A.W.}},
\bauthor{\bsnm{Neild}, \binits{S.A.}}:
\batitle{Robust identification of backbone curves using control-based
  continuation}.
\bjtitle{Journal of Sound and Vibration}
\bvolume{367},
\bfpage{145}--\blpage{158}
(\byear{2016})
\doiurl{10.1016/j.jsv.2015.12.035}
\end{barticle}
\endbibitem

\bibitem[\protect\citeauthoryear{Renson et~al.}{2019}]{Renson2019}
\begin{barticle}
\bauthor{\bsnm{Renson}, \binits{L.}},
\bauthor{\bsnm{Shaw}, \binits{A.D.}},
\bauthor{\bsnm{Barton}, \binits{D.A.W.}},
\bauthor{\bsnm{Neild}, \binits{S.A.}}:
\batitle{Application of control-based continuation to a nonlinear structure
  with harmonically coupled modes}.
\bjtitle{Mechanical Systems and Signal Processing}
\bvolume{120},
\bfpage{449}--\blpage{464}
(\byear{2019})
\doiurl{10.1016/j.ymssp.2018.10.008}
{\href{https://arxiv.org/abs/1808.01865}{{arXiv:1808.01865}}}
\end{barticle}
\endbibitem

\bibitem[\protect\citeauthoryear{Abeloos et~al.}{2021}]{Abeloos2021}
\begin{barticle}
\bauthor{\bsnm{Abeloos}, \binits{G.}},
\bauthor{\bsnm{Renson}, \binits{L.}},
\bauthor{\bsnm{Collette}, \binits{C.}},
\bauthor{\bsnm{Kerschen}, \binits{G.}}:
\batitle{Stepped and swept control-based continuation using adaptive
  filtering}.
\bjtitle{Nonlinear Dynamics}
\bvolume{104}(\bissue{4}),
\bfpage{3793}--\blpage{3808}
(\byear{2021})
\doiurl{10.1007/s11071-021-06506-z}
\end{barticle}
\endbibitem

\bibitem[\protect\citeauthoryear{Beregi}{2022}]{Beregi2022}
\begin{barticle}
\bauthor{\bsnm{Beregi}, \binits{S.}}:
\batitle{Nonlinear analysis of the delayed tyre model with control-based
  continuation}.
\bjtitle{Nonlinear Dynamics}
\bvolume{110}(\bissue{4}),
\bfpage{3151}--\blpage{3165}
(\byear{2022})
\doiurl{10.1007/s11071-022-07796-7}
\end{barticle}
\endbibitem

\bibitem[\protect\citeauthoryear{Renson et~al.}{2019}]{Renson2019b}
\begin{barticle}
\bauthor{\bsnm{Renson}, \binits{L.}},
\bauthor{\bsnm{Sieber}, \binits{J.}},
\bauthor{\bsnm{Barton}, \binits{D.A.W.}},
\bauthor{\bsnm{Shaw}, \binits{A.D.}},
\bauthor{\bsnm{Neild}, \binits{S.A.}}:
\batitle{Numerical continuation in nonlinear experiments using local gaussian
  process regression}.
\bjtitle{Nonlinear Dynamics}
\bvolume{98}(\bissue{4}),
\bfpage{2811}--\blpage{2826}
(\byear{2019})
\doiurl{10.1007/s11071-019-05118-y}
\end{barticle}
\endbibitem

\bibitem[\protect\citeauthoryear{Beregi et~al.}{2021}]{Beregi2021}
\begin{barticle}
\bauthor{\bsnm{Beregi}, \binits{S.}},
\bauthor{\bsnm{Barton}, \binits{D.A.W.}},
\bauthor{\bsnm{Rezgui}, \binits{D.}},
\bauthor{\bsnm{Neild}, \binits{S.A.}}:
\batitle{Robustness of nonlinear parameter identification in the presence of
  process noise using control-based continuation}.
\bjtitle{Nonlinear Dynamics}
\bvolume{104}(\bissue{2}),
\bfpage{885}--\blpage{900}
(\byear{2021})
\doiurl{10.1007/s11071-021-06347-w}
{\href{https://arxiv.org/abs/2001.11008}{{arXiv:2001.11008}}}
\end{barticle}
\endbibitem

\bibitem[\protect\citeauthoryear{Beregi et~al.}{2023}]{Beregi2023}
\begin{barticle}
\bauthor{\bsnm{Beregi}, \binits{S.}},
\bauthor{\bsnm{Barton}, \binits{D.A.W.}},
\bauthor{\bsnm{Rezgui}, \binits{D.}},
\bauthor{\bsnm{Neild}, \binits{S.}}:
\batitle{Using scientific machine learning for experimental bifurcation
  analysis of dynamic systems}.
\bjtitle{Mechanical Systems and Signal Processing}
\bvolume{184}(\bissue{June 2022}),
\bfpage{109649}
(\byear{2023})
\doiurl{10.1016/j.ymssp.2022.109649}
{\href{https://arxiv.org/abs/2110.11854}{{arXiv:2110.11854}}}
\end{barticle}
\endbibitem

\bibitem[\protect\citeauthoryear{Lee et~al.}{2023}]{Lee2023}
\begin{barticle}
\bauthor{\bsnm{Lee}, \binits{K.H.}},
\bauthor{\bsnm{Barton}, \binits{D.A.W.}},
\bauthor{\bsnm{Renson}, \binits{L.}}:
\batitle{Modelling of physical systems with a hopf bifurcation using
  mechanistic models and machine learning}.
\bjtitle{Mechanical Systems and Signal Processing}
\bvolume{191}(\bissue{October 2022}),
\bfpage{110173}
(\byear{2023})
\doiurl{10.1016/j.ymssp.2023.110173}
\end{barticle}
\endbibitem

\bibitem[\protect\citeauthoryear{Barton}{2017}]{Barton2017}
\begin{barticle}
\bauthor{\bsnm{Barton}, \binits{D.A.W.}}:
\batitle{Control-based continuation: Bifurcation and stability analysis for
  physical experiments}.
\bjtitle{Mechanical Systems and Signal Processing}
\bvolume{84},
\bfpage{54}--\blpage{64}
(\byear{2017})
\doiurl{10.1016/j.ymssp.2015.12.039}
{\href{https://arxiv.org/abs/1506.04052}{{arXiv:1506.04052}}}
\end{barticle}
\endbibitem

\bibitem[\protect\citeauthoryear{Dittus et~al.}{2023}]{Dittus2023}
\begin{barticle}
\bauthor{\bsnm{Dittus}, \binits{A.}},
\bauthor{\bsnm{Kruse}, \binits{N.}},
\bauthor{\bsnm{Barke}, \binits{I.}},
\bauthor{\bsnm{Speller}, \binits{S.}},
\bauthor{\bsnm{Starke}, \binits{J.}}:
\batitle{Detecting stability and bifurcation points in control-based
  continuation for a physical experiment of the zeeman catastrophe machine}.
\bjtitle{SIAM Journal on Applied Dynamical Systems}
\bvolume{22}(\bissue{2}),
\bfpage{1275}--\blpage{1299}
(\byear{2023})
\doiurl{10.1137/22M1503245}
\end{barticle}
\endbibitem

\bibitem[\protect\citeauthoryear{Tartaruga et~al.}{2019}]{Tartaruga2019}
\begin{bchapter}
\bauthor{\bsnm{Tartaruga}, \binits{I.}},
\bauthor{\bsnm{Barton}, \binits{D.A.W.}},
\bauthor{\bsnm{Rezgui}, \binits{D.}},
\bauthor{\bsnm{Neild}, \binits{S.A.}}:
\bctitle{Experimental bifurcation analysis of a wing profile}.
In: \bbtitle{International Forum on Aeroelasticity and Structural Dynamics :
  IFASD 2019},
\bconflocation{Savannah}
(\byear{2019}).
\burl{https://hdl.handle.net/1983/6774dc99-50b6-4e51-ba30-f06207581973}
\end{bchapter}
\endbibitem

\bibitem[\protect\citeauthoryear{{De Cesare} et~al.}{2022}]{DeCesare2022}
\begin{barticle}
\bauthor{\bsnm{{De Cesare}}, \binits{I.}},
\bauthor{\bsnm{Salzano}, \binits{D.}},
\bauthor{\bsnm{{Di Bernardo}}, \binits{M.}},
\bauthor{\bsnm{Renson}, \binits{L.}},
\bauthor{\bsnm{Marucci}, \binits{L.}}:
\batitle{Control-based continuation: A new approach to prototype synthetic gene
  networks}.
\bjtitle{ACS Synthetic Biology}
\bvolume{11}(\bissue{7}),
\bfpage{2300}--\blpage{2313}
(\byear{2022})
\doiurl{10.1021/acssynbio.1c00632}
\end{barticle}
\endbibitem

\bibitem[\protect\citeauthoryear{Hayashi et~al.}{2024}]{Hayashi2024c}
\begin{bchapter}
\bauthor{\bsnm{Hayashi}, \binits{S.}},
\bauthor{\bsnm{Gutschmidt}, \binits{S.}},
\bauthor{\bsnm{Murray}, \binits{R.}},
\bauthor{\bsnm{Krauskopf}, \binits{B.}}:
\bctitle{Control-based continuation of an externally excited mems
  self-oscillator}.
In: \bbtitle{ENOC 2024},
\bconflocation{Delft}
(\byear{2024})
\end{bchapter}
\endbibitem

\bibitem[\protect\citeauthoryear{Panagiotopoulos
  et~al.}{2022}]{Panagiotopoulos2022}
\begin{barticle}
\bauthor{\bsnm{Panagiotopoulos}, \binits{I.}},
\bauthor{\bsnm{Starke}, \binits{J.}},
\bauthor{\bsnm{Just}, \binits{W.}}:
\batitle{Control of collective human behavior: Social dynamics beyond
  modeling}.
\bjtitle{Physical Review Research}
\bvolume{4}(\bissue{4}),
\bfpage{043190}
(\byear{2022})
\doiurl{10.1103/PhysRevResearch.4.043190}
\end{barticle}
\endbibitem

\bibitem[\protect\citeauthoryear{Panagiotopoulos
  et~al.}{2023}]{Panagiotopoulos2023}
\begin{barticle}
\bauthor{\bsnm{Panagiotopoulos}, \binits{I.}},
\bauthor{\bsnm{Starke}, \binits{J.}},
\bauthor{\bsnm{Sieber}, \binits{J.}},
\bauthor{\bsnm{Just}, \binits{W.}}:
\batitle{Continuation with noninvasive control schemes: Revealing unstable
  states in a pedestrian evacuation scenario}.
\bjtitle{SIAM Journal on Applied Dynamical Systems}
\bvolume{22}(\bissue{1}),
\bfpage{1}--\blpage{36}
(\byear{2023})
\doiurl{10.1137/22M1482032}
{\href{https://arxiv.org/abs/2203.02484}{{arXiv:2203.02484}}}
\end{barticle}
\endbibitem

\bibitem[\protect\citeauthoryear{Babitsky}{1995}]{Babitsky1995}
\begin{barticle}
\bauthor{\bsnm{Babitsky}, \binits{V.I.}}:
\batitle{Autoresonant mechatronic systems}.
\bjtitle{Mechatronics}
\bvolume{5}(\bissue{5}),
\bfpage{483}--\blpage{495}
(\byear{1995})
\doiurl{10.1016/0957-4158(95)00026-2}
\end{barticle}
\endbibitem

\bibitem[\protect\citeauthoryear{Sokolov and Babitsky}{2001}]{SOKOLOV2001}
\begin{barticle}
\bauthor{\bsnm{Sokolov}, \binits{I.J.}},
\bauthor{\bsnm{Babitsky}, \binits{V.I.}}:
\batitle{Phase control of self-sustained vibration}.
\bjtitle{Journal of Sound and Vibration}
\bvolume{248}(\bissue{4}),
\bfpage{725}--\blpage{744}
(\byear{2001})
\doiurl{10.1006/jsvi.2001.3810}
\end{barticle}
\endbibitem

\bibitem[\protect\citeauthoryear{Mojrzisch et~al.}{2012}]{Mojrzisch2012}
\begin{barticle}
\bauthor{\bsnm{Mojrzisch}, \binits{S.}},
\bauthor{\bsnm{Wallaschek}, \binits{J.}},
\bauthor{\bsnm{Bremer}, \binits{J.}}:
\batitle{An experimental method for the phase controlled frequency response
  measurement of nonlinear vibration systems}.
\bjtitle{PAMM}
\bvolume{12}(\bissue{1}),
\bfpage{253}--\blpage{254}
(\byear{2012})
\doiurl{10.1002/pamm.201210117}
\end{barticle}
\endbibitem

\bibitem[\protect\citeauthoryear{Davis and Bucher}{2018}]{Davis2018}
\begin{barticle}
\bauthor{\bsnm{Davis}, \binits{S.}},
\bauthor{\bsnm{Bucher}, \binits{I.}}:
\batitle{Automatic vibration mode selection and excitation; combining modal
  filtering with autoresonance}.
\bjtitle{Mechanical Systems and Signal Processing}
\bvolume{101},
\bfpage{140}--\blpage{155}
(\byear{2018})
\doiurl{10.1016/j.ymssp.2017.08.009}
\end{barticle}
\endbibitem

\bibitem[\protect\citeauthoryear{Scheel}{2022}]{Scheel2022}
\begin{barticle}
\bauthor{\bsnm{Scheel}, \binits{M.}}:
\batitle{Nonlinear modal testing of damped structures: Velocity feedback vs.
  phase resonance}.
\bjtitle{Mechanical Systems and Signal Processing}
\bvolume{165}(\bissue{June 2021}),
\bfpage{108305}
(\byear{2022})
\doiurl{10.1016/j.ymssp.2021.108305}
{\href{https://arxiv.org/abs/2108.06189}{{arXiv:2108.06189}}}
\end{barticle}
\endbibitem

\bibitem[\protect\citeauthoryear{Abramovitch}{2002}]{Abramovitch2002}
\begin{bchapter}
\bauthor{\bsnm{Abramovitch}, \binits{D.}}:
\bctitle{Phase-locked loops: a control centric tutorial}.
In: \bbtitle{Proceedings of the 2002 American Control Conference (IEEE Cat.
  No.CH37301)},
vol. \bseriesno{1},
pp. \bfpage{1}--\blpage{15}.
\bpublisher{American Automatic Control Council},
\blocation{Anchorage}
(\byear{2002}).
\doiurl{10.1109/ACC.2002.1024769} .
\burl{http://ieeexplore.ieee.org/document/1024769/}
\end{bchapter}
\endbibitem

\bibitem[\protect\citeauthoryear{Balas}{1978}]{Balas1978}
\begin{barticle}
\bauthor{\bsnm{Balas}, \binits{M.J.}}:
\batitle{Active control of flexible systems}.
\bjtitle{Journal of Optimization Theory and Applications}
\bvolume{25}(\bissue{3}),
\bfpage{415}--\blpage{436}
(\byear{1978})
\doiurl{10.1007/BF00932903}
\end{barticle}
\endbibitem

\bibitem[\protect\citeauthoryear{Niezrecki and Cudney}{1997}]{Niezrecki1997}
\begin{barticle}
\bauthor{\bsnm{Niezrecki}, \binits{C.}},
\bauthor{\bsnm{Cudney}, \binits{H.H.}}:
\batitle{Structural control using analog phase-locked loops}.
\bjtitle{Journal of Vibration and Acoustics}
\bvolume{119}(\bissue{1}),
\bfpage{104}--\blpage{109}
(\byear{1997})
\doiurl{10.1115/1.2889677}
\end{barticle}
\endbibitem

\bibitem[\protect\citeauthoryear{Connally and Brown}{1993}]{Connally1993}
\begin{barticle}
\bauthor{\bsnm{Connally}, \binits{J.A.}},
\bauthor{\bsnm{Brown}, \binits{S.B.}}:
\batitle{Micromechanical fatigue testing}.
\bjtitle{Experimental Mechanics}
\bvolume{33}(\bissue{2}),
\bfpage{81}--\blpage{90}
(\byear{1993})
\doiurl{10.1007/BF02322482}
\end{barticle}
\endbibitem

\bibitem[\protect\citeauthoryear{Kern and Seemann}{2010}]{Kern2010}
\begin{barticle}
\bauthor{\bsnm{Kern}, \binits{D.}},
\bauthor{\bsnm{Seemann}, \binits{W.}}:
\batitle{Tracking of mechanical system parameters by phase‐locked loops}.
\bjtitle{PAMM}
\bvolume{10}(\bissue{1}),
\bfpage{613}--\blpage{614}
(\byear{2010})
\doiurl{10.1002/pamm.201010299}
\end{barticle}
\endbibitem

\bibitem[\protect\citeauthoryear{Kern et~al.}{2012}]{Kern2012}
\begin{barticle}
\bauthor{\bsnm{Kern}, \binits{D.}},
\bauthor{\bsnm{Brack}, \binits{T.}},
\bauthor{\bsnm{Seemann}, \binits{W.}}:
\batitle{Resonance tracking of continua using self-sensing actuators}.
\bjtitle{Journal of Dynamic Systems, Measurement, and Control}
\bvolume{134}(\bissue{5}),
\bfpage{1}--\blpage{9}
(\byear{2012})
\doiurl{10.1115/1.4006224}
\end{barticle}
\endbibitem

\bibitem[\protect\citeauthoryear{Denis et~al.}{2018}]{Denis2018}
\begin{barticle}
\bauthor{\bsnm{Denis}, \binits{V.}},
\bauthor{\bsnm{Jossic}, \binits{M.}},
\bauthor{\bsnm{Giraud-Audine}, \binits{C.}},
\bauthor{\bsnm{Chomette}, \binits{B.}},
\bauthor{\bsnm{Renault}, \binits{A.}},
\bauthor{\bsnm{Thomas}, \binits{O.}}:
\batitle{Identification of nonlinear modes using phase-locked-loop experimental
  continuation and normal form}.
\bjtitle{Mechanical Systems and Signal Processing}
\bvolume{106},
\bfpage{430}--\blpage{452}
(\byear{2018})
\doiurl{10.1016/j.ymssp.2018.01.014}
\end{barticle}
\endbibitem

\bibitem[\protect\citeauthoryear{Peter et~al.}{2018}]{Peter2018}
\begin{barticle}
\bauthor{\bsnm{Peter}, \binits{S.}},
\bauthor{\bsnm{Scheel}, \binits{M.}},
\bauthor{\bsnm{Krack}, \binits{M.}},
\bauthor{\bsnm{Leine}, \binits{R.I.}}:
\batitle{Synthesis of nonlinear frequency responses with experimentally
  extracted nonlinear modes}.
\bjtitle{Mechanical Systems and Signal Processing}
\bvolume{101},
\bfpage{498}--\blpage{515}
(\byear{2018})
\doiurl{10.1016/j.ymssp.2017.09.014}
\end{barticle}
\endbibitem

\bibitem[\protect\citeauthoryear{Scheel et~al.}{2018}]{Scheel2018}
\begin{barticle}
\bauthor{\bsnm{Scheel}, \binits{M.}},
\bauthor{\bsnm{Peter}, \binits{S.}},
\bauthor{\bsnm{Leine}, \binits{R.I.}},
\bauthor{\bsnm{Krack}, \binits{M.}}:
\batitle{A phase resonance approach for modal testing of structures with
  nonlinear dissipation}.
\bjtitle{Journal of Sound and Vibration}
\bvolume{435},
\bfpage{56}--\blpage{73}
(\byear{2018})
\doiurl{10.1016/j.jsv.2018.07.010}
{\href{https://arxiv.org/abs/2011.08500}{{arXiv:2011.08500}}}
\end{barticle}
\endbibitem

\bibitem[\protect\citeauthoryear{Scheel et~al.}{2020}]{Scheel2020b}
\begin{barticle}
\bauthor{\bsnm{Scheel}, \binits{M.}},
\bauthor{\bsnm{Weigele}, \binits{T.}},
\bauthor{\bsnm{Krack}, \binits{M.}}:
\batitle{Challenging an experimental nonlinear modal analysis method with a new
  strongly friction-damped structure}.
\bjtitle{Journal of Sound and Vibration}
\bvolume{485},
\bfpage{115580}
(\byear{2020})
\doiurl{10.1016/j.jsv.2020.115580}
{\href{https://arxiv.org/abs/2011.08527}{{arXiv:2011.08527}}}
\end{barticle}
\endbibitem

\bibitem[\protect\citeauthoryear{Peter et~al.}{2019}]{Peter2019}
\begin{barticle}
\bauthor{\bsnm{Peter}, \binits{S.}},
\bauthor{\bsnm{Schreyer}, \binits{F.}},
\bauthor{\bsnm{Leine}, \binits{R.I.}}:
\batitle{A method for numerical and experimental nonlinear modal analysis of
  nonsmooth systems}.
\bjtitle{Mechanical Systems and Signal Processing}
\bvolume{120},
\bfpage{793}--\blpage{807}
(\byear{2019})
\doiurl{10.1016/j.ymssp.2018.11.009}
\end{barticle}
\endbibitem

\bibitem[\protect\citeauthoryear{Scheel et~al.}{2020}]{Scheel2020}
\begin{barticle}
\bauthor{\bsnm{Scheel}, \binits{M.}},
\bauthor{\bsnm{Kleyman}, \binits{G.}},
\bauthor{\bsnm{Tatar}, \binits{A.}},
\bauthor{\bsnm{Brake}, \binits{M.R.W.}},
\bauthor{\bsnm{Peter}, \binits{S.}},
\bauthor{\bsnm{No{\"{e}}l}, \binits{J.-P.}},
\bauthor{\bsnm{Allen}, \binits{M.S.}},
\bauthor{\bsnm{Krack}, \binits{M.}}:
\batitle{Experimental assessment of polynomial nonlinear state-space and
  nonlinear-mode models for near-resonant vibrations}.
\bjtitle{Mechanical Systems and Signal Processing}
\bvolume{143},
\bfpage{106796}
(\byear{2020})
\doiurl{10.1016/j.ymssp.2020.106796}
\end{barticle}
\endbibitem

\bibitem[\protect\citeauthoryear{Debeurre et~al.}{2024}]{Debeurre2024}
\begin{barticle}
\bauthor{\bsnm{Debeurre}, \binits{M.}},
\bauthor{\bsnm{Benacchio}, \binits{S.}},
\bauthor{\bsnm{Grolet}, \binits{A.}},
\bauthor{\bsnm{Grenat}, \binits{C.}},
\bauthor{\bsnm{Giraud-Audine}, \binits{C.}},
\bauthor{\bsnm{Thomas}, \binits{O.}}:
\batitle{Phase resonance testing of highly flexible structures: Measurement of
  conservative nonlinear modes and nonlinear damping identification}.
\bjtitle{Mechanical Systems and Signal Processing}
\bvolume{215}(\bissue{February}),
\bfpage{111423}
(\byear{2024})
\doiurl{10.1016/j.ymssp.2024.111423}
\end{barticle}
\endbibitem

\bibitem[\protect\citeauthoryear{Nagesh et~al.}{2022}]{Nagesh2022}
\begin{bchapter}
\bauthor{\bsnm{Nagesh}, \binits{M.}},
\bauthor{\bsnm{Allemang}, \binits{R.J.}},
\bauthor{\bsnm{Phillips}, \binits{A.W.}}:
\bctitle{Characterization of nonlinearities in a structure using nonlinear
  modal testing methods}.
In: \bbtitle{Conference Proceedings of the Society for Experimental Mechanics
  Series}
vol. \bseriesno{1},
pp. \bfpage{167}--\blpage{178}
(\byear{2022}).
\doiurl{10.1007/978-3-030-77135-5_19} .
\burl{https://link.springer.com/10.1007/978-3-030-77135-5{\_}19}
\end{bchapter}
\endbibitem

\bibitem[\protect\citeauthoryear{Givois et~al.}{2020a}]{Givois2020}
\begin{barticle}
\bauthor{\bsnm{Givois}, \binits{A.}},
\bauthor{\bsnm{Giraud-Audine}, \binits{C.}},
\bauthor{\bsnm{De{\"{u}}}, \binits{J.-F.}},
\bauthor{\bsnm{Thomas}, \binits{O.}}:
\batitle{Experimental analysis of nonlinear resonances in piezoelectric plates
  with geometric nonlinearities}.
\bjtitle{Nonlinear Dynamics}
\bvolume{102}(\bissue{3}),
\bfpage{1451}--\blpage{1462}
(\byear{2020})
\doiurl{10.1007/s11071-020-05997-6}
\end{barticle}
\endbibitem

\bibitem[\protect\citeauthoryear{Givois et~al.}{2020b}]{Givois2020b}
\begin{barticle}
\bauthor{\bsnm{Givois}, \binits{A.}},
\bauthor{\bsnm{Tan}, \binits{J.-J.}},
\bauthor{\bsnm{Touz{\'{e}}}, \binits{C.}},
\bauthor{\bsnm{Thomas}, \binits{O.}}:
\batitle{Backbone curves of coupled cubic oscillators in one-to-one internal
  resonance: bifurcation scenario, measurements and parameter identification}.
\bjtitle{Meccanica}
\bvolume{55}(\bissue{3}),
\bfpage{481}--\blpage{503}
(\byear{2020})
\doiurl{10.1007/s11012-020-01132-2}
\end{barticle}
\endbibitem

\bibitem[\protect\citeauthoryear{Schwarz et~al.}{2020}]{Schwarz2020}
\begin{barticle}
\bauthor{\bsnm{Schwarz}, \binits{S.}},
\bauthor{\bsnm{Kohlmann}, \binits{L.}},
\bauthor{\bsnm{Hartung}, \binits{A.}},
\bauthor{\bsnm{Gross}, \binits{J.}},
\bauthor{\bsnm{Scheel}, \binits{M.}},
\bauthor{\bsnm{Krack}, \binits{M.}}:
\batitle{Validation of a turbine blade component test with frictional contacts
  by phase-locked-loop and force-controlled measurements}.
\bjtitle{Journal of Engineering for Gas Turbines and Power}
\bvolume{142}(\bissue{5}),
\bfpage{1}--\blpage{8}
(\byear{2020})
\doiurl{10.1115/1.4044772}
\end{barticle}
\endbibitem

\bibitem[\protect\citeauthoryear{Tang et~al.}{2020}]{Tang2020}
\begin{botherref}
\oauthor{\bsnm{Tang}, \binits{M.}},
\oauthor{\bsnm{Stephan}, \binits{C.}},
\oauthor{\bsnm{B{\"{o}}swald}, \binits{M.}}:
Phase resonance method for nonlinear mechanical structures with phase locked
  loop control.
Proceedings of ISMA 2020 - International Conference on Noise and Vibration
  Engineering and USD 2020 - International Conference on Uncertainty in
  Structural Dynamics,
1805--1818
(2020)
\end{botherref}
\endbibitem

\bibitem[\protect\citeauthoryear{M{\"{u}}ller et~al.}{2022}]{Muller2022}
\begin{barticle}
\bauthor{\bsnm{M{\"{u}}ller}, \binits{F.}},
\bauthor{\bsnm{Woiwode}, \binits{L.}},
\bauthor{\bsnm{Gross}, \binits{J.}},
\bauthor{\bsnm{Scheel}, \binits{M.}},
\bauthor{\bsnm{Krack}, \binits{M.}}:
\batitle{Nonlinear damping quantification from phase-resonant tests under base
  excitation}.
\bjtitle{Mechanical Systems and Signal Processing}
\bvolume{177}(\bissue{January}),
\bfpage{109170}
(\byear{2022})
\doiurl{10.1016/j.ymssp.2022.109170}
\end{barticle}
\endbibitem

\bibitem[\protect\citeauthoryear{Abeloos et~al.}{2022}]{Abeloos2022}
\begin{barticle}
\bauthor{\bsnm{Abeloos}, \binits{G.}},
\bauthor{\bsnm{M{\"{u}}ller}, \binits{F.}},
\bauthor{\bsnm{Ferhatoglu}, \binits{E.}},
\bauthor{\bsnm{Scheel}, \binits{M.}},
\bauthor{\bsnm{Collette}, \binits{C.}},
\bauthor{\bsnm{Kerschen}, \binits{G.}},
\bauthor{\bsnm{Brake}, \binits{M.R.W.}},
\bauthor{\bsnm{Tiso}, \binits{P.}},
\bauthor{\bsnm{Renson}, \binits{L.}},
\bauthor{\bsnm{Krack}, \binits{M.}}:
\batitle{A consistency analysis of phase-locked-loop testing and control-based
  continuation for a geometrically nonlinear frictional system}.
\bjtitle{Mechanical Systems and Signal Processing}
\bvolume{170}(\bissue{October 2021}),
\bfpage{108820}
(\byear{2022})
\doiurl{10.1016/j.ymssp.2022.108820}
\end{barticle}
\endbibitem

\bibitem[\protect\citeauthoryear{Peter}{2018}]{Peter2018Thesis}
\begin{botherref}
\oauthor{\bsnm{Peter}, \binits{S.}}:
Nonlinear experimental modal analysis and its application to the identification
  of nonlinear structures.
PhD thesis,
University of Stuttgart
(2018).
\doiurl{10.2370/9783844061369} .
\url{https://www.researchgate.net/publication/329756057}
\end{botherref}
\endbibitem

\bibitem[\protect\citeauthoryear{Scheel}{2022}]{Scheel2022Thesis}
\begin{botherref}
\oauthor{\bsnm{Scheel}, \binits{M.}}:
Experimental nonlinear modal analysis - method development with particular
  focus on nonlinear damping.
PhD thesis,
University of Stuttgart
(2022).
\url{https://www.dr.hut-verlag.de/978-3-8439-5187-6.html}
\end{botherref}
\endbibitem

\bibitem[\protect\citeauthoryear{Link et~al.}{2011}]{Link2011}
\begin{botherref}
\oauthor{\bsnm{Link}, \binits{M.}},
\oauthor{\bsnm{Boeswald}, \binits{M.}},
\oauthor{\bsnm{Laborde}, \binits{S.}},
\oauthor{\bsnm{Weiland}, \binits{M.}},
\oauthor{\bsnm{Calvi}, \binits{A.}}:
Non-linear experimental modal analysis and application to satellite vibration
  test data.
ECCOMAS Thematic Conference - COMPDYN 2011: 3rd International Conference on
  Computational Methods in Structural Dynamics and Earthquake Engineering: An
  IACM Special Interest Conference, Programme
(May 2014)
(2011)
\end{botherref}
\endbibitem

\bibitem[\protect\citeauthoryear{Karaağa{\c{c}}lı and
  {\"{O}}zg{\"{u}}ven}{2021}]{Karaagacl2021}
\begin{barticle}
\bauthor{\bsnm{Karaağa{\c{c}}lı}, \binits{T.}},
\bauthor{\bsnm{{\"{O}}zg{\"{u}}ven}, \binits{H.N.}}:
\batitle{Experimental modal analysis of nonlinear systems by using
  response-controlled stepped-sine testing}.
\bjtitle{Mechanical Systems and Signal Processing}
\bvolume{146},
\bfpage{107023}
(\byear{2021})
\doiurl{10.1016/j.ymssp.2020.107023}
\end{barticle}
\endbibitem

\bibitem[\protect\citeauthoryear{B{\'{e}}liveau et~al.}{1986}]{Beliveau1986}
\begin{barticle}
\bauthor{\bsnm{B{\'{e}}liveau}, \binits{J.-G.}},
\bauthor{\bsnm{Vigneron}, \binits{F.R.}},
\bauthor{\bsnm{Soucy}, \binits{Y.}},
\bauthor{\bsnm{Draisey}, \binits{S.}}:
\batitle{Modal parameter estimation from base excitation}.
\bjtitle{Journal of Sound and Vibration}
\bvolume{107}(\bissue{3}),
\bfpage{435}--\blpage{449}
(\byear{1986})
\doiurl{10.1016/S0022-460X(86)80117-1}
\end{barticle}
\endbibitem

\bibitem[\protect\citeauthoryear{Carrella}{2020}]{Carrella2020}
\begin{botherref}
\oauthor{\bsnm{Carrella}, \binits{A.}}:
Introduction to Environmental Testing.
Siemens
(2020).
\url{https://community.sw.siemens.com/s/article/simcenter-testlab-vibration-control}
\end{botherref}
\endbibitem

\bibitem[\protect\citeauthoryear{Szempli{\'{n}}ska-Stupnicka}{1979}]{Szemplinska-Stupnicka1979}
\begin{barticle}
\bauthor{\bsnm{Szempli{\'{n}}ska-Stupnicka}, \binits{W.}}:
\batitle{The modified single mode method in the investigations of the resonant
  vibrations of non-linear systems}.
\bjtitle{Journal of Sound and Vibration}
\bvolume{63}(\bissue{4}),
\bfpage{475}--\blpage{489}
(\byear{1979})
\doiurl{10.1016/0022-460X(79)90823-X}
\end{barticle}
\endbibitem

\bibitem[\protect\citeauthoryear{Karaağa{\c{c}}lı and
  {\"{O}}zg{\"{u}}ven}{2022}]{Karaagacl2022}
\begin{barticle}
\bauthor{\bsnm{Karaağa{\c{c}}lı}, \binits{T.}},
\bauthor{\bsnm{{\"{O}}zg{\"{u}}ven}, \binits{H.N.}}:
\batitle{Experimental quantification and validation of modal properties of
  geometrically nonlinear structures by using response-controlled stepped-sine
  testing}.
\bjtitle{Experimental Mechanics}
\bvolume{62}(\bissue{2}),
\bfpage{199}--\blpage{211}
(\byear{2022})
\doiurl{10.1007/s11340-021-00784-9}
\end{barticle}
\endbibitem

\bibitem[\protect\citeauthoryear{G{\"{u}}rb{\"{u}}z et~al.}{2024}]{Gurbuz2024}
\begin{bchapter}
\bauthor{\bsnm{G{\"{u}}rb{\"{u}}z}, \binits{M.F.}},
\bauthor{\bsnm{Karaağa{\c{c}}lı}, \binits{T.}},
\bauthor{\bsnm{{\"{O}}zer}, \binits{M.B.}},
\bauthor{\bsnm{{\"{O}}zg{\"{u}}ven}, \binits{H.N.}}:
\bctitle{Bypassing the repeatability issue in nonlinear experimental modal
  analysis of jointed structures by using the {RCT-HFS} framework}.
In: \bbtitle{Society for Experimental Mechanics Annual Conference and
  Exposition},
pp. \bfpage{75}--\blpage{80}.
\bpublisher{Springer},
\blocation{Cham}
(\byear{2024}).
\doiurl{10.1007/978-3-031-36999-5_10} .
\burl{https://link.springer.com/10.1007/978-3-031-36999-5{\_}10}
\end{bchapter}
\endbibitem

\bibitem[\protect\citeauthoryear{Koyuncu et~al.}{2022}]{Koyuncu2022}
\begin{barticle}
\bauthor{\bsnm{Koyuncu}, \binits{A.}},
\bauthor{\bsnm{Karaağa{\c{c}}lı}, \binits{T.}},
\bauthor{\bsnm{Şahin}, \binits{M.}},
\bauthor{\bsnm{{\"{O}}zg{\"{u}}ven}, \binits{H.N.}}:
\batitle{Experimental modal analysis of nonlinear amplified piezoelectric
  actuators by using response-controlled stepped-sine testing}.
\bjtitle{Experimental Mechanics}
\bvolume{62}(\bissue{9}),
\bfpage{1579}--\blpage{1594}
(\byear{2022})
\doiurl{10.1007/s11340-022-00878-y}
\end{barticle}
\endbibitem

\bibitem[\protect\citeauthoryear{Karaağa{\c{c}}lı and {Nevzat
  {\"{O}}zg{\"{u}}ven}}{2024}]{Karaagacl2024}
\begin{bchapter}
\bauthor{\bsnm{Karaağa{\c{c}}lı}, \binits{T.}},
\bauthor{\bsnm{{Nevzat {\"{O}}zg{\"{u}}ven}}, \binits{H.}}:
\bctitle{Experimental modal analysis of structures with high nonlinear damping
  by using response-controlled stepped-sine testing}.
In: \beditor{\bsnm{Brake}, \binits{M.R.}},
\beditor{\bsnm{Renson}, \binits{L.}},
\beditor{\bsnm{Kuether}, \binits{R.J.}},
\beditor{\bsnm{Tiso}, \binits{P.}} (eds.)
\bbtitle{Nonlinear Structures {\&} Systems, Volume 1. SEM 2023. Conference
  Proceedings of the Society for Experimental Mechanics Series.},
pp. \bfpage{125}--\blpage{132}.
\bpublisher{Springer},
\blocation{Cham}
(\byear{2024}).
\bcomment{Chap. 19}.
\doiurl{10.1007/978-3-031-36999-5_17} .
\burl{https://link.springer.com/10.1007/978-3-031-36999-5{\_}17}
\end{bchapter}
\endbibitem

\bibitem[\protect\citeauthoryear{Karaaga{\c{c}}lı}{2020}]{Karaagacl2020thesis}
\begin{botherref}
\oauthor{\bsnm{Karaaga{\c{c}}lı}, \binits{T.}}:
Nonlinear system identification and nonlinear experimental modal analysis by
  using response controlled stepped sine testing.
PhD thesis,
Middle East Technical University
(2020).
\url{https://hdl.handle.net/11511/89603}
\end{botherref}
\endbibitem

\bibitem[\protect\citeauthoryear{Jones and Trefan}{2001}]{Jones2001}
\begin{barticle}
\bauthor{\bsnm{Jones}, \binits{B.K.}},
\bauthor{\bsnm{Trefan}, \binits{G.}}:
\batitle{The {D}uffing oscillator: A precise electronic analog chaos
  demonstrator for the undergraduate laboratory}.
\bjtitle{American Journal of Physics}
\bvolume{69}(\bissue{4}),
\bfpage{464}--\blpage{469}
(\byear{2001})
\doiurl{10.1119/1.1336838}
\end{barticle}
\endbibitem

\bibitem[\protect\citeauthoryear{Srinivasan et~al.}{2009}]{Srinivasan2009}
\begin{barticle}
\bauthor{\bsnm{Srinivasan}, \binits{K.}},
\bauthor{\bsnm{Thamilmaran}, \binits{K.}},
\bauthor{\bsnm{Venkatesan}, \binits{A.}}:
\batitle{Effect of nonsinusoidal periodic forces in {D}uffing oscillator:
  Numerical and analog simulation studies}.
\bjtitle{Chaos, Solitons {\&} Fractals}
\bvolume{40}(\bissue{1}),
\bfpage{319}--\blpage{330}
(\byear{2009})
\doiurl{10.1016/j.chaos.2007.07.090}
\end{barticle}
\endbibitem

\bibitem[\protect\citeauthoryear{Raze}{2024}]{Raze2024}
\begin{botherref}
\oauthor{\bsnm{Raze}, \binits{G.}}:
An electronic {D}uffing oscillator
(2024).
\url{https://github.com/GhislainRaze/Electronic-Duffing}
\end{botherref}
\endbibitem

\bibitem[\protect\citeauthoryear{No{\"{e}}l and Kerschen}{2013}]{Noel2013}
\begin{barticle}
\bauthor{\bsnm{No{\"{e}}l}, \binits{J.P.}},
\bauthor{\bsnm{Kerschen}, \binits{G.}}:
\batitle{Frequency-domain subspace identification for nonlinear mechanical
  systems}.
\bjtitle{Mechanical Systems and Signal Processing}
\bvolume{40}(\bissue{2}),
\bfpage{701}--\blpage{717}
(\byear{2013})
\doiurl{10.1016/j.ymssp.2013.06.034}
\end{barticle}
\endbibitem

\bibitem[\protect\citeauthoryear{Hippold et~al.}{2024}]{Hippold2024}
\begin{barticle}
\bauthor{\bsnm{Hippold}, \binits{P.}},
\bauthor{\bsnm{Scheel}, \binits{M.}},
\bauthor{\bsnm{Renson}, \binits{L.}},
\bauthor{\bsnm{Krack}, \binits{M.}}:
\batitle{Robust and fast backbone tracking via phase-locked loops}.
\bjtitle{Mechanical Systems and Signal Processing}
\bvolume{220}(\bissue{June}),
\bfpage{111670}
(\byear{2024})
\doiurl{10.1016/j.ymssp.2024.111670}
{\href{https://arxiv.org/abs/2403.06639}{{arXiv:2403.06639}}}
\end{barticle}
\endbibitem

\bibitem[\protect\citeauthoryear{Jain et~al.}{2019}]{Jain2019}
\begin{barticle}
\bauthor{\bsnm{Jain}, \binits{S.}},
\bauthor{\bsnm{Breunung}, \binits{T.}},
\bauthor{\bsnm{Haller}, \binits{G.}}:
\batitle{Fast computation of steady-state response for high-degree-of-freedom
  nonlinear systems}.
\bjtitle{Nonlinear Dynamics}
\bvolume{97}(\bissue{1}),
\bfpage{313}--\blpage{341}
(\byear{2019})
\doiurl{10.1007/s11071-019-04971-1}
\end{barticle}
\endbibitem

\bibitem[\protect\citeauthoryear{Robbins et~al.}{2023}]{Robbins2023}
\begin{barticle}
\bauthor{\bsnm{Robbins}, \binits{E.}},
\bauthor{\bsnm{Kuether}, \binits{R.J.}},
\bauthor{\bsnm{Pacini}, \binits{B.R.}},
\bauthor{\bsnm{Moreu}, \binits{F.}}:
\batitle{Stabilizing a strongly nonlinear structure through shaker dynamics in
  fixed frequency voltage control tests}.
\bjtitle{Mechanical Systems and Signal Processing}
\bvolume{190}(\bissue{December 2022}),
\bfpage{110118}
(\byear{2023})
\doiurl{10.1016/j.ymssp.2023.110118}
\end{barticle}
\endbibitem

\bibitem[\protect\citeauthoryear{Zhou and
  Kerschen}{2024}]{zhou2024identification}
\begin{botherref}
\oauthor{\bsnm{Zhou}, \binits{T.}},
\oauthor{\bsnm{Kerschen}, \binits{G.}}:
Identification of secondary resonances of nonlinear systems using phase-locked
  loop testing.
Journal of Sound and Vibration,
118549
(2024)
\doiurl{10.1016/j.jsv.2024.118549}
\end{botherref}
\endbibitem

\bibitem[\protect\citeauthoryear{Peeters et~al.}{2004}]{Peeters2004}
\begin{barticle}
\bauthor{\bsnm{Peeters}, \binits{B.}},
\bauthor{\bsnm{{Van der Auweraer}}, \binits{H.}},
\bauthor{\bsnm{Guillaume}, \binits{P.}},
\bauthor{\bsnm{Leuridan}, \binits{J.}}:
\batitle{The polymax frequency-domain method: A new standard for modal
  parameter estimation?}
\bjtitle{Shock and Vibration}
\bvolume{11}(\bissue{3-4}),
\bfpage{395}--\blpage{409}
(\byear{2004})
\doiurl{10.1155/2004/523692}
\end{barticle}
\endbibitem

\bibitem[\protect\citeauthoryear{Ewins}{2009}]{Ewins2009}
\begin{bbook}
\bauthor{\bsnm{Ewins}, \binits{D.J.}}:
\bbtitle{Modal Testing: Theory, Practice and Application},
(\byear{2009})
\end{bbook}
\endbibitem

\bibitem[\protect\citeauthoryear{Preumont}{2011}]{Preumont2011}
\begin{bbook}
\bauthor{\bsnm{Preumont}, \binits{A.}}:
\bbtitle{Vibration Control of Active Structures},
\bedition{3}rd edn.
\bsertitle{Solid Mechanics and Its Applications},
vol. \bseriesno{179}.
\bpublisher{Springer},
\blocation{Dordrecht}
(\byear{2011}).
\doiurl{10.1007/978-94-007-2033-6} .
\burl{http://link.springer.com/10.1007/978-94-007-2033-6}
\end{bbook}
\endbibitem

\bibitem[\protect\citeauthoryear{Quaegebeur et~al.}{2023}]{Quaegebeur2023}
\begin{barticle}
\bauthor{\bsnm{Quaegebeur}, \binits{S.}},
\bauthor{\bsnm{Raze}, \binits{G.}},
\bauthor{\bsnm{Cheng}, \binits{L.}},
\bauthor{\bsnm{Kerschen}, \binits{G.}}:
\batitle{A virtual acoustic black hole on a cantilever beam}.
\bjtitle{Journal of Sound and Vibration}
\bvolume{554}(\bissue{March}),
\bfpage{117697}
(\byear{2023})
\doiurl{10.1016/j.jsv.2023.117697}
{\href{https://arxiv.org/abs/2212.05939}{{arXiv:2212.05939}}}
\end{barticle}
\endbibitem

\bibitem[\protect\citeauthoryear{Goldberg et~al.}{1998}]{Goldberg1998}
\begin{botherref}
\oauthor{\bsnm{Goldberg}, \binits{P.W.}},
\oauthor{\bsnm{Williams}, \binits{C.K.I.}},
\oauthor{\bsnm{Bishop}, \binits{C.M.}}:
Regression with input-dependent noise a gaussian process treatment.
Advances in Neural Information Processing Systems,
493--499
(1998)
\end{botherref}
\endbibitem

\bibitem[\protect\citeauthoryear{Li and Dankowicz}{2021}]{Li2021}
\begin{barticle}
\bauthor{\bsnm{Li}, \binits{Y.}},
\bauthor{\bsnm{Dankowicz}, \binits{H.}}:
\batitle{Adaptive control designs for control-based continuation of periodic
  orbits in a class of uncertain linear systems}.
\bjtitle{Nonlinear Dynamics}
\bvolume{103}(\bissue{3}),
\bfpage{2563}--\blpage{2579}
(\byear{2021})
\doiurl{10.1007/s11071-021-06216-6}
\end{barticle}
\endbibitem

\bibitem[\protect\citeauthoryear{Li and Dankowicz}{2023}]{Li2023}
\begin{barticle}
\bauthor{\bsnm{Li}, \binits{Y.}},
\bauthor{\bsnm{Dankowicz}, \binits{H.}}:
\batitle{Model-free continuation of periodic orbits in certain nonlinear
  systems using continuous-time adaptive control}.
\bjtitle{Nonlinear Dynamics}
\bvolume{111}(\bissue{6}),
\bfpage{4945}--\blpage{4957}
(\byear{2023})
\doiurl{10.1007/s11071-022-08059-1}
{\href{https://arxiv.org/abs/2203.10306}{{arXiv:2203.10306}}}
\end{barticle}
\endbibitem

\bibitem[\protect\citeauthoryear{Rezaee and Renson}{2023}]{Rezaee2023}
\begin{botherref}
\oauthor{\bsnm{Rezaee}, \binits{H.}},
\oauthor{\bsnm{Renson}, \binits{L.}}:
Noninvasive adaptive control of a class of nonlinear systems with unknown
  parameters,
1--21
(2023)
{\href{https://arxiv.org/abs/2307.09806}{{arXiv:2307.09806}}}
\end{botherref}
\endbibitem

\bibitem[\protect\citeauthoryear{Kruse et~al.}{2024}]{Kruse2024}
\begin{botherref}
\oauthor{\bsnm{Kruse}, \binits{N.}},
\oauthor{\bsnm{Wallner}, \binits{H.}},
\oauthor{\bsnm{Dittus}, \binits{A.}},
\oauthor{\bsnm{Lukas}, \binits{B.}}:
{Large basins of attraction for control-based continuation of unstable periodic
  states},
1--24
(2024)
\end{botherref}
\endbibitem

\bibitem[\protect\citeauthoryear{Lim and Epureanu}{2011}]{Lim2011}
\begin{barticle}
\bauthor{\bsnm{Lim}, \binits{J.}},
\bauthor{\bsnm{Epureanu}, \binits{B.I.}}:
\batitle{{Forecasting a class of bifurcations: Theory and experiment}}.
\bjtitle{Physical Review E}
\bvolume{83}(\bissue{1}),
\bfpage{016203}
(\byear{2011})
\doiurl{10.1103/PhysRevE.83.016203}
\end{barticle}
\endbibitem

\bibitem[\protect\citeauthoryear{Habib}{2023}]{Habib2023}
\begin{barticle}
\bauthor{\bsnm{Habib}, \binits{G.}}:
\batitle{{Predicting saddle-node bifurcations using transient dynamics: a
  model-free approach}}.
\bjtitle{Nonlinear Dynamics}
\bvolume{111}(\bissue{22}),
\bfpage{20579}--\blpage{20596}
(\byear{2023})
\doiurl{10.1007/s11071-023-08941-6}
\end{barticle}
\endbibitem

\bibitem[\protect\citeauthoryear{Shen et~al.}{2021}]{Shen2021b}
\begin{barticle}
\bauthor{\bsnm{Shen}, \binits{J.}},
\bauthor{\bsnm{Groh}, \binits{R.M.J.}},
\bauthor{\bsnm{Schenk}, \binits{M.}},
\bauthor{\bsnm{Pirrera}, \binits{A.}}:
\batitle{Experimental path-following of equilibria using newton's method. part
  ii: Applications and outlook}.
\bjtitle{International Journal of Solids and Structures}
\bvolume{213},
\bfpage{25}--\blpage{40}
(\byear{2021})
\doiurl{10.1016/j.ijsolstr.2020.11.038}
\end{barticle}
\endbibitem

\bibitem[\protect\citeauthoryear{Magnevall et~al.}{2006}]{Magnevall2006}
\begin{botherref}
\oauthor{\bsnm{Magnevall}, \binits{M.}},
\oauthor{\bsnm{Josefsson}, \binits{A.}},
\oauthor{\bsnm{Ahlin}, \binits{K.}}:
Experimental verification of a control algorithm for nonlinear systems.
Conference Proceedings of the Society for Experimental Mechanics Series
(2006)
\end{botherref}
\endbibitem

\bibitem[\protect\citeauthoryear{Novak et~al.}{2018}]{Novak2018}
\begin{barticle}
\bauthor{\bsnm{Novak}, \binits{A.}},
\bauthor{\bsnm{Simon}, \binits{L.}},
\bauthor{\bsnm{Lotton}, \binits{P.}}:
\batitle{A simple predistortion technique for suppression of nonlinear effects
  in periodic signals generated by nonlinear transducers}.
\bjtitle{Journal of Sound and Vibration}
\bvolume{420},
\bfpage{104}--\blpage{113}
(\byear{2018})
\doiurl{10.1016/j.jsv.2018.01.038}
\end{barticle}
\endbibitem

\bibitem[\protect\citeauthoryear{Raze et~al.}{2024}]{Raze2024fig}
\begin{botherref}
\oauthor{\bsnm{Raze}, \binits{G.}},
\oauthor{\bsnm{Abeloos}, \binits{G.}},
\oauthor{\bsnm{Kerschen}, \binits{G.}}:
Experimental data from "experimental continuation in nonlinear dynamics: recent
  advances and future challenges"
(2024)
\doiurl{10.6084/m9.figshare.26412187}
\end{botherref}
\endbibitem

\end{thebibliography}

\section*{Statements and Declarations}


\begin{itemize}
\item \textbf{Funding:} Gaëtan Abeloos and Ghislain Raze are a FRIA Grantee and a Postdoctoral Researcher, respectively, of the Fonds de la Recherche
Scientifique - FNRS, which is gratefully acknowledged.
\item \textbf{Competing interests:} The authors have no relevant financial or non-financial interests to disclose.
\item \textbf{Author Contributions:} 
Ghislain Raze: literature search, formal analysis, data analysis, investigation and writing - original draft preparation. Gaëtan Abeloos: conceptualization, methodology and writing - review and editing. Gaëtan Kerschen: literature search, funding acquisition, supervision and writing - review and editing.
\item \textbf{Data Availability:}  The data collected and used
in this research is available in~\cite{Raze2024fig}.
\end{itemize}

\end{document}